\documentclass[11pt,twoside]{article}
\usepackage{latexsym}
\usepackage{rotating}
\input{amssym.def}
\input{amssym}
\newfont{\fra}{eufm10 scaled 1095}
\newfont{\Bb}{msbm10 scaled 1095}
\newfont{\Bbg}{msbm10 scaled 1280}
\newcommand\CC{{\mbox{\Bb C}}}
\newcommand\RR{{\mbox{\Bb R}}}
\newcommand\NN{{\mbox{\Bb N}}}
\newcommand\ZZ{{\mbox{\Bb Z}}}

\newcommand\fg{{\frak{g}}}
\newcommand\fh{{\frak h}}
\newcommand\fri{{\frak i}}
\newcommand\fj{{\frak j}}
\newcommand\fl{{\frak l}}
\newcommand\fm{{\frak m}}
\newcommand\fn{{\frak n}}
\newcommand\fa{{\frak a}}
\newcommand\fb{{\frak b}}
\newcommand\fd{{\frak d}}

\newcommand\fk{{\frak k}}
\newcommand\fr{{\frak r}}
\newcommand\fs{{\frak s}}
\newcommand\ft{{\frak t}}
\newcommand\fz{{\frak z}}
\newcommand\fS{{\frak S}}
\newcommand\cC{{\cal C}}
\newcommand\cZ{{\cal Z}}
\newcommand\cH{{\cal H}}
\newcommand\cM{{\cal M}}
\newcommand\cA{{\cal A}}
\newcommand\cB{{\cal B}}
\newcommand{\gl}{\mathop{{\frak g \frak l}}}
\newcommand{\fsl}{\mathop{{\frak s \frak l}}}
\newcommand{\fsu}{\mathop{{\frak s \frak u}}}

\newcommand{\so}{\mathop{{\frak s \frak o}}}

\newcommand{\Aut}{\mathop{{\rm Aut}}}
\newcommand{\GL}{\mathop{{\it GL}}}
\newcommand{\SL}{\mathop{{\it SL}}}
\newcommand{\Hom}{\mathop{{\rm Hom}}}

\newcommand{\Id}{\mathop{{\rm Id}}}
\newcommand{\id}{\mathop{{\rm Id}}}
\newcommand{\ad}{\mathop{{\rm ad}}}
\newcommand{\tr}{\mathop{{\rm tr}}}
\newcommand{\Ad}{\mathop{{\rm Ad}}}
\newcommand{\adR}{\mathop{{\rm ad}_{0}}}

\newcommand{\sgn}{\mathop{{\rm sgn}}}

\newcommand{\Ker}{\mathop{{\rm ker}}}
\newcommand{\im}{\mathop{{\rm im}}}

\renewcommand{\Re}{\mathop{{\rm Re}}}
\renewcommand{\Im}{\mathop{{\rm Im}}}

\newcommand{\Span}{\mathop{{\rm span}}}
\newcommand{\mod}{\mathop{{\rm mod}}}
\newcommand{\proj}{{\rm proj}}
\newcommand{\cyclsum}{\mathop{\sum_{\mbox{\footnotesize{\rm cycl}}}}}
\newcommand\ip{\mbox{$\langle\cdot \,,\cdot \rangle$}}
\newcommand\ipa{{\langle\cdot \,,\cdot \rangle_\fa}}
\newcommand\lb{{[\cdot\,,\cdot]}}

\newcommand\dd{\fd_{\alpha,\gamma}(\fl,\fa,\rho)}
\newcommand\proof{{\sl Proof. }}
\newcommand{\qed}{\hspace*{\fill}\hbox{$\Box$}\vspace{2ex}}
\newtheorem{theo}{Theorem}[section]
\newtheorem{pr}{Proposition}[section]
\newtheorem{de}{Definition}[section]

\newtheorem{re}{Remark}[section]
\newtheorem{co}{Corollary}[section]
\newtheorem{lm}{Lemma}[section]
\begin{document}
\title{Metric Lie algebras and quadratic extensions}
\author{Ines Kath and Martin Olbrich}
\maketitle
\begin{abstract}
\noindent
The present
paper contains a systematic study of the structure of metric Lie algebras, i.e.,
finite-dimensional real Lie algebras equipped with a non-degenerate invariant symmetric bilinear form. We show that any metric Lie algebra
$\fg$ without simple ideals has the structure of a so called balanced
quadratic extension of an auxiliary Lie algebra $\fl$ by an orthogonal
$\fl$-module $\fa$ in a canonical way. Identifying equivalence classes
of quadratic extensions of $\fl$ by $\fa$ with a certain cohomology set $\cH^2_Q(\fl,\fa)$
we obtain a classification scheme for general
metric Lie algebras
and a complete classification of metric
Lie algebras of index $3$.
\end{abstract}
\tableofcontents
%%%%%%%%%%%%%%%%%%%%%%%%%%%%%%%%%%%%%%%%%%%%%%%%%%%%%%%%%%%%%%%%%%%%%%%%%%%%%%%%

\section{Introduction}

The present paper is an attempt towards the classification
of metric Lie algebras up to isomorphism. Here a metric Lie algebra is a finite-dimensional real Lie algebra equipped with an
invariant non-degenerate symmetric bilinear form.
An isomorphism
of metric Lie algebras is by definition a Lie algebra isomorphism
which is in addition an isometry with respect
to the given inner products.

A metric Lie algebra $(\fg,\ip)$ is called decomposable if
it contains a proper ideal $\fk$ which is non-degenerate (i.e. $\ip_{|\fk\times\fk}$ is
non-degenerate), and indecomposable otherwise.
Then any metric Lie algebra
is the orthogonal direct sum of uniquely determined indecomposable metric Lie algebras. Note that indecomposable metric Lie algebras are
either simple or do not contain simple ideals.
Therefore, one is lead to investigate
indecomposable metric Lie algebras without simple ideals.

Our study is mainly motivated by the wish of understanding general
pseudo-Riemannian symmetric spaces. Indeed, the Lie algebra of the
transvection group of a pseudo-Riemannian symmetric space has the
structure of a metric Lie algebra (\cite{CP80}, Prop. 1.6).
Moreover,
the local geometry of a pseudo-Riemannian symmetric space is completely
determined by this metric Lie algebra together with the isometric involutive automorphism of it induced by the geodesic symmetry.
Our theory of metric Lie algebras is designed in such a way that
the incorporation of an involutive automorphism into the structure does
not present serious additional difficulties. Thus there is a theory
of pseudo-Riemannian symmetric spaces completely parallel to the theory
of metric Lie algebras developed in the present paper. The details will
appear in a forthcoming paper.

For further motivation and remarks
on the history of the subject we refer to the introduction of \cite{KO03}.
In the paper \cite{KO03} we studied the relatively special class of
metric Lie
algebras $(\fg,\ip)$ which satisfy $\fg''\subset \fz(\fg)$. One consequence
of this investigation was the complete classification of indecomposable metric Lie algebras
of index 2  (this classification has been already announced in \cite{BK02}, the classification of
metric Lie algebras of index 1 is due to Medina \cite{M85}).

In the present paper we are able to carry over the approach of \cite{KO03}
to general metric Lie algebras. Let us outline the main ideas and results
of the paper.

The basic construction, which goes back to an idea
of Berard Bergery used in his unpublished work on pseudo-Riemannian
symmetric spaces \cite{BB2}, is the following:
Let $(\fg,\ip)$ be a metric Lie algebra. Assume that we have an isotropic
ideal $\fri\subset \fg$ such that $\fri^\perp/\fri$ is abelian. Here
$\fri^\perp$ denotes the orthogonal ``complement'' of $\fri$. Set
$\fl=\fg/\fri^\perp$ and $\fa=\fri^\perp/\fri$. Then $\fa$ inherits
an inner product from $\fg$ and an $\fl$-action respecting this inner product, i.e., it inherits the structure of an orthogonal $\fl$-module. Moreover, $\fri\cong \fl^*$ as an $\fl$-module, and $\fg$ can be represented
as the result of two subsequent extensions of Lie algebras with abelian kernel
\begin{equation}
0\rightarrow \fa{\longrightarrow}\fg/\fri
{\longrightarrow} \fl\rightarrow 0\ ,\qquad
0\rightarrow \fl^*{\longrightarrow}\fg
{\longrightarrow} \fg/\fri\rightarrow 0\ .\label{lena}
\end{equation}
Vice versa, given a Lie algebra $\fl$, an orthogonal $\fl$-module
$\fa$ and two extensions as in (\ref{lena}) which in addition satisfy
certain natural compatibility conditions, the resulting Lie algebra $\fg$
has a distinguished invariant inner product. This construction of metric Lie algebras,
being a  relative of the double extension method of Medina and Revoy \cite{MR85},
will be formalised into the notion of a quadratic extension of $\fl$ by
an orthogonal $\fl$-module $\fa$ in Subsection \ref{leonce}. In particular, there is a natural
equivalence relation on the set of quadratic extensions of $\fl$ by $\fa$.
The cocycles defining the extensions in
(\ref{lena}) represent an element in a certain cohomology set $\cH^2_Q(\fl,\fa)$, and it turns out that there is a bijection
between equivalence classes of such quadratic extensions and $\cH^2_Q(\fl,\fa)$ (Theorem \ref{brandt}). We introduce these cohomology
sets and study their basic functorial properties in Section \ref{qc}.
Cohomology sets of this kind were studied first by Grishkov \cite{grishkov}.

What makes the theory of quadratic extensions so useful is the fact
that any metric Lie algebra $(\fg,\ip)$ without simple ideals has a canonical
isotropic ideal $\fri=\fri(\fg)$ (Definitions \ref{meer} and \ref{schwein}) such that $\fri(\fg)^\perp/\fri(\fg)$
is abelian. In other words, $(\fg,\ip)$ has a canonical structure
of a quadratic extension of $\fl=\fg/\fri(\fg)^\perp$ by $\fa=\fri(\fg)^\perp/\fri(\fg)$ (Proposition \ref{schmidt}).
It has the
property that the orthogonal $\fl$-module $\fa$
is semi-simple. However, not every
quadratic extension of a Lie algebra $\fl$ by a semi-simple orthogonal $\fl$-module $\fa$ arises in this way. The obvious condition that
the image of $\fl^*$ in $\fg$ in (\ref{lena}) should be equal to $\fri(\fg)$ is
not always satisfied. If it is satisfied we call the quadratic extension
balanced and the corresponding cohomology class in $\cH^2_Q(\fl,\fa)$
admissible. The main result of Section \ref{baal} is Theorem \ref{bebe}
which characterises the admissible cohomology classes in $\cH^2_Q(\fl,\fa)$.

It is now easy to decide which admissible cohomology classes correspond
to decomposable metric Lie algebras (see
Section \ref{isodec}, in particular Definition \ref{kosel} for the definition of indecomposable cohomology classes). We denote the set of indecomposable admissible cohomology classes by $\cH^2_Q(\fl,\fa)_0$.
Moreover, it turns out that elements of $\cH^2_Q(\fl,\fa)_0$ correspond
to isomorphic Lie algebras if and only if they can be transformed into each other by the induced action of the automorphisms group $G_{\fl,\fa}$ of
the pair $(\fl,\fa)$.
Thus we can identify the set of isomorphism classes of non-simple indecomposable metric Lie algebras with the union of orbit spaces
\begin{equation}\label{class}
 \coprod_{\fl,\fa} \cH^2_Q(\fl,\fa)_0/G_{\fl,\fa}\ ,
\end{equation}
where the union is taken over a set of representatives of isomorphism
classes of pairs $(\fl,\fa)$ consisting of real finite-dimensional Lie algebras $\fl$ and a semi-simple orthogonal $\fl$-module $\fa$ (Theorem \ref{Pwieder}). This is the classification scheme we aimed at.
In particular, it says that an arbitrary non-simple indecomposable
Lie algebra can be constructed as the quadratic extension corresponding
to an element of $\cH^2_Q(\fl,\fa)_0$ for some pair $(\fl,\fa)$ and it
clearly indicates when two metric Lie algebras constructed in this way
are isomorphic.

In order to approach a true classification one has to evaluate (\ref{class})
further. One  first observes (see Section \ref{SM}) that the class
of Lie algebras $\fl$ which really occurs in (\ref{class}), i.e.,
$\cH^2_Q(\fl,\fa)_0\ne\emptyset$ for some orthogonal $\fl$-module $\fa$,
is a proper subclass of all Lie algebras. The critical condition is that
there has to exist an admissible cohomology class for some orthogonal
$\fl$-module $\fa$ (in contrast, the indecomposability condition is harmless
and gives no restrictions on the Lie algebra $\fl$).
It is an open question how large this class of Lie algebras
really is; see the comment at the end of Section \ref{isodec}.
Some partial results in this direction are obtained in Section \ref{SM}.

By construction, the dimension of the Lie algebra $\fl$ associated
with a metric Lie algebra $(\fg,\ip)$ is bounded by the index of $\ip$.
For this reason the bijection (\ref{class}) is extremely useful for a concrete classification of metric Lie algebras with
small index, where all the ingredients of (\ref{class}) can be explicitely computed. This
is demonstrated in Section \ref{i3} for index $3$.

\noindent
{\it Acknowledgements}:
Parts of the results were obtained during a stay of both authors at the
Institut \' Elie Cartan in Nancy.
We gratefully acknowledge the hospitality of the institute
as well as the lively working atmosphere there. We are particularly indebted
to L. Berard Bergery for sharing his insight into the subject with us.

\section{Quadratic cohomology}\label{qc}

Let $\fl$ be a finite-dimensional real Lie algebra, and let $\rho: 
\fl\rightarrow \gl(\fa)$
be a representation of $\fl$ on a finite-dimensional real vector
space $\fa$ equipped with an invariant symmetric bilinear form
$\ip$, i.e. for $L\in\fl$ and
$v,w\in\fa$
$$ \langle \rho(L)v,w \rangle +\langle v,\rho(L)w\rangle = 0\ .$$
In addition, we require $\ip$ to be non-degenerate.
Then the triple $(\rho,\fa,\ip)$ is called an orthogonal $\fl$-module.
Often $\rho$ or $\ip$ will be omitted in the notation.

The goal of this
section is to associate with an orthogonal $\fl$-module a sequence of sets
$$ \cH^p_Q(\fl,\fa)\ ,\quad p\in 2\NN_0\ ,$$
the quadratic cohomology sets of $\fl$ with coefficients in
$(\rho,\fa,\ip)$. Though
they
are specializations of the cohomology sets associated by Grishkov 
\cite{grishkov} to a cochain complex with a cup product taking values
in a second cochain complex we prefer to present a self-contained 
treatment
here.

Let us recall the construction of usual Lie algebra cohomology. For 
any representation $\rho:\fl\rightarrow \gl(\fa)$ of a Lie algebra 
$\fl$
on a vector space $\fa$ we have the standard Lie algebra cochain 
complex
$(C^*(\fl,\fa),d)$, where $C^p(\fl,\fa)=\Hom(\bigwedge^p\fl,\fa)$, 
and     
for $\tau\in C^p(\fl,\fa)$
\begin{eqnarray*}
d\tau(L_1,\dots,L_{p+1})&=&\sum_{i=1}^{p+1} (-1)^{i-1} \rho(L_i) 
\tau(L_1,\dots,\hat L_i,\dots,L_{p+1}) \\
&& +\sum_{i<j} (-1)^{i+j}\tau([L_i,L_j],L_1,\dots,\hat L_i,\dots,\hat 
L_j,\dots,L_{p+1})\ .
\end{eqnarray*}
We denote the groups of cocycles and coboundaries of $C^p(\fl,\fa)$ by
$Z^p(\fl,\fa)$ and $B^p(\fl,\fa)$, respectively.
Then the quotients
$$ H^p(\fl,\fa):= Z^p(\fl,\fa)/B^p(\fl,\fa)$$ 
constitute
the cohomology groups of $\fl$ with coefficients in $\fa$.

Now let $(\rho,\fa,\ip)$ be an orthogonal $\fl$-module.
By $C^*(\fl)$ we denote the cochain complex
associated to the trivial one-dimensional representation.
Then the composition of maps
$$ C^p(\fl,\fa)\times C^q(\fl,\fa)\stackrel{\wedge}{\longrightarrow}
C^{p+q}(\fl,\fa\otimes\fa)\stackrel{\ip}{\longrightarrow} 
C^{p+q}(\fl)$$
defines a bilinear multiplication
$$ C^p(\fl,\fa)\times C^q(\fl,\fa) \rightarrow C^{p+q}(\fl)$$
which we will denote by $(\alpha,\tau)\mapsto \langle 
\alpha\wedge\tau\rangle$.
In concrete terms
\begin{eqnarray*}\lefteqn{\langle 
\alpha\wedge\tau\rangle(L_1,\dots,L_{p+q})}\\
&=& \sum_{[\sigma]\in\fS_{p+q}/{\fS_p\times\fS_q}}
\sgn(\sigma)\langle \alpha(L_{\sigma(1)},
\dots,L_{\sigma(p)}),\tau(L_{\sigma(p+1)},
\dots,L_{\sigma(p+q)})\rangle\ ,
\end{eqnarray*}
where $\fS_k$ denotes the symmetric group in $k$ letters. As a consequence
of the invariance of $\ip$ we have for $\alpha\in C^p(\fl,\fa)$ and $\tau \in
C^q(\fl,\fa)$
\begin{equation}\label{vita}
d\langle \alpha\wedge\tau\rangle=\langle d\alpha\wedge\tau\rangle
+(-1)^p \langle \alpha\wedge d\tau\rangle\ .
\end{equation}
It follows that $\langle\cdot\wedge\cdot\rangle$
induces a kind of cup product on the cohomology groups
\begin{equation}\label{stokes}
\cup: H^p(\fl,\fa)\times H^q(\fl,\fa) \longrightarrow H^{p+q}(\fl)\
,\quad [\alpha]\cup [\tau]:=[\langle \alpha\wedge\tau\rangle ]\ .
\end{equation}
Set $n=\dim\fl$. If $\fl$ is unimodular which is equivalent
to $H^n(\fl)\ne\{0\}$, then we obtain a non-degenerate pairing
(Poincar\'e duality)
\begin{equation}\label{poincare}
\cup: H^p(\fl,\fa)\times H^{n-p}(\fl,\fa) \longrightarrow H^{n}(\fl)\cong \RR .
\end{equation}

\begin{de}\label{unmoti}
Let $p$ be even. We define a group structure on 
$\cC^{p-1}_Q(\fl,\fa):=C^{p-1}(\fl,\fa)\oplus
C^{2p-2}(\fl)$ by
$$ (\tau_1,\sigma_1)*(\tau_2,\sigma_2):=
(\tau_1+\tau_2,\sigma_1+\sigma_2+
\textstyle\frac{1}{2}\langle\tau_1\wedge\tau_2\rangle)\ .$$
We call $\cC^{p-1}_Q(\fl,\fa)$ the group of quadratic 
$(p-1)$-cochains.
The set of quadratic $p$-cocycles is given by 
$$\cZ^p_Q(\fl,\fa):= \{(\alpha,\gamma)\in C^p(\fl,\fa)\oplus
C^{2p-1}(\fl)\:|\: d\alpha=0,\, d\gamma= 
\textstyle\frac{1}{2}\langle\alpha\wedge\alpha\rangle\}\ .$$
\end{de}
Clearly, $*$ is associative. The inverse of $(\tau,\sigma)$ is given 
by
$(-\tau, \textstyle\frac{1}{2}\langle\tau\wedge\tau\rangle -\sigma)$.
Thus $*$ indeed defines a group structure on $\cC^{p-1}_Q(\fl,\fa)$.

\begin{lm}\label{pepsi}
Let $p$ be even. Let $(\alpha,\gamma)\in C^p(\fl,\fa)\oplus 
C^{2p-1}(\fl)$ 
and $(\tau,\sigma)\in C^{p-1}(\fl,\fa)\oplus
C^{2p-2}(\fl)$. Then the formula
\begin{equation}\label{unmotiviert}
(\alpha,\gamma)(\tau,\sigma):=\left(\alpha+d\tau,\gamma+d\sigma+\langle
(\alpha+\textstyle\frac{1}{2}d\tau)\wedge\tau \rangle \right)
\end{equation}
defines a right action of the
group $\cC^{p-1}_Q(\fl,\fa)$ on $C^p(\fl,\fa)\oplus
C^{2p-1}(\fl)$.
This action leaves the set of $p$-cocycles $\cZ^p_Q(\fl,\fa)\subset 
C^p(\fl,\fa)\oplus
C^{2p-1}(\fl)$ invariant.
\end{lm}
\proof
Let $(\tau_i,\sigma_i)\in C^{p-1}(\fl,\fa)\oplus
C^{2p-2}(\fl)$, $i=1,2$.
Since $p$ is even we have by (\ref{vita}) that
$$ d\langle\tau_1\wedge\tau_2\rangle =\langle d\tau_1\wedge 
\tau_2\rangle -\langle d\tau_2\wedge\tau_1\rangle \ .$$
We obtain
$$\textstyle\frac{1}{2}\left(\langle d\langle 
\tau_1\wedge\tau_2\rangle +\langle d\tau_1\wedge \tau_2\rangle 
+\langle d\tau_2\wedge\tau_1\rangle \right)=\langle 
d\tau_1\wedge\tau_2\rangle \ .$$
Therefore we have for $(\alpha,\gamma)\in C^p(\fl,\fa)\oplus 
C^{2p-1}(\fl)$ 
\begin{eqnarray*}
\lefteqn{(\alpha,\gamma)((\tau_1,\sigma_1)*(\tau_2,\sigma_2))}\\
\hspace{-3pt}&=&\hspace{-3pt}(\alpha,\gamma)+\left( 
d(\tau_1+\tau_2),
d(\sigma_1+\sigma_2+\textstyle\frac{1}{2}\langle\tau_1\wedge\tau_2\rangle)
+\langle 
(\alpha+\textstyle\frac{1}{2}d(\tau_1+\tau_2))\wedge(\tau_1+\tau_2) 
\rangle\right)\\
\hspace{-3pt}
&=&\hspace{-4pt}\left(\alpha+d\tau_1,\gamma+d\sigma_1+\langle 
(\alpha+\textstyle\frac{1}{2}d\tau_1)\wedge\tau_1\rangle\right)
\hspace{-2.4pt}+\hspace{-2.4pt}
\left(d\tau_2,d\sigma_2+\langle 
(\alpha+d\tau_1+\textstyle\frac{1}{2}d\tau_2)\wedge\tau_2\rangle\right)\\
\hspace{-3pt}
&=&\hspace{-3pt} ((\alpha,\gamma)(\tau_1,\sigma_1))(\tau_2,\sigma_2)\ 
.
\end{eqnarray*}
This proves the first assertion of the lemma.

Now let $(\alpha,\gamma)\in \cZ^p_Q(\fl,\fa)$, i.e.
$ d\alpha=0$, $ d\gamma= \frac{1}{2}\langle\alpha\wedge\alpha\rangle$.
Using again Equation (\ref{vita}) we find for $(\tau,\sigma)\in
\cC^{p-1}_Q(\fl,\fa)$ \begin{eqnarray*}
d(\gamma+d\sigma+\langle 
(\alpha+\textstyle\frac{1}{2}d\tau)\wedge\tau \rangle)
&=& \textstyle\frac{1}{2}\langle\alpha\wedge\alpha\rangle
+\langle 
(\alpha+\textstyle\frac{1}{2}d\tau)\wedge d\tau \rangle\\
&=&\textstyle\frac{1}{2}\langle(\alpha+d\tau)\wedge(\alpha+d\tau)\rangle\ 
.
\end{eqnarray*}
This implies $(\alpha,\gamma)(\tau,\sigma)\in \cZ^p_Q(\fl,\fa)$.
The proof of the lemma is now complete.
\qed

\begin{re} {\rm The subgroup $Z^{2p-2}(\fl)\subset 
\cC^{p-1}(\fl,\fa)$ acts
trivially on $C^p(\fl,\fa)\oplus
C^{2p-1}(\fl)$. Thus the above action amounts to an action of the 
group
$$\cB^{p}_Q(\fl,\fa):=(C^{p-1}(\fl,\fa)\oplus B^{2p-1}(\fl), *)\ ,$$
where the multiplication $*$ is given by}
$$ (\tau_1,\sigma_1)*(\tau_2,\sigma_2):=
(\tau_1+\tau_2,\sigma_1+\sigma_2+
\textstyle\frac{1}{2}d\langle\tau_1\wedge\tau_2\rangle)\ .$$
\end{re} 

\begin{de}
For any $p\in 2\NN_0$ we define the $p$-th cohomology set 
$\cH^p_Q(\fl,\fa)$  
of $\fl$ with coefficients in $(\fa,\ip)$ as the quotient space
of the action of $\cC^{p-1}_Q(\fl,\fa)$ (or $\cB^p_Q(\fl,\fa)$) on 
$\cZ^p_Q(\fl,\fa)$.
If $(\alpha,\gamma)\in \cZ^p_Q(\fl,\fa)$, then we denote the 
corresponding
cohomology class in $\cH^p_Q(\fl,\fa)$ by $[\alpha,\gamma]$. 
\end{de}

The sets $\cH^p_Q(\fl,\fa)$ appear among the cohomology sets
introduced and
studied by Grish\-kov in \cite{grishkov}. There they are denoted by 
$H^{2p-1}_\Delta$,
where $\Delta$ is the multiplication given by 
$\frac{1}{2}\langle\cdot\wedge\cdot\rangle$. Grishkov also defines 
corresponding
cohomology sets for odd $p$ in an analogous way. 
We have avoided to discuss them here because
they do not give anything new. In fact, they turn out to be
canonically isomorphic to
$H^p(\fl,\fa)\oplus H^{2p-1}(\fl)$.

Note that $\cH^0_Q(\fl,\fa)$ is equal to the set of isotropic 
invariant
vectors in $\fa$. In the present paper the set 
$\cH^2_Q(\fl,\fa)$ will be the main technical tool for studying 
metric Lie algebras. For abelian $\fl$ this set has been already 
intensively
investigated in \cite{KO03}.

We will need the basic functorial properties of the cohomology sets 
$\cH^{p}_Q(\fl,\fa)$,
$p\in 2\NN_0$. Similar as usual Lie algebra cohomology groups they turn out to be 
contravariant
with respect to Lie algebras and covariant with respect to modules. 

We consider the category $\cal{LO}$ of pairs $(\fl,\fa)$ of Lie algebras $\fl$
and orthogonal $\fl$-modules. Morphisms $(\fl_1,\fa_1)\mapsto (\fl_2,\fa_2)$
are given by pairs $(S,U)$, where $S:\fl_1\rightarrow \fl_2$ is a Lie algebra
homomorphism and $U:\fa_2\rightarrow \fa_1$ is an isometric embedding such that
\begin{equation}\label{opair}
U\circ\rho_2(S(L))=\rho_1(L)\circ U\ .
\end{equation}
Notice that $U$ maps in the reverse direction.
The condition (\ref{opair}) means that $U\in\Hom_{\fl_1}(\fa_2,\fa_1)$, where
the $\fl_1$-module structure of $\fa_2$ is given by $S^*\rho_2:=\rho_2\circ S$.
We call morphisms in $\cal{LO}$ morphisms of pairs.

Let $(S,U): (\fl_1,\fa_1)\rightarrow (\fl_2,\fa_2)$ be a morphism of pairs.
For all $p\in\NN_0$ we then have the
pull back maps
$$ S^*: C^p(\fl_2)\longrightarrow C^p(\fl_1) \ ,\quad 
S^*\gamma(L_1,\dots,L_p):=\gamma(S(L_1),\dots,S(L_p)) $$
and
$$(S,U)^*:
 C^p(\fl_2,\fa_2)\longrightarrow C^p(\fl_1,\fa_1)\ ,\quad 
(S,U)^*\alpha(L_1,\dots,L_p):=U\circ\alpha(S(L_1),\dots,S(L_p))\ . $$
$S^*$ and $(S,U)^*$ commute with the differentials. Moreover we have
$$\langle (S,U)^*\alpha\wedge (S,U)^*\tau \rangle = S^*\langle \alpha\wedge \tau\rangle \
.$$
For this property it is crucial that $U$ is an isometric embedding.
Thus the direct sum $(S,U)^*\oplus S^*$ maps the space of
cocycles $\cZ^p_Q(\fl_2,\fa_2)\subset C^p(\fl_2,\fa_2)\oplus C^{2p-1}(\fl_2)$
to $\cZ^p_Q(\fl_1,\fa_1)$ and defines a group homomorphism from
$\cC^{p-1}_Q(\fl_2,\fa_2)$ to $\cC^{p-1}_Q(\fl_1,\fa_1)$. Moreover, $(S,U)^*\oplus S^*$
intertwines the $\cC^{p-1}_Q$-actions defined in Lemma \ref{pepsi}:
$$ ((S,U)^*\oplus S^*) ((\alpha,\gamma)(\tau,\sigma))=
((S,U)^*\alpha,S^*\gamma)((S,U)^*\tau,S^*\sigma)\ .
$$
We obtain

\begin{pr}\label{pull}
Let $F=(S,U):(\fl_1,\fa_1)\rightarrow (\fl_2,\fa_2)$ be a morphism of pairs.
Then for $p\in 2\NN_0$ there is a pull back map
$$ F^*: \cH^p_Q(\fl_2,\fa_2)\longrightarrow \cH^p_Q(\fl_1,\fa_1)$$
given by
$$ F^*[\alpha,\gamma]:= [(S,U)^*\alpha, S^*\gamma]\ ,\quad (\alpha,\gamma)\in
\cZ^p_Q(\fl_2,\fa_2)\ .$$
Therefore the assignment $(\fl,\fa)\leadsto\cH^p_Q(\fl,\fa)$ is a contravariant
functor from the category $\cal{LO}$ to the category of sets.
\end{pr}

For any orthogonal $\fl$-module $(\rho,\fa)$ and any $L\in \fl$ the element
$I_L=(e^{\ad(L)},e^{-\rho(L)})$ is an automorphism of the pair $(\fl,\fa)$.
By definition, the group of inner automorphisms of $(\fl,\fa)$ is the group generated
by the elements $I_L$, $L\in\fl$. It is well known that inner
automorphisms act trivially on usual Lie algebra cohomology.
We now show that the same is true for quadratic cohomology.

\begin{pr}\label{inner}
Let $(\rho,\fa)$ be an orthogonal $\fl$-module. Then for any $L\in \fl$ and $p\in 2\NN_0$
$$ I_L^*=\textstyle{\id_{\cH^p_Q(\fl,\fa)}}\ .$$
\end{pr}
\proof We look at the one parameter group
$$t\mapsto
\Phi_t=((e^{-\ad(tL)},e^{\rho(tL)})^*,(e^{-\ad(tL)})^*)$$
acting on
$\cC^p_Q(\fl,\fa)$.
Let $X_L$ be the corresponding vector field on $\cC^p_Q(\fl,\fa)$,
$X_L(\alpha,\gamma):={\frac{d}{dt}|}_{t=0} \Phi_t(\alpha,\gamma)$.
Inserting $L$ in the first place of a cochain defines maps $i_L: C^q(\fl,\fa)\rightarrow 
C^{q-1}(\fl,\fa)$.
By the well known homotopy formula for the Lie derivative we have
$$ X_L(\alpha,\gamma)=((d\circ i_L+i_L\circ d)\alpha, (d\circ i_L+i_L\circ d)\gamma) \ .$$
Set $\tau_L=i_L(\alpha)$, $\sigma_L=i_L(\gamma)$.
If $(\alpha,\gamma)\in\cZ^p_Q(\fl,\fa)$, then
$$ X_L(\alpha,\gamma)=(d\tau_L,d\sigma_L+\textstyle\frac{1}{2}i_L\langle \alpha\wedge\alpha 
\rangle)
=(d\tau_L,d\sigma_L+\langle \tau_L\wedge\alpha \rangle)\ .$$
This formula shows that $X_L$ is tangential to $\cC^{p-1}_Q(\fl,\fa)$-orbits
in $\cZ^p_Q(\fl,\fa)$. Therefore, $\Phi_t$ maps each orbit to itself.
In other words, the induced action of $\Phi_t$ on cohomology is trivial.
Since $I_L^*$ is induced by $\Phi_{-1}$ the proposition follows.
\qed

\begin{de}\label{pro}
The direct sum of two pairs $(\fl_i,\fa_i)$, $i=1,2$, is defined by
$$ (\fl,\fa)=(\fl_1,\fa_1)\oplus (\fl_2,\fa_2):= (\fl_1\oplus \fl_2,\fa_1\oplus \fa_2)\ ,$$
where the direct sum of $\fa_1$ and $\fa_2$ is orthogonal, and for $i\ne j$ the Lie algebra
$\fl_i$ acts trivially
on $\fa_j$. A direct sum is called non-trivial if both summands are
different from the pair $(0,0)$.
\end{de}
The following lemma can be easily verified.

\begin{lm}\label{pru}
Let $(\fl,\fa)=(\fl_1,\fa_1)\oplus (\fl_2,\fa_2)$ be the direct sum of two pairs.
Let $q_i:\fl\rightarrow \fl_i$ be the projection, and let $j_i: \fa_i\rightarrow \fa$
be the injection.
Then there is a natural injective map
$$  +\,:\,\left( (q_1,j_1)^*\cH^p_Q(\fl_1,\fa_1)\right) \times \left( 
(q_2,j_2)^*\cH^p_Q(\fl_2,\fa_2)\right) \longrightarrow \cH^p_Q(\fl,\fa)$$
induced by addition in $\cZ^p_Q(\fl,\fa)\subset Z^p(\fl,\fa)\oplus C^{2p-1}(\fl)$.
\end{lm}

We conclude this section by clarifying the relationship between
$\cH^p_Q(\fl,\fa)$ and the cohomology groups $H^*(\fl,\fa)$ and
$H^{*}(\fl)$. Recall the definition (\ref{stokes}) of the cup product $\cup$.
We consider the subvariety $H^{p}_\cup(\fl,\fa)\subset
H^p(\fl,\fa)$
defined by
$$ H^{p}_\cup(\fl,\fa):=\{a \in H^p(\fl,\fa)\:|\: a\cup a= 0\}\ .$$
The map $\tilde p: \cZ^p_Q(\fl,\fa)\rightarrow H^{p}(\fl,\fa)$,
$(\alpha,\gamma)\mapsto [\alpha]$, is constant along
$\cC^{p-1}_Q(\fl,\fa)$-orbits and has image $H^{p}_\cup(\fl,\fa)$.
Thus it induces a surjective map
$ p: \cH^{p}_Q(\fl,\fa)\rightarrow H^{p}_\cup(\fl,\fa)$.

\begin{pr}\label{club}
The map $p: \cH^{p}_Q(\fl,\fa)\rightarrow H^{p}_\cup(\fl,\fa)$
gives rise to a partition of $\cH^{p}_Q(\fl,\fa)$
$$ \cH^p_Q(\fl,\fa)=\coprod_{a\in H^p_\cup(\fl,\fa)} p^{-1}(a)$$
into affine spaces $p^{-1}(a)$ with associated vector spaces
$H^{2p-1}(\fl)/(a\cup H^{p-1}(\fl,\fa))$.
The affine structure of $p^{-1}(a)$ is given by the formula
\begin{equation}\label{affin}
p^{-1}(a)\times H^{2p-1}(\fl)/(a\cup H^{p-1}(\fl,\fa))\ni 
([\alpha,\gamma],[\delta]) \longmapsto 
[\alpha,\gamma+\delta] \in 
p^{-1}(a)\ .
\end{equation}
Here $(\alpha,\gamma)\in \tilde p^{-1}(a)$ and 
$\delta \in Z^{2p-1}(\fl)$.
\end{pr}
\proof
We have to check that the action (\ref{affin}) is well-defined and simply 
transitive. First we observe that the abelian group
$Z^{2p-1}(\fl)$ acts on $\tilde p^{-1}(a)$ by 
\begin{equation}\label{affe}
(\alpha,\gamma)\delta\mapsto
(\alpha,\gamma+\delta)\ .
\end{equation} 
This action commutes with the action of $\cC^{p-1}_Q(\fl,\fa)$. The resulting 
action of 
$\cC^{p-1}_Q(\fl,\fa)\times Z^{2p-1}(\fl)$ on $\tilde p^{-1}(a)$ is transitive.
Therefore we obtain a well-defined transitive action of $Z^{2p-1}(\fl)$ on
$p^{-1}(a)$. We compute its kernel. We have
$$ [\alpha,\gamma+\delta]=[\alpha,\gamma]\ \Leftrightarrow\ 
\exists\, \tau\in Z^{p-1}(\fl,\fa), \sigma\in C^{2p-2}(\fl)\, \mbox{ s.th. }
\delta=\langle \alpha\wedge\tau \rangle + d\sigma\ .$$
Thus (\ref{affe}) induces a simply transitive action of
$$ Z^{2p-1}(\fl)/(\langle \alpha\wedge Z^{p-1}(\fl,\fa)\rangle+ B^{2p-1}(\fl))
\cong  H^{2p-1}(\fl)/(a\cup H^{p-1}(\fl,\fa)) $$
on $p^{-1}(a)$ which coincides with (\ref{affin}).
\qed 
 
%%%%%%%%%%%%%%%%%%%%%%%%%%%%%%%%%%%%%%%%%%%%%%%%%%%%%%%%%%%%%%%%%%%%%%%%%%%%%%%%

\section{Quadratic extensions}
\label{Sdd}
In this section we study a two step extension procedure, called 
quadratic extension, of Lie algebras
by orthogonal modules resulting in metric Lie algebras. For a fixed 
Lie algebra $\fl$ and an orthogonal module $\fa$ we establish a 
bijection between
equivalence classes of quadratic extensions of $\fl$ by $\fa$ and 
the cohomology set $\cH^2_Q(\fl,\fa)$. In particular, for any 
$(\alpha,\gamma)\in \cZ^2_Q(\fl,\fa)$ we construct
a metric Lie algebra which
has the structure of a quadratic extension.

\subsection{Definition}\label{leonce}

Let $\fl$ be a Lie algebra and $(\rho,\fa,\ipa)$ an orthogonal
$\fl$-module. We will also consider $\fa$ as an abelian
metric Lie algebra.

\begin{de}
A quadratic extension of $\fl$ by $\fa$ is given by a quadrupel
$(\fg,\fri,i,p)$, where
\begin{itemize}
\item $\fg$ is a metric Lie algebra
\item $\fri\subset\fg$ is an isotropic ideal and
\item $i$ and $p$ are Lie algebra homomorphisms constituting
an exact sequence of Lie algebras
\begin{equation}\label{cruc}
 0\rightarrow \fa\stackrel{i}{\longrightarrow}\fg/\fri
\stackrel{p}{\longrightarrow} \fl\rightarrow 0\ ,
\end{equation}
which is consistent with the representation $\rho$ of $\fl$ on $\fa$ in the following sense:
\begin{equation}\label{Ei}
    i(\rho(L)A)=[\tilde L,i(A)]\  \in\ i(A)
\end{equation} 
holds for all $\tilde L \in\fg/\fri$ with $p(\tilde L)=L$.
In addition we require that $\im i=\fri^\perp/\fri$ and
that
  $i: \fa\rightarrow \fri^\perp/\fri$ is an isometry.
  \end{itemize}
\end{de}

Recall that an isotropic ideal of a metric Lie algebra is always abelian, hence $\fri\subset\fg$ is abelian.

If $\fg$ is a metric Lie algebra with an isotropic
ideal $\fri\in\fg$ such that $\fri^\perp/\fri$ is abelian, then
the sequence
\begin{equation}\label{ifix} 
0\rightarrow \fri^\perp/\fri\stackrel{i}{\longrightarrow}\fg/\fri
\stackrel{p}{\longrightarrow} \fg/\fri^\perp\rightarrow 0 
\end{equation}
defines a quadratic extension of $\fg/\fri^\perp$ by the orthogonal
module $\fri^\perp/\fri$. We call (\ref{ifix}) the canonical extension
associated with $(\fg,\fri)$.

Let $\tilde p:\fg\rightarrow\fl$ be the composition 
of the
natural projection $\fg\rightarrow \fg/\fri$ with $p$. 
Now let $p^*:=\fl^*\rightarrow\fg$ be the dual map of $\tilde p$, where we identify $\fg^*$ 
with $\fg$ using the non-degenerate inner product on $\fg$. This homomorphism is injective
since
$\tilde p$ is surjective. Its image
equals $(\ker \tilde p)^\perp=\fri$.  

Using the homomorphism $p^*$ we see that
a quadratic extension
determines a second exact sequence of Lie algebras
\begin{equation}\label{ofox}
0\rightarrow \fl^*{\longrightarrow}\fg
{\longrightarrow} \fg/\fri\rightarrow 0\ ,
\end{equation}
where we consider $\fl^*$ as abelian Lie algebra.

Equations (\ref{cruc}) and (\ref{ofox}) show that a metric Lie algebra which is
a quadratic extension of $\fl$
by $\fa$ can be considered as the result of two subsequent extensions
of Lie algebras which satisfy certain compatibility conditions: first
an extension of $\fl$ by $\fa$ and second an extension of the 
resulting
Lie algebra by $\fl^*$. This is the point of view taken in 
\cite{KO03},
where the equivalent notion of a twofold extension was studied (for 
abelian
$\fl$).

\begin{lm}\label{rabe}
Let $(\fg,\fri,i,p)$ be a quadratic extension of $\fl$ by $\fa$.
Then $\fg$ does not contain a simple ideal.
\end{lm} 
\proof
Assume that $\fs$ is a simple ideal in $\fg$. Then $\fs$ does not contain non-zero abelian ideals. Hence,
$\fri\cap\fs=0$. In particular, $[\fri,\fs]=0$. This implies
$$\langle \fri,\fs\rangle = \langle\fri,[\fs,\fs]\rangle = \langle[\fri,\fs],\fs\rangle =0,$$
thus $\fs\subset\fri^\perp$. We conclude $\fs\subset\fri^\perp/\fri$, which contradicts the 
assumption that $\fa$ is abelian.
\qed

\begin{de}
Two quadratic extensions $(\fg_j,\fri_j,i_j,p_j)$, $j=1,2$, of $\fl$
by $\fa$ are called to be equivalent if there exists an isomorphism
of metric Lie algebras
$\Psi: \fg_1 \rightarrow \fg_2$
which maps $\fri_1$ onto $\fri_2$ and satisfies
$$ \overline{\Psi}\circ i_1=i_2\qquad\mbox{ and }\qquad
p_2\circ\overline{\Psi}=p_1\ ,$$
where $\overline{\Psi}:\fg_1/\fri_1\rightarrow\fg_2/\fri_2$ is the
induced map.
\end{de}

\subsection{The standard model}

Let $\fl$ be a Lie algebra and let $(\rho,\fa,\ipa)$ be an orthogonal $\fl$-module. We 
choose $\alpha\in C^{2}(\fl,\fa)$ and $\gamma\in
C^{3}(\fl)$. We consider the vector space $\fd:=\fl^{*}\oplus \fa
\oplus \fl$ and define an inner product $\ip$ on $\fd$ by
$$\langle Z_1+A_1+L_1, Z_2+A_2+L_2\rangle=\langle A_1,A_2\rangle_{\fa}+ Z_1(L_2) + 
Z_2(L_1)$$
for all $Z_1,Z_2\in\fl^*$, $A_1,A_2\in\fa$ and  $L_1,L_2\in\fl$.
\begin{pr}\label{Lemma1}
    \begin{itemize}
    \item[(i)]
    There exists a unique antisymmetric bilinear map $\lb:
    \fd\times\fd \rightarrow\fd$ which satisfies
    \begin{eqnarray}
        \label{c1}
        && \ip \mbox{ is invariant, i.e. } \langle[X,Z],Y\rangle =
        \langle X,[Z,Y]\rangle \mbox{ for all } X,Y,Z\in\fd, \\
        \label{c2}
        && [\fd,\fl^{*}]\subset \fl^{*},\ [\fl^{*},\fl^{*}]=0,\
        [\fa,\fa]\subset\fl^{*},\\
        \label{c3}
        && [L_{1},L_{2}]=\gamma(L_{1},L_2,\cdot) +
        \alpha(L_{1},L_{2}) + [L_{1},L_{2}]_\fl \quad\mbox{ for all } L_1,L_2\in\fl,\\
        \label{c4}
        && \langle [L,A_{1}], A_{2}\rangle = \langle
        \rho(L)A_1,A_2\rangle\quad \mbox{ for all } L\in\fl,\ A_1,A_2\in\fa.
    \end{eqnarray} 
    \item[(ii)]
    The triple $\dd :=(\,\fd,\,\lb,\,\ip\,)$ is a metric Lie algebra if and
    only if $(\alpha,\gamma)\in\cZ_{Q}^{2}(\fl,\fa).$
    \end{itemize}
\end{pr}
\proof
Assume $\lb$ is an antisymmetric bilinear map satisfying Equations
(\ref{c1}) to (\ref{c4}). By $[\fl^{*},\fd]\subset \fl^{*}$ and the
invariance of $\ip$ we have $[(\fl^{*})^{\perp},\fd]\subset
(\fl^{*})^{\perp}$, thus $[\fl^{*}+\fa,\fd]\subset \fl^{*}+\fa$.
Hence, $[L,A]$ is in $\fl^{*}+\fa$ for all $L\in\fl$ and $A\in\fa$.
By the invariance of $\ip$ and Equation (\ref{c3}) we obtain
$$\langle [L,A],L'\rangle = -\langle A, [L,L']\rangle =-\langle
A,\alpha(L,L')\rangle$$
for all $L'\in \fl$.
Together with Equation (\ref{c4}) this yields
\begin{equation}
    \label{c5}
    [L,A]=\rho(L)A - \langle A,\alpha(L,\cdot)\rangle.
\end{equation}
Similarly, we obtain
\begin{eqnarray}
        \label{c6}
        &&[L,Z]=\ad{}^{*}(L)Z\\
        \label{c7}
        &&[A_{1},A_{2}]=\langle \rho(\cdot) A_{1},A_{2}\rangle\\
        \label{c8}
        &&[A,\sigma]=0
\end{eqnarray}
for all $L\in\fl,\ Z\in\fl^{*}$ and $\ A,A_{1},A_{2}\in\fa$.
Hence, if $\lb$ is an antisymmetric bilinear map satisfying
(\ref{c1}) to (\ref{c4}), then it is uniquely determined. On the
other hand we can define an antisymmetric bilinear map by Equations
(\ref{c3}) to (\ref{c8}) and $[\fl^{*},\fl^{*}]=0$. Then $\ip$ is
invariant and Equation (\ref{c2}) holds. This proves the first
assertion of the lemma.

The triple $\dd$ is a metric Lie algebra if and only if $\lb$ as
defined in {\it (i)} satisfies the Jacobi identity. Obviously, for any
$\alpha\in C^{2}(\fl,\fa)$ and $\gamma\in C^{3}(\fl)$ the map $\lb$
defined above satisfies
$$[\fl,[\fl^{*}+\fa,\fl^{*}]]=[\fl^{*}+\fa,[\fl^{*},\fl]]=
[\fl^{*},[\fl,\fl^{*}+\fa]]=0$$
and
$$[\fl^{*}+\fa,[\fl^{*}+\fa,\fl^{*}+\fa]]=0.$$
Furthermore, it satisfies
\begin{eqnarray*}
    \lefteqn{[L_{1},[L_{2},Z]]+[L_{2},[Z,L_{1}]]+
    [Z,[L_{1},L_{2}]]}\\
    &&= \ad{}^{*}(L_{1})\ad{}^{*}(L_{2})Z -
\ad{}^{*}(L_{2})\ad{}^{*}(L_{1})Z -\ad{}^{*}([L_{1},L_{2}])Z =0
\end{eqnarray*}
and
\begin{eqnarray*}
    \lefteqn{[A_{1},[A_{2},L]]+[A_{2},[L,A_{1}]]+
    [L,[A_{1},A_{2}]]}\\
    &&= \langle \rho(\cdot)A_{1},-\rho(L)A_{2}\rangle+
    \langle \rho(\cdot)A_{2},\rho(L)A_{1}\rangle
    -\langle \rho([L,\cdot\,])A_{1},A_{2}\rangle=0
\end{eqnarray*}
since $\rho$ is an orthogonal representation. We have to prove that
the remaining identities
\begin{eqnarray}
    &&\cyclsum[L_{1},[L_2,L_{3}]]=0 \label{z1}\\
    &&[A,[L_{1},L_{2}]]+[L_{1},[L_{2},A]]+[L_{2},[A,L_{1}]]=0\label{z2}
\end{eqnarray}
for $L_{1},L_{2},L_{3}\in \fl$ and $A\in\fa$ are equivalent to
the condition $(\alpha,\gamma)\in\cZ_{Q}^{p}.$  Here
$\cyclsum$ denotes the sum over all cyclic permutations of
$L_{1},L_{2}$ and $L_3$. Because of
\begin{eqnarray*}
    \lefteqn{\cyclsum[L_{1},[L_2,L_{3}]]
    =\cyclsum[L_{1},\gamma (L_{2},L_{3},\cdot\,)
    +\alpha(L_{2},L_{3}) +[L_{2},L_{3}]_{\fl}]}\\
&&   =\cyclsum\Big( -\gamma(L_{2},L_{3},[L_{1},\cdot\,]) +
\rho(L_{1})\alpha(L_{2},L_{3}) -
\langle\alpha(L_{2},L_{3}),\alpha(L_{1},\cdot\,)\rangle\\
&&\hspace{2.4cm} +
\gamma(L_{1},[L_{2},L_{3}]_{\fl},\cdot\,)
+\alpha(L_{1},[L_{2},L_{3}]_{\fl}) +
[L_{1},[L_{2},L_{3}]_{\fl}]_{\fl}\,\Big)\\
&&= (d\gamma-\frac  12 \langle \alpha\wedge\alpha\rangle )
(L_{1},L_{2},L_{3},\cdot\,) + d\alpha(L_{1},L_{2},L_{3})
\ \in\,\fl^{*}\oplus\fa
\end{eqnarray*}
Equation (\ref{z1}) is equivalent with $d\gamma=\frac  12 \langle
\alpha\wedge\alpha\rangle$ and $d\alpha=0$. Similarly one proves
$$    [A,[L_{1},L_{2}]]+[L_{1},[L_{2},A]]+
    [L_{2},[A,L_{1}]]
    =\langle  A, -d\alpha(L_{1},L_{2},\cdot\,)\rangle.$$
Hence, Equation (\ref{z2}) is equivalent to $d\alpha=0$. This proves
the second assertion of the lemma.
\qed

We identify $\fd/\fl^{*}$ with $\fa\oplus\fl$ and denote by
$i:\fa\rightarrow \fa\oplus\fl$ the injection and by $p:\fa\oplus\fl\rightarrow
\fl$ the projection. Then the following proposition is obvious.
\begin{pr} If $(\alpha,\gamma) \in \cZ^{2}_{Q}(\fl,\fa)$, then
    the quadrupel  $(\dd,\fl^{*},i,p)$ is a quadratic extension of 
    $\fl$ by $\fa$.
\end{pr}
We will denote the quadratic extension $(\dd,\fl^{*},i,p)$ also by
$\dd$.
\begin{re}\label{dprime}
{\rm Let $(\alpha,\gamma)\in\cZ^{2}_{Q}(\fl,\fa)$ be a cocycle and
$\dd=(\,\fd,\lb,\ip\,)$ the associated metric Lie algebra constructed above.
Furthermore, let $\ip_{\fl}$ be an invariant not necessarily non-degenerate
inner product on $\fl$. We define a new scalar product $\ip'$ on $\fd$ by
$\ip'=\ip\oplus\ip_\fl$. Then
$\fd'_{\alpha,\gamma}(\fl,\ip_{\fl},\fa,\rho)=(\,\fd,\lb,\ip'\,)$ is also a metric 
Lie algebra.
Moreover, $(\fd'_{\alpha,\gamma}(\fl,\ip_{\fl},\fa,\rho),\fl^{*},i,p)$ with $i$ and
$p$ as above is a quadratic extension of $\fl$ by $\fa$
and
\begin{eqnarray*}
\fd'_{\alpha,\gamma}(\fl,\ip_{\fl},\fa,\rho)&\longrightarrow&
\fd_{\alpha,\gamma-\frac12\langle[\cdot\,,\cdot]_\fl,\cdot\,\rangle}
(\fl,\fa,\rho)\\
Z+A+L&\longmapsto&Z+A+L+\frac12\langle L,\cdot\,\rangle_\fl
\end{eqnarray*}
for $Z\in\fl^*$, $A\in\fa$, $L\in\fl$
is an equivalence of quadratic extensions.

}
\end{re}
\subsection{Classification by cohomology}

\begin{pr}\label{P}
    For $(\alpha_{1},\gamma_{1}),\ (\alpha_{2},\gamma_{2})\in 
    \cZ^{2}_{Q}(\fl,\fa)$ the quadratic extensions
    $\fd_{\alpha_{1},\gamma_{1}}(\fl,\fa,\rho)$ and
    $\fd_{\alpha_{2},\gamma_{2}}(\fl,\fa,\rho)$ of $\fl$ by $\fa$ are
    equivalent 
    if and  only if $[\alpha_{1},\gamma_{1}]=[\alpha_{2},\gamma_{2}]\in 
    \cH^{2}_{Q}(\fl,\fa)$.
\end{pr}
\proof A linear map $\Psi:\fd_{\alpha_{1},\gamma_{1}}(\fl,\fa,\rho)
\rightarrow\fd_{\alpha_{2},\gamma_{2}}(\fl,\fa,\rho)$ is an
equivalence of quadratic extensions if and only if
\begin{itemize}
    \item[(i)]$\Psi(\fl^{*})=\fl^{*}$, $\proj_{\fa}\Psi|_{\fa}=\Id$, 
$\proj_{\fl}\Psi|_{\fl}=\Id$,
    \item[(ii)]$\Psi$ is an isometry, and
    \item[(iii)]$\Psi$ is a Lie algebra isomorphism.
\end{itemize}    
Condition (i) is equivalent to
\begin{equation}\label{EM}
    \Psi=\left(
\begin{array}{ccc}
    \psi & \eta & \xi  \\
    0 & \Id & \tau  \\
    0 & 0 & \Id
\end{array}\right) :\fl^{*}\oplus\fa\oplus\fl \longrightarrow 
\fl^{*}\oplus\fa\oplus\fl
\end{equation}
for linear maps $\psi:\fl^{*}\rightarrow\fl^{*},\ \eta:\fa 
\rightarrow\fl^{*},\ \xi:\fl\rightarrow\fl^{*},\ \tau:\fl\rightarrow\fa$.

Conditions (i) and (ii) are satisfied if and only if $\Psi$ is as 
above and the equations
\begin{eqnarray*}
\langle Z,L\rangle \ =&\langle\Psi Z,\Psi L\rangle&=\
\langle\psi Z,\xi L+\tau L +L\rangle = \langle \psi
Z,L\rangle\\
0\ = &\langle\Psi A,\Psi L\rangle&= \ \langle \eta A +A,\xi L+\tau L 
+L\rangle \,=\, \langle \eta A,L\rangle + \langle A,\tau L\rangle\\
0\ = &\langle\Psi L_{1},\Psi L_{2}\rangle&=  \ \langle \xi L_{1}+\tau  
L_{1} +L_{1},\xi L_{2}+\tau L_{2}
+L_{2}\rangle \\&&=\ \langle \xi L_{1},L_{2}\rangle + \langle L_{1},\xi 
L_{2}\rangle +\langle\tau L_{1},\tau L_{2}\rangle.
\end{eqnarray*}    
hold. Let $\xi^{*}:\fl\rightarrow \fl^{*}$ and $\tau^{*}:\fa\rightarrow 
\fl^{*}$ be the dual maps of $\xi$ and $\tau$, respectively. These 
maps are given by $\langle \xi^{*}L_{1},L_{2}\rangle  = \langle 
L_{1},\xi L_{2}\rangle$ and $\langle \tau^{*}A,L\rangle  = \langle 
A,\tau L\rangle$ for $A\in\fa$, $L,L_{1},L_{2}\in\fl$. 
Then the last equation says that the selfdual part $\frac 12
(\xi+\xi^{*})$ of $\xi$ equals $-\frac 12 \tau^{*}\tau$.

Consequently, 
Conditions (i) and (ii) are satisfied if and only if $\Psi$ is as 
in (\ref{EM}) with $\psi=\Id$, $\eta=-\tau^{*}$ and $\xi=\bar \sigma -
\frac 12 \tau^{*}\tau$ for an anti-selfdual map $\bar \sigma:\fl
\rightarrow \fl^{*}$, i.e. if and 
only if
\begin{equation}\label{EM1}
    \Psi=\Psi(\tau,\sigma):=\left(
\begin{array}{ccc}
    \Id & -\tau^{*} & \bar\sigma -\frac 12\tau^{*}\tau  \\
    0 & \Id & \tau  \\
    0 & 0 & \Id
\end{array}\right) :\fl^{*}\oplus\fa\oplus\fl \longrightarrow 
\fl^{*}\oplus\fa\oplus\fl,
\end{equation}
where $\tau\in C^{1}(\fl,\fa)$ and  $\sigma (\cdot\,,\cdot)= \langle\bar\sigma
(\cdot),\cdot\rangle \in C^{2}(\fl)$.

Now we consider Condition (iii). We denote the Lie brackets on 
$\fd_{\alpha_{1}, \gamma_{1}}(\fl,\fa,\rho)$ and
$\fd_{\alpha_{2},\gamma_{2}}(\fl,\fa,\rho)$ by $\lb_{1}$ and
$\lb_{2}$, respectively. Assume $\Psi$ is given as in (\ref{EM1}). 
Then $\Psi$  is a Lie algebra
isomorphism if and only if the Lie bracket $\lb'$ defined by 
$[X,Y]':=\Psi^{-1}[\Psi 
X,\Psi Y]_{2}$ for all $X,Y\in\fd$ is equal to $\lb_{1}$. By 
Proposition \ref{Lemma1} this is the case if and only if $\lb'$ 
satisfies 
\begin{eqnarray*}
      && \langle[X,Z]',Y\rangle = 
        \langle X,[Z,Y]'\rangle \mbox{ for all } X,Y,Z\in\fd, \\
      && [\fl^{*},\fd]'\subset \fl^{*},\ [\fl^{*},\fl^{*}]'=0,\
        [\fa,\fa]'\subset\fl^{*},\\
      && [L_{1},L_{2}]'=\gamma_{1}(L_{1},L_2,\cdot) + 
        \alpha_{1}(L_{1},L_{2}) + [L_{1},L_{2}]_\fl\\
      && \langle [L,A_{1}]', A_{2}\rangle = \langle 
        \rho(L)A_1,A_2\rangle   
\end{eqnarray*}
for all $A_{1},A_{2}\in\fa$ and $L,L_{1},L_{2}\in\fl$. Obviously, the first two of 
these conditions and the last one are always satisfied 
if we choose  $\Psi$ as in (\ref{EM1}). The third condition is equivalent to the 
following equations
\begin{eqnarray*}
    \langle  [L_{1},L_{2}]',L_{3}\rangle\ =& \langle [\Psi L_{1},\Psi 
    L_{2}]_{2},\Psi L_{3} \rangle &=\ \gamma_1(L_{1},L_{2},L_{3})\\
    \langle  [L_{1},L_{2}]',A\rangle\ =& \langle [\Psi L_{1},\Psi 
    L_{2}]_{2},\Psi A \rangle &=\ \langle \alpha_1(L_{1},L_{2}),A\rangle\\
    \langle  [L_{1},L_{2}]',Z\rangle\ =& \langle [\Psi L_{1},\Psi
    L_{2}]_{2},\Psi Z \rangle &=\ \langle
    [L_{1},L_{2}]_{\fl},Z\rangle
\end{eqnarray*}
for all $Z\in\fl^{*}$, $A\in\fa$ and $L_{1},L_{2},L_{3}\in\fl$.
The third equation is always satisfied if $\Psi$ is an isometry as 
in (\ref{EM1}). Hence $\Psi$ is a Lie algebra isomorphism if and only 
if 
\begin{eqnarray}
    \label{Ea}
   \langle [\Psi L_{1},\Psi 
    L_{2}]_{2},\Psi L_{3} \rangle &=& \gamma_1(L_{1},L_{2},L_{3})\\
    \label{Eb}
    \langle [\Psi L_{1},\Psi 
    L_{2}]_{2},\Psi A \rangle &=& \langle \alpha_1(L_{1},L_{2}),A\rangle 
\end{eqnarray}    
for all $A\in\fa$ and $L_{1},L_{2},L_{3}\in\fl$. By definition of 
$\Psi$ Equation (\ref{Ea}) is equivalent to
\begin{eqnarray*}
    \lefteqn{\gamma_{1}(L_{1},L_{2},L_{3})=}\\
    &&=\Big\langle\,[\bar\sigma(L_{1})-\frac
    12\tau^{*}\tau(L_{1})+\tau (L_{1})+L_{1}, \bar\sigma (L_{2})-\frac
    12\tau^{*}\tau(L_{2})+\tau( L_{2})+L_{2}]_{2},\\
    &&\hspace{7.3cm}\bar\sigma(L_{3})-\frac
    12\tau^{*}\tau(L_{3})+\tau( L_{3})+L_{3}\,\Big\rangle\\
    &&=\Big\langle \,
    \sigma(L_1,[L_2,\cdot\,]_{\fl})-\sigma(L_2,[L_1,\cdot\,]_{\fl}) - \frac 12
    \tau^{*}\tau(L_{1})([L_{2},\cdot\,]_{\fl})+ \frac 12 
    \tau^{*}\tau(L_{2})([L_{1},\cdot\,]_{\fl}) \\
    &&\hspace{0.7cm} +\langle\rho(\cdot)\tau(L_{1}),\tau(L_2)\rangle 
    +\langle\tau (L_{1}),\alpha_{2}(L_2,\cdot\,) \rangle -\langle\tau 
    (L_{2}),\alpha_{2}(L_1,\cdot\,) \rangle \\
    &&\hspace{0.7cm}+L_{1}\tau (L_{2}) -L_2\tau 
    (L_{1}) +\gamma_{2}(L_{1},L_{2},\cdot\,) + \alpha_{2}(L_{1},L_{2}) 
    +[L_{1},L_{2}]_{\fl},\,\\
    &&\hspace{7.3cm}
    \bar\sigma(L_{3})-\frac
    12\tau^{*}\tau(L_{3})+\tau( L_{3})+L_{3}\,\Big\rangle\\
    &&=\gamma_{2}(L_{1},L_{2},L_{3}) +\cyclsum 
    \sigma(L_{1},[L_{2},L_{3}]_{\fl})  +\cyclsum \langle
    \alpha_{2}(L_{1},L_{2}),\tau(L_{3})\rangle\\
    &&\hspace{0.7cm} +\frac 12 \cyclsum \Big(\,\langle \tau (L_{1}),L_2\tau 
    (L_{3})\rangle - \langle \tau (L_{1}),L_3\tau (L_{2})\rangle -
    \langle \tau (L_1),[L_2,L_{3}]_{\fl}\rangle\,\Big)\\
    &&= (\gamma_{2}+d\sigma +\langle \alpha_{2}\wedge\tau\rangle
    +\frac 12\langle 
    \tau\wedge d\tau\rangle)(L_{1},L_{2},L_{3}).
\end{eqnarray*}   
Moreover, using the previous calculation of $[\Psi (L_{1}),\Psi 
(L_{2})]_{2}$ we see that (\ref{Eb}) is equivalent to
\begin{eqnarray*}
    \langle \alpha_{1}(L_{1},L_{2}),A\rangle &=&\langle [\Psi (L_{1}),\Psi 
(L_{2})]_{2}, \Psi(A) \rangle\,=\, \langle [\Psi (L_{1}),\Psi
(L_{2})]_{2}, A-\tau^{*}A \rangle\\
&=&\langle L_{1}\tau(L_{2})-L_{2}\tau(L_{1}) 
+\alpha_{2}(L_{1},L_{2}),A\rangle - \langle 
[L_{1},L_{2}]_{\fl},\tau^{*}(A)\rangle\\
&=&\langle  \alpha_{2}(L_{1},L_{2})+d\tau(L_{1},L_{2}),A\rangle.
\end{eqnarray*} 
Hence, $\Psi$ is a Lie algebra isomorphism if and only if 
$(\alpha_{1},\gamma_{1})=(\alpha_{2},\gamma_{2})(\tau,\sigma)$.

Now let us finish the proof of the proposition. If 
$\fd_{\alpha_{1},\gamma_{1}} (\fl,\fa,\rho)$ and
$\fd_{\alpha_{2},\gamma_{2}}(\fl,\fa,\rho)$ are equivalent, then we
can choose an equivalence map $\Psi$, which can be written as in
(\ref{EM1}) and $(\alpha_{1},\gamma_{1})=(\alpha_{2},\gamma_{2})(\tau,\sigma)$
holds. Thus $[\alpha_{1},\gamma_{1}]=[\alpha_{2},\gamma_{2}]\in 
\cH^{2}_{Q}(\fl,\fa)$. Conversely, if there exists an element  $(\tau,\sigma)$
in $\cC^{1}_{Q}(\fl,\fa)$ such that $(\alpha_{1},\gamma_{1})= 
(\alpha_{2},\gamma_{2})(\tau,\sigma)$, then 
$\Psi(\tau,\sigma):\fd_{\alpha_{1},\gamma_{1}}(\fl,\fa,\rho)
\rightarrow\fd_{\alpha_{2},\gamma_{2}}(\fl,\fa,\rho)$  (see (\ref{EM1})) is an equivalence.
\qed

\begin{re}
{\rm
The map
\begin{eqnarray*}
    \overline\Psi:\quad \cC^{1}_{Q}(\fl,\fa)&\longrightarrow& O(\fd,\ip)\\
    (\tau,\sigma)&\longmapsto&\Psi(\tau,\sigma)
\end{eqnarray*}
is an injective group homomorphism. This motivates Definition \ref{unmoti}.
As we have seen above the image of $\overline 
\Psi$ acts on the set of quadratic extensions of $\fl$ by $\fa$ of the form $\dd$ 
by equivalences.
The orbits of this action are exactly the equivalence classes of such quadratic
extensions. Moreover, the map $\cZ^{2}_{Q}(\fl,\fa)\ni(\alpha,\gamma)\mapsto\dd$ is equivariant with respect to the action of $\cC^{1}_{Q}(\fl,\fa)$ on
$\cZ^2_{Q}(\fl,\fa)$ and the action
of $\cC^{1}_{Q}(\fl,\fa)$ on the set $\{\dd\mid (\alpha,\gamma)\in\cZ^2_{Q}(\fl,\fa)\}$ defined by $\overline\Psi$. This motivates Formula (\ref{unmotiviert}).
}
\end{re}

Let $(\fg,\fri,i,p)$ be a quadratic extension of $\fl$ by $\fa$. Then
$\fg\cong \fri\oplus\fa\oplus\fl\cong\fri^{\perp}\oplus\fl$ as vector spaces. 
Since,
furthermore, $\fri$ is isotropic we can choose an isotropic 
complement $V_{\fl}$ of $\fri^\perp$ in $\fg$ and an isomorphism $s:\fl\rightarrow 
V_{\fl}$
such that $\tilde p\circ s=\Id$ holds. Here $\tilde p:\fg\rightarrow\fl$ is the composition of the natural projection $\fg\rightarrow \fg/\fri$ with $p$ as already defined above.

We define  $\alpha\in C^2(\fl,\fa)$ and $\gamma\in C^3(\fl)$ by
\begin{eqnarray}       
i(\alpha(L_1,L_2))&:=&[s(L_1),s(L_2)]-s([L_1,L_2])\, +\, \fri\quad 
\in\,
\fg/\fri
\label{Ealpha}\\
\gamma(L_1,L_2,L_3)&:=& \langle\, [s(L_1),s(L_2)], s(L_3)\rangle \ .
\label{Egamma}
\end{eqnarray}
\begin{pr}\label{Pequiv}
We have that $(\alpha,\gamma)\in \cZ^2_Q(\fl,\fa)$. The quadratic 
extension  $(\fg,\fri,i,p)$ is equivalent to 
$\fd_{\alpha,\gamma}(\fl,\fa,\rho)$.
\end{pr}    
\proof
Here we will denote the inner product and the Lie bracket on $\fg$ by 
$\ip_{\fg}$ and $\lb_{\fg}$, respectively. Let $V_\fa$ be the orthogonal 
complement of $\fri\oplus s(\fl)$ in 
$\fg$. Then $\fri^{\perp}=\fri\oplus V_{\fa}$ and we can define a 
linear map $t:\fa\longrightarrow V_{\fa}$ by
$$i(A)=t(A)+\fri\in \fg/\fri .$$ 
Since $i:\fa\rightarrow \fri^{\perp}/\fri$ is an isometry, also $t$ is 
an isometry, and because of (\ref{Ei}) we have
$$ t(\rho(L)A)\equiv[s(L),t(A)]_\fg\quad \mod \fri. $$
Recall that $p^*:\fl^*\rightarrow\fri$ is an isomorphism
satisfying
\begin{equation}\label{bill}\langle p^*(Z),s(L)\rangle_{\fg}=\langle Z,\tilde p\circ s(L)\rangle_\fg=Z(L)
\end{equation}
for all $Z\in\fl^{*}$ and $L\in\fl$.

Now we consider the triple $\dd=(\,\fd,\lb, \ip\,)$ for $\alpha\in
C^2(\fl,\fa)$ and $\gamma\in C^3(\fl)$ as defined in (\ref{Ealpha})
and (\ref{Egamma}). We define
$$\Psi=p^*+ t+ s: \
\fd=\fl^{*}\oplus\fa\oplus\fl\longrightarrow \fg.$$
By construction $\Psi:(\fd,\ip)\rightarrow (\fg,\ip_{\fg})$ is an isometry.
Next we will show that $\Psi[X,Y]=[\Psi( X),\Psi (Y)]_{\fg}$ holds for
all $X,Y\in\fd$. To do that we define a bilinear map $\lb'$ on
$\fd\times\fd$ by $[X,Y]'=\Psi^{-1}[\Psi(X),\Psi(Y)]_{\fg}$ and prove
that $\lb'$ satisfies (\ref{c1})--(\ref{c4}). Clearly, $\ip$ is
invariant with respect to $\lb'$ since $\ip_{\fg}$ is invariant and
$\Psi$ is an isometry. By construction of $\Psi$ we have
\begin{eqnarray}
    &&[\fl^{*},\fd]'=\Psi^{-1}[\Psi\fl^{*},\Psi\fd]_{\fg}
    =\Psi^{-1}[\fri,\fg]_{\fg}
        \subset\Psi^{-1}(\fri)=\fl^{*}\nonumber\\
    &&[\fl^{*},\fl^{*}]'=\Psi^{-1}[\Psi\fl^{*},\Psi\fl^{*}]_{\fg} =
     \Psi^{-1}[\fri,\fri]_{\fg} =0\label{E2.0}\\
    &&[\fa,\fa]'=\Psi^{-1}[\Psi \fa,\Psi \fa]_{\fg}\subset
    \Psi^{-1}[\fri^{\perp},\fri^{\perp}]_\fg \subset
    \Psi^{-1}(\fri)=\fl^{*},\nonumber
\end{eqnarray}
where we have used the inclusion $[\fri^{\perp},\fri^{\perp}]_\fg\subset \fri$,
which holds since $\fri^{\perp}/\fri\cong \fa$ and $\fa$ is abelian.
Moreover, (\ref{Ealpha}) gives
\begin{equation}\label{nix1}
[s(L_1),s(L_2)]_\fg-t(\alpha(L_1,L_2))-s([L_1,L_2])\in\fri
\end{equation}
and by (\ref{Egamma}) and (\ref{bill}) we have
\begin{eqnarray}
\label{nix2}
&\langle[s(L_1),s(L_2)]_\fg-t(\alpha(L_1,L_2))-s([L_1,L_2]_\fl),s(L_3)\rangle_\fg=
\gamma(L_1,L_2,L_3)&\\
&=\langle p^*(\gamma(L_1,L_2,\cdot)),s(L_3)\rangle_\fg\,.&\nonumber
\end{eqnarray}
Since $p^*(\gamma(L_1,L_2,\cdot))\in\fri$ and $\ip_\fg|_{\fri\times s(\fl)}$ is a
non-degenerate pairing we obtain by (\ref{nix1}) and (\ref{nix2})
\begin{eqnarray*}
[s(L_1),s(L_2)]_\fg&=&p^*(\gamma(L_1,L_2,\cdot))+t(\alpha(L_1,L_2))+s([L_1,L_2]_\fl
)\\
&=&\Psi(\gamma(L_1,L_2,\cdot)+\alpha(L_1,L_2)+[L_1,L_2]_\fl)\,.
\end{eqnarray*}
This yields
\begin{eqnarray}\nonumber
    [L_{1},L_{2}]'&=&  \Psi^{-1}([\Psi(L_1),\Psi(L_2)]_\fg
)=\Psi^{-1}([s(L_1),s(L_2)]_\fg)\\
    &=&\gamma(L_{1},L_{2},\cdot\,)
    +\alpha(L_{1},L_{2}) +[L_{1},L_{2}]_{\fl}.
    \label{E2}
\end{eqnarray}
Finally, we have
\begin{eqnarray*}
    \langle [L,A_{1}]',A_{2}\rangle &=& \langle \Psi(
    [L,A_{1}]'),\Psi(A_{2}) \rangle_\fg = \langle
    [s(L),t(A_{1})]_{\fg},t(A_{2}) \rangle_\fg\\
    &=& \langle t(\rho(L)(A_{1})),t(A_{2}) \rangle_\fg=
    \langle \rho(L)(A_{1}),A_{2} \rangle.
\end{eqnarray*}
By Proposition \ref{Lemma1} this equation together with (\ref{E2.0})
and (\ref{E2}) implies $\lb'=\lb$, thus $\Psi[X,Y]=
[\Psi( X),\Psi (Y)]_{\fg}$ for all $X,Y\in\fd$. In particular, $(\fd,\lb)$
is a Lie algebra since $(\fg,\lb_{\fg})$ is a Lie algebra. Proposition
\ref{Lemma1} now implies $(\alpha,\gamma)\in\cZ^{2}_{Q}(\fl,\fa)$.
Moreover, we conclude that $\Psi:\fd\rightarrow\fg$ is an equivalence
of the quadratic extensions $\fd_{\alpha,\gamma}(\fl,\fa,\rho)$ and
$(\fg,\fri,i,p)$.
\qed

\begin{co}
The cohomology class $[\alpha,\gamma]\in\cH_Q^2(\fl,\fa)$ does not depend
on the choice of $s$.
\end{co}
\proof
Let
$s_{i}:\fl\rightarrow \fg$, $i=1,2$,
be two linear maps with isotropic image and $\tilde p\circ s_{i}=\Id$.
Consider
$(\alpha_{i},\gamma_{i})\in\cZ^{2}_{Q}(\fl,\fa)$, $i=1,2$, as above.
By Proposition \ref{Pequiv} the quadratic extensions
$\fd_{\alpha_{i},\gamma_{i}}(\fl,\fa,\rho)$, $i=1,2$, are equivalent
since both are equivalent to  $(\fg,\fri,i,p)$. Proposition \ref{P}
now implies $[\alpha_{1},\gamma_{1}]=[\alpha_{2},\gamma_{2}]\in \cH_Q^2(\fl,\fa)$.
\qed

We can summarize the results of Section \ref{Sdd} as follows.
\begin{theo}\label{brandt}
The equivalence classes of quadratic extensions of a Lie algebra $\fl$
by an orthogonal module $\fa$ are in one-to-one correspondence with
elements of $\cH_Q^2(\fl,\fa)$.
\end{theo}
%%%%%%%%%%%%%%%%%%%%%%%%%%%%%%%%%%%%%%%%%%%%%%%%%%%%%%%%%%%%%%%%%%%%%%
%%%%%%%%%%%%%%%%%%%%%%%%%%%%%%%%%%%%%%%%%%%%%%%%%%%%%%%%%%%%%%%%%%%%%%

\section{Balanced extensions}\label{baal}

In this section we equip any metric Lie algebra $\fg$ without simple
ideals with the structure of a quadratic extension in a canonical
way, i.e., we construct a canonical isotropic ideal
$\fri(\fg)\subset \fg$ such that $\fri(\fg)^\perp/\fri(\fg)$ is abelian.
A quadratic extension $(\fg,\fri,i,p)$ will be called balanced,
if
$\fri=\fri(\fg)$. The main result
of this section is Theorem \ref{bebe} which describes the subset
of $\cH^2_Q(\fl,\fa)$ corresponding to balanced extensions.

In order to construct the desired canonical ideal we need a little preparation.
Let $\fg$ be a (real, finite-dimensional) Lie algebra, and let $V$ be a
finite-dimensional $\fg$-module.
The socle $S(V)\subset V$ is by definition the maximal submodule of $V$ on
which $\fg$ acts semi-simply. It is equal to the sum of all irreducible submodules
of $V$. There is the dual notion of the radical $R(V)\subset V$
which is the minimal submodule such that $\fg$ acts semi-simply on
$V/R(V)$. For later use we collect the basic computation rules for the
functors $S$ and $R$.
Let $U,V\subset W$ be $\fg$-submodules. Then
\begin{equation}\label{basic}
\begin{array}{lcl}
S(W)=0\Rightarrow W=0&& R(W)=W\Rightarrow W=0\\[0.5ex]
S(S(W))=S(W)&& R(S(W))=0\\[0.5ex]
S(U\cap V)=S(U)\cap V&& R(U\cap V)\subset R(U)\cap R(V)  \\[0.5ex]
S(U+V)\supset S(U)+S(V)&&R(U+V)= R(U)+R(V)\\[0.5ex]
S(W/U)\subset (S(W)+U)/U&&R(W/U)= (R(W)+U)/U \ .
\end{array}
\end{equation}

\begin{de}\label{sockel}
We define the higher socles $S_k(V)\subset V$ and radicals $R_k(V)\subset V$,
$k\in \NN$,
inductively by
\begin{eqnarray*}
S_0(V)&:=&\{0\}\ ,\quad S_k(V):=(p_{k-1})^{-1}(S(V/S_{k-1}(V)))\ ,\\
R_0(V)&:=&V  \ ,\quad R_k(V):= R(R_{k-1}(V))\ ,
\end{eqnarray*}
where $p_{k-1}: V\rightarrow V/S_{k-1}(V)$ is the natural projection.
\end{de}
Clearly, $S_1(V)=S(V)$ and $R_1(V)=R(V)$.

If $V^*$ is the dual $\fg$-module, then for $k\in \NN$
\begin{equation}\label{dual}
S_k(V^*)= R_k(V)^\perp\quad \mbox{and}\quad R_k(V^*)=S_k(V)^\perp\ .
\end{equation}

We are particularly interested in the case of the adjoint representation
$V=\fg$ of $\fg$. In this case Definition \ref{sockel} provides
an increasing and a decreasing chain of ideals of $\fg$
\begin{eqnarray*}
\{0\}=S_0(\fg)\subset S_1(\fg)\subset S_2(\fg)\subset&\dots&\subset
S_{l_+}(\fg)=\fg \\
\fg=R_0(\fg)\supset R_1(\fg)\supset R_2(\fg)\supset&\dots&\supset
R_{l_-}(\fg)=\{0\} \ .
\end{eqnarray*}
We call $R(\fg)$ the radical of nilpotency of $\fg$ in order to distinguish
it from the (solvable) radical $\fr$ and the nilpotent radical (= maximal nilpotent
ideal) $\fn$. Note that for $k>1$ the radical of nilpotency of the Lie algebra
$R_{k-1}(\fg)$ may be larger than
$R_k(\fg)$ since $R_k(\fg)$ is defined
in terms of the $\fg$-module structure of $R_{k-1}(\fg)$. It is a consequence of
Lie's Theorem that (see \cite{bourbaki})
\begin{equation}\label{bour}
R(\fg)=\fr\cap \fg'=[\fr,\fg']\subset \fn
\end{equation}
and that $R(\fg)$ acts trivially on any semi-simple $\fg$-module $V$.
The last property implies (consider $V=S_k(\fg)/S_{k-1}(\fg)$)
\begin{equation}
[R(\fg), S_{k}(\fg)]\subset S_{k-1}(\fg)\ . \label{bour2}
\end{equation}
We will also need the relation of $\fz(\fg)$ with $S(\fg)$.
Of course, $\fz(\fg)\subset S(\fg)$. We will formulate a more precise
result for the case that $\fg$ does not contain simple ideals which is
relevant for quadratic extensions (see Lemma \ref{rabe}). Of course, the
general case does not present essential difficulties since any Lie algebra splits
into a
direct sum of a semi-simple ideal and an ideal which does not contain simple
ideals.
\begin{lm}\label{ralf}
If $\fg$ does not contain a simple ideal, then
\begin{enumerate}
\item[{\rm (a)}] $\displaystyle S(\fg)= \fz(\fg)+S(\fg)\cap R(\fg)$.
\item[{\rm (b)}] $\displaystyle S(\fg)\subset R(\fg)$ if and only if
$\fz(\fg)\subset\fg^\prime$.
\end{enumerate}
\end{lm}
\proof
Let $\fs$ be a $\fg$-invariant complement of $S(\fg)\cap R(\fg)$ in $S(\fg)$. Note
that $\fs$ is a reductive Lie algebra. Let $p:\fg\rightarrow \fg/R(\fg)$ be the
projection, and let $\ft$ be a $\fg$-invariant complement of $p(S(\fg))=p(\fs)$ in
$\fg/R(\fg)$.
We obtain a decomposition of $\fg$ into a direct sum of ideals
$$
\fg = \fs\oplus  p^{-1}(\ft) = \fz(\fs)\oplus \fs^\prime\oplus p^{-1}(\ft)\ .
$$
The decomposition shows that $\fz(\fs)\subset \fz(\fg)$. Now, if $\fg$ does
not contain simple ideals, then $\fs^\prime={0}$.
We obtain $S(\fg)=\fs \oplus S(\fg)\cap R(\fg)\subset \fz(\fg)+ S(\fg)\cap R(\fg)$.
Since the opposite inclusion is obvious this proves (a).

Assume that $S(\fg)\subset R(\fg)$. Using (\ref{bour})
we find $\fz(\fg)\subset S(\fg)\subset R(\fg)\subset \fg^\prime$.
If $\fz(\fg)\subset \fg^\prime$, then $\fz(\fg)\subset \fr\cap \fg^\prime=R(\fg)$.
Now (a) implies that $S(\fg)\subset R(\fg)$.
\qed

\begin{de}\label{meer}
Let $\fg$ be a Lie algebra. We define characteristic ideals
$\fri(\fg)\subset\fj(\fg)\subset\fg$ by
$$\fri(\fg):=\sum_{k=1}^\infty R_k(\fg)\cap S_k(\fg)\quad
\mbox{ and }\ \quad
\fj(\fg):=\bigcap_{k=1}^\infty (R_k(\fg)+ S_k(\fg))\ .$$
\end{de}
Of course, the sum and the intersection are only formally infinite.
For all $j\le k\le l$ we have $R_k(\fg)\cap S_k(\fg)\subset R_j(\fg)$
and $R_k(\fg)\cap S_k(\fg)\subset S_l(\fg)$. Thus $R_k(\fg)\cap S_k(\fg)\subset
R_l(\fg)+S_l(\fg)$ for all $l$. This shows that indeed $\fri(\fg)\subset\fj(\fg)$.

\begin{lm}\label{tau}
The ideals $\fri(\fg)$ and $\fj(\fg)$ satisfy
\begin{enumerate}
\item[{\rm (a)}] $\displaystyle\fj(\fg)=S(\fg)+
\sum_{k=1}^\infty R_k(\fg)\cap S_{k+1}(\fg).$
\item[\rm (b)] The natural representation of $\fg$ on the quotient
$\fj(\fg)/\fri(\fg)$ is semi-simple.
\item[\rm (c)] $[R(\fg), \fj(\fg)]\subset\fri(\fg)$.
\item[\rm (d)] If $\fg$ does not contain a simple ideal, then the Lie algebra
$\fj(\fg)/\fri(\fg)$ is abelian.
\end{enumerate}
\end{lm}
\proof We compute
\begin{eqnarray*}
\fj(\fg)&=&\bigcap_{k=1}^\infty (S_k(\fg)+ R_k(\fg))
= S(\fg) + R(\fg)\cap \bigcap_{k=2}^\infty (S_k(\fg)+ R_k(\fg))\\
&=& S(\fg) + R(\fg) \cap \left( S_2(\fg) + R_2(\fg)\cap \bigcap_{k=3}^\infty
(S_k(\fg)+ R_k(\fg))\right)\\
&=& S(\fg) + R(\fg)\cap S_2(\fg)+R_2(\fg)\cap \bigcap_{k=3}^\infty (S_k(\fg)+
R_k(\fg))\\
&=& \dots\\
&=& S(\fg)+
\sum_{k=1}^\infty R_k(\fg)\cap S_{k+1}(\fg)\ .
\end{eqnarray*}
This shows (a).

In order to prove (b) we show that
$R(\fj(\fg))\subset \fri(\fg)$,
where the radical is taken with respect to the $\fg$-module structure.
Using (a) and the rules (\ref{basic}) we find
\begin{eqnarray*}
R(\fj(\fg))&=& R(S(\fg))+\sum_{k=1}^\infty R\left(R_k(\fg)\cap S_{k+1}(\fg)\right)
\subset \sum_{k=1}^\infty R(R_k(\fg))\cap S_{k+1}(\fg)\\
&=& \sum_{k=1}^\infty R_{k+1}(\fg)\cap S_{k+1}(\fg)\subset \fri(\fg)\ .
\end{eqnarray*}
This proves (b) and
implies that $R(\fg)$ acts trivially on $\fj(\fg)/\fri(\fg)$. Alternatively, (c)
could be shown directly using (a) and (\ref{bour2}).

By (a) we have $\fj(\fg)\subset S(\fg)+R(\fg)$. If $\fg$ does not contain simple
ideals, then Lemma \ref{ralf} implies that $\fj(\fg)\subset \fz(\fg)+R(\fg)$.
Now (d) follows from (c).
\qed

{}From now on let $\fg$ be a metric Lie algebra.  Then by (\ref{dual})
we have
$$
R_k(\fg)^\perp= S_k(\fg) \ .
%\label{dual2}
$$
This implies that $\fri(\fg)^\perp=\fj(\fg)$. In particular, $\fri(\fg)$ is
isotropic.

\begin{de}\label{schwein}
If $\fg$ is a metric Lie algebra, then we call $\fri(\fg)$ its canonical
isotropic ideal.
\end{de}

Observe that both $\fri(\fg)$ and $\fri(\fg)^\perp=\fj(\fg)$ are completely
determined by the Lie algebra structure of $\fg$ and do not depend
on the particular form of the inner product on $\fg$.

Let $\fl$ be a Lie algebra, and let $\fa$ be an
orthogonal $\fl$-module.

\begin{de}
A quadratic extension $(\fg,\fri,i,p)$ of $\fl$ by $\fa$ is called balanced if
$\fri=\fri(\fg)$.
We call the extension regularly balanced, if in addition
$\fz(\fg)\subset\fg^\prime$.
\end{de}

The study of metric Lie algebras is easily reduced to the study of metric
Lie algebras without simple ideals. Indeed, each
simple ideal of a metric Lie algebra is non-degenerate (see e.g. \cite{KO03}, Lemma 2.2). It
follows that any metric Lie algebra
is the direct sum of a semi-simple metric Lie algebra and a metric Lie
algebra without simple ideals.

\begin{pr}\label{schmidt}
    Any metric Lie algebra $\fg$ without simple ideals has the structure of a
balanced
    quadratic extension in a canonical way. It is regularly balanced if and
    only if $\fg$ does not contain a non-degenerate abelian ideal.
\end{pr}

\proof By Part (d) of Lemma \ref{tau} the quotient  $\fri(\fg)^\perp/\fri(\fg)$
is abelian. Thus $\fri(\fg)\subset \fg$ defines a canonical quadratic extension
(see (\ref{ifix})) which is balanced by the very definition.

Any non-degenerate abelian ideal of $\fg$ is central. Since $\fz(\fg)^\perp=\fg'$
such an ideal does not
intersect $\fz(\fg)\cap \fg'$. Vice versa, any complement of $\fz(\fg)\cap \fg'$ in
$\fz(\fg)$ is a non-degenerate abelian ideal. This proves the second assertion.
\qed

\begin{re}\label{georg}
{\rm Slightly different canonical isotropic ideals (and corresponding
notions of balanced extensions) can be obtained
from other canonical decreasing chains of ideals of $\fg$ (like the derived
series of $\fr$ or the lower central series of $\fn$) and the corresponding
increasing chains of their orthogonal complements. In fact, in an earlier
stage of our work
we investigated metric Lie algebras based on
the series
$$ \fg\supset\fr\supset\fr'\supset (\fr')^2\dots\supset (\fr')^{k}\supset\dots
\ .$$
For related constructions
compare also \cite{bordemann} and \cite{N03}. We are indebted to L. Berard
Bergery who drawed our attention to the series used here, which is
also the basis of his investigations of pseudo-Riemanniann holonomy representations
and symmetric spaces (compare \cite{BB1} and \cite{BB2}) and
which seems to be most appropriate for the study of metric Lie algebras.}
\end{re}

Now we are going to derive necessary and sufficient conditions for
a quadratic extension to be (regularly) balanced in terms of the
characterizing data $\fl$, $\rho$, and $[\alpha,\gamma]\in
\cH^2_Q(\fl,\fa)$. We will work with the standard model $\dd$,
$(\alpha,\gamma)\in \cZ^2_Q(\fl,\fa)$, of such an extension.

A first necessary condition is given by Part (b) of Lemma \ref{tau}.

\begin{co}\label{to}
If $\dd$ is balanced, then the representation $\rho$ of $\fl$ on $\fa$ is
semi-simple. In particular, $\rho|_{R(\fl)}=0$.
\end{co}

We look at the decreasing chain of ideals
$$ \fl=R_0(\fl)\supset
R_1(\fl)\supset R_2(\fl)\supset\dots\supset R_k(\fl)\supset\dots
\ .$$
Then we can form corresponding quadratic extensions
$$ \fd_k:=\fd_{\alpha_k,\gamma_k}(\fa,R_k(\fl),\rho_k)\ ,\quad k\ge
0\ ,$$
where $\rho_k$, $\alpha_k$ and $\gamma_k$ are obtained from $\rho$, $\alpha$,
$\gamma$ by restriction
to $R_k(\fl)$ . As a vector space we
have
$$ \fd_k= R_k(\fl)^*\oplus \fa\oplus R_k(\fl)\ .$$
The Lie algebra $\fd_k$ is equipped with a natural action of $\fd=\dd$ by
antisymmetric
derivations.
Indeed, the subspace
$$ \fh_k:=\fl^*\oplus \fa\oplus R_k(\fl)\subset \fd $$
is a coisotropic ideal of $\fd$. Observe that $\fh_k^\perp$ is equal to the annihilator
$R_k(\fl)^\perp=S_k(\fl^*)$ of $R_k(\fl)$ in
$\fl^*$. It follows
$\fd_k \cong \fh_k/\fh_k^\perp$
as $\fd$-module.

The subspace
$$ M_k:= R_k(\fl)^*\oplus\fa\subset \fd_k$$
is a $\fd$-submodule and the projection onto the second summand
$$ pr_\fa: M_k\rightarrow \fa$$
is $\fd$-equivariant.

Now we consider the following conditions (all socles are taken with respect
to the $\fd$-module structure):
\begin{eqnarray*}
(a_k)&& S(\fd_k)\subset M_k\\
(b_k)&& pr_\fa(S(M_k))\subset \fa\ \mbox{ is
non-degenerate w.r.t. } \ip_\fa\\
(b_0^\prime)&& pr_\fa (S(M_0))=0\ .
\end{eqnarray*}

Then we have the following
\begin{lm}\label{loeffel}
The quadratic extension $\dd$ is balanced if and only if the
conditions $(a_k)$ and $(b_k)$ hold for all $k\ge 0$.
The assertion remains true if we replace ``balanced'' by ``regularly
balanced'' and $(b_0)$ by $(b_0^\prime)$.
\end{lm}
\proof
For $k\ge 0$ we introduce the following ideals of $\fd$:
\begin{eqnarray*}
\fri_k&:=&\sum_{l=1}^{k} R_l(\fd)\cap S_l(\fd) =S_k(\fd)\cap \fri(\fd)\ ,\\
\fj_k&:=& \fri_k^\perp= R_k(\fd)+\fj(\fd)=R_k(\fd)+\fm_k\ ,\ \mbox{ where}\\
\fm_k&:=& \sum_{l=1}^{k} R_{l-1}(\fd)\cap S_l(\fd)\supset \fri_k\ .
\end{eqnarray*}
We are interested in the socle and the radical of the $\fd$-module $\fj_k/\fri_k$.
We will frequently use the rules (\ref{basic}).
First we have
$$R(\fm_k)= \sum_{l=1}^{k} R(R_{l-1}(\fd)\cap S_l(\fd))
\subset \sum_{l=1}^{k} R_{l}(\fd)\cap R(S_l(\fd))\subset \fri_k\ .$$
This implies
$$ R(\fj_k/\fri_k)=(R(\fj_k)+\fri_k)/\fri_k= (R_{k+1}(\fd)+R(\fm_k)+\fri_k)/\fri_k
= (R_{k+1}(\fd)+\fri_k)/\fri_k$$
and
\begin{equation}\label{socle1}
S(\fj_k/\fri_k)=R(\fj_k/\fri_k)^\perp=(S_{k+1}(\fd)\cap\fj_k)/\fri_k
=\fm_{k+1}/\fri_k\ .
\end{equation}
It follows that
$$
S(\fj_k/\fri_k)+R(\fj_k/\fri_k)=(R_{k+1}(\fd)+\fm_{k+1})/\fri_k=\fj_{k+1}/\fri_k
$$
and
\begin{equation}\label{socle3}
S(\fj_k/\fri_k)\cap R(\fj_k/\fri_k)=
(S(\fj_k/\fri_k)+R(\fj_k/\fri_k))^\perp=\fri_{k+1}/\fri_k\ .
\end{equation}

For $k\ge 0$ we consider the condition
\begin{eqnarray*}
(c_k) && \fh_k=\fj_k\ .
\end{eqnarray*}
Note that $(c_0)$ is trivially satisfied.
We now claim
\begin{equation}\label{claim}
(a_{k}), (b_{k}), (c_{k}) \Rightarrow (c_{k+1})
\end{equation}

Let us prove (\ref{claim}). We fix $k\ge 0$ and assume
$(a_{k})$, $(b_{k})$, and $(c_{k})$.
Condition $(c_k)$ implies that $\fh_k^\perp=\fri_k$ and $\fd_k=\fh_k/\fh_k^\perp=\fj_k/\fri_k$.
This
together with $(a_k)$ yields
$$S(\fj_k/\fri_k)=S(\fd_k)\subset M_k= (\fl^*\oplus\fa)/\fri_k\ .$$
By (\ref{socle1}) we obtain
\begin{equation}\label{edgar}
\fm_{k+1}\subset \fl^*\oplus\fa\ .
\end{equation}
By $(a_k)$ we have $S(M_k)=S(\fd_k)$. Now $(b_k)$ tells us that
$$\fl^*/\fri_k\supset S(\fd_k)\cap S(\fd_k)^\perp=S(\fd_k)\cap
R(\fd_k)=S(\fj_k/\fri_k)\cap R(\fj_k/\fri_k)\ .$$
By (\ref{socle3}) we obtain
\begin{equation}\label{allan}
\fri_{k+1}\subset \fl^*\ .
\end{equation}
Taking orthogonal complements this gives
\begin{equation}\label{poe}
\fj_{k+1}\supset \fl^*\oplus \fa\ .
\end{equation}
Let $pr_\fl:\fd\rightarrow \fl$ be the natural projection.
Using (\ref{edgar}) we obtain
$$ pr_\fl(\fj_{k+1})=pr_\fl(R_{k+1}(\fd)+\fm_{k+1})=pr_\fl(R_{k+1}(\fd))\ .$$
Now $pr_\fl(R_{k+1}(\fd))=R_{k+1}(\fl)$. This together with (\ref{poe})
shows that
$$\fj_{k+1}=\fl^*\oplus\fa\oplus R_{k+1}(\fl)=\fh_{k+1}\ .$$
This is $(c_{k+1})$, thus we have proved the claim (\ref{claim}).

Now assume that $\fd$ satisfies $(a_k)$ and $(b_k)$ for all $k\in \NN_0$.
Then by (\ref{claim}) and the triviality of $(c_0)$ Condition $(c_k)$ holds for any
$k$. If $k$ is sufficiently large, then $\fj_k=\fj(\fd)$ and
$\fh_k= \fl^*\oplus\fa$. We obtain
$\fl^*\oplus\fa =\fj(\fd)$ or, equivalently, $\fl^*=\fri(\fd)$. Thus $\fd$ is
balanced.

For the opposite direction we first recall that $\fm_{k+1}\subset\fj(\fd)$
and $\fri_{k+1}\subset\fri(\fd)$. Therefore, if $\fd$ is balanced, then
(\ref{edgar}) and (\ref{allan}) hold.
We assume in addition that $(c_k)$ holds. Then
$$ S(\fd_k)=S(\fj_k/\fri_k)=\fm_{k+1}/\fri_k$$
and
$$ S(\fd_k)\cap R(\fd_k)=S(\fj_k/\fri_k)\cap
S(\fj_k/\fri_k)^\perp=\fri_{k+1}/\fri_k \ .$$
Now (\ref{edgar}) implies $(a_k)$, in particular $S(M_k)=S(\fj_k/\fri_k)$.
This together with (\ref{allan}) yields $S(M_k)\cap S(M_k)^\perp\subset \fl^*/\fri_k$, hence
$(b_k)$.
Thus for any $k\ge 0$ the following implication is true
$$ \dd \mbox{ balanced and } (c_k) \Rightarrow (a_k), (b_k) \ .$$
Since $(c_0)$ is the empty condition we obtain using (\ref{claim})
that for a balanced quadratic extension $(a_k)$, $(b_k)$ hold for all $k$.

Assume now that the conditions $(a_k)$, $(b_k)$, and the strengthend
version $(b_0^\prime)$ of $(b_0)$ hold. Then $\fd$ is balanced
and $S(\fd)\subset \fl^*$ which implies that $S(\fd)$ is isotropic,
hence $S(\fd)\subset R(\fd)$. By Assertion (b) of Lemma \ref{ralf}
the extension $\fd$ is regularly balanced.
Vice versa, if $\fd$ is regularly balanced,
then $S(\fd)\subset R(\fd)$. Hence $S(\fd)=S(\fd)\cap R(\fd)\subset
\fri(\fg)=\fl^*$, i.e.,
Condition $(b_0^\prime)$ holds. This finishes the proof of the lemma.
\qed

\begin{re}\label{verdi} \rm{For sufficiently large $k$ condition $(b_k)$ simply
says
that $S(\fa)$ is non-degenerate. Thus $\fa=S(\fa)\oplus S(\fa)^\perp$.
Taking socles we obtain $S\left(S(\fa)^\perp\right)=0$, hence $S(\fa)^\perp=0$.
It follows that $\fa$ is a semi-simple $\fl$-module. We just recover
Corollary
\ref{to}}.
\end{re}

\begin{theo}\label{bebe}
Let $\fl$ be a Lie algebra, let $(\rho,\fa,\ipa)$ be an
orthogonal $\fl$-module, and let $(\alpha,\gamma)\in \cZ^2_Q(\fl,\fa)$.
If $\rho$ is semi-simple, then $\fa=\fa^\fl\oplus \rho(\fl)\fa$, and we have a
corresponding decomposition $\alpha=\alpha_0+\alpha_1$. In this case we consider
the following conditions
\begin{enumerate}
\item[$(A_0)$]
Let $L_0\in \fz(\fl)\cap \ker \rho$ be such that there exist
elements
$A_0\in \fa$ and $Z_0\in \fl^*$ satisfying
for all $L\in\fl$
\begin{enumerate}
\item[(i)] $\alpha(L,L_0)=\rho(L) A_0 $,
\item[(ii)] $\gamma(L,L_0,\cdot)=-\langle A_0,\alpha(L,\cdot)\rangle_\fa +\langle
Z_0, [L,\cdot]_\fl\rangle$ as an element of $\fl^*$,
\end{enumerate}
then $L_0=0$.
\item[$(B_0)$] The subspace $\alpha_0(\ker \lb_\fl)\subset
\fa^\fl$ is non-degenerate.
\item[$(B_0^\prime)$] $\alpha_0(\ker \lb_\fl)=\fa^\fl$.
\item[$(A_k)$] $(k\ge 1)$\\
Let $\fk\subset S(\fl)\cap R_k(\fl)$ be an $\fl$-ideal such that there exist
elements
$\Phi_1\in \Hom(\fk,\fa)$ and $\Phi_2\in \Hom(\fk,R_k(\fl)^*)$ satisfying
for all $L\in\fl$ and $K\in\fk$
\begin{enumerate}
\item[(i)] $\alpha(L,K)=\rho(L)\Phi_1(K)-\Phi_1([L,K]_\fl)$,
\item[(ii)] $\gamma(L,K,\cdot)=-\langle \Phi_1(K),\alpha(L,\cdot)\rangle_\fa +\langle
\Phi_2(K), [L,\cdot]_\fl\rangle +\langle \Phi_2([L,K]_\fl), \cdot \rangle$ as an
element of $R_k(\fl)^*$,
\end{enumerate}
then $\fk=0$.
\item[$(B_k)$] $(k\ge 1)$\\
Let $\fb_k\subset\fa$ be the maximal submodule such that the system of equations
$$ \langle\alpha(L,K), B\rangle_\fa=\langle\rho(L)\Phi(K)-\Phi([L,K]_\fl),B\rangle_\fa\ ,
\quad L\in\fl, K\in R_k(\fl), B\in\fb_k, $$
has a solution $\Phi \in \Hom(R_k(\fl),\fa)$. Then $\fb_k$ is non-degenerate.
\end{enumerate}
Let $m$ be such that $R_{m+1}(\fl)=0$. Then
the quadratic extension $\dd$ is balanced if and only if
$\rho$ is semi-simple and the
conditions $(A_k)$ and $(B_k)$ hold for all $0\le k\le m$.
The assertion remains true if we replace ``balanced'' by ``regularly
balanced'' and $(B_0)$ by $(B_0^\prime)$.
\end{theo}

\proof
For $k\ge m+1$ condition $(a_k)$ is trivially satisfied and $(b_k)$ is equivalent
to the semi-simplicity of $\rho$ (see Remark \ref{verdi}). For $0\le k\le m$ we
will show assuming $\rho$ to be semi-simple that $(A_k)$, $(B_0')$, and $(B_k)$ are
equivalent to $(a_k)$, $(b_0')$, and $(b_k)$, respectively. The theorem then
follows from
Lemma \ref{loeffel}.

First we consider the case $k=0$. Note that for any Lie algebra $\fg$
$$S(\fg)=\fz(\fg)\oplus [\fg, S(\fg)]\ .$$
\begin{lm}
Let $\fg$ be a metric Lie algebra and $\fj\subset \fg$ be a nilpotent ideal.
Then $[\fg, S(\fg)]\subset \fj^\perp$.
\end{lm}
\proof The Lie algebra $\fj$ acts nilpotently and semi-simply, hence trivially, on
$S(\fg)$.
We obtain
$$ \langle  [\fg, S(\fg)], \fj\rangle = \langle \fg, [\fj, S(\fg)]  \rangle = 0\
.$$
\qed

Applying the lemma to the nilpotent ideal $M_0\subset \fd=\dd$ we obtain
$\fz(\fd)\subset S(\fd)\subset \fz(\fd)+\fl^*$. Note that $S(M_0)=S(\fd)\cap M_0$.
Therefore we can reformulate the conditions for $k=0$
in terms of $\fz(\fd)$:
\begin{eqnarray*}
(a_0)&& \fz(\fd)\subset M_0\\
(b_0)&& pr_\fa(\fz(\fd)\cap M_0)\subset \fa\ \mbox{ is
non-degenerate w.r.t. } \ip_\fa\\
(b_0^\prime)&& pr_\fa (\fz(\fd)\cap M_0)=0\ .
\end{eqnarray*}
Using the commutator formulas (\ref{c3}) to (\ref{c8}) we find by straightforward
computation that an element $Z_0-A_0+L_0\in \fl^*\oplus\fa\oplus\fl=\fd$ is central
if and only if the equations
$(i), (ii)$ of condition $(A_0)$ are satisfied and $L_0\in \fz(\fl)\cap
\ker(\rho)$. This shows the equivalence of $(a_0)$ and $(A_0)$. Moreover, it
implies
that $pr_\fa(\fz(\fd)\cap M_0)$ consists of all elements $A_0\in \fa^\fl$ such
that the linear functional on $\Lambda^2(\fl)$ given by $\langle
A_0,\alpha(\cdot,\cdot)\rangle=\langle A_0,\alpha_0(\cdot,\cdot)\rangle$ factors
over the commutator map $\lb_\fl:\Lambda^2(\fl)\rightarrow\fl$, i.e., it is
given by the composition of a linear functional $Z_0$ on $\fl$ with $\lb_\fl$. It
follows
that $pr_\fa(\fz(\fd)\cap M_0)\perp \alpha_0(\ker \lb_\fl)$. Vice versa, if
$A_0\in \alpha_0(\ker \lb_\fl)^{\perp_{\fa^\fl}}$, then we can unambigously
define an element $Z_1\in (\fl')^*$ by
$$ \langle Z_1, [L_1,L_2] \rangle = \langle A_0,\alpha_0(L_1,L_2)\rangle_\fa\ ,
\quad L_1, L_2 \in \fl\ , $$
which can be extended to $\fl$ in an arbitrary manner. We conclude that
$$ \left( pr_\fa(\fz(\fd)\cap M_0)\right)^{\perp_{\fa^\fl}}=\alpha_0(\ker \lb_\fl)\
.$$
This equation implies the equivalence of $(b_0)$ and $(B_0)$ and of $(b_0')$ and
$(B_0')$.

Next we discuss condition $(a_k)$ for $k\ge 1$. We consider the projection
of $\fd$-modules
$p_k: \fd_k\rightarrow R_k(\fl)$. Then $(a_k)$ is equivalent to $p_k(S(\fd_k))=0$.
It turns out to be useful to express this condition in the following awkward way:
Any $\fd$-submodule $\fk\subset p_k(S(\fd_k))$ vanishes.

We have $p_k(S(\fd_k))\subset S(\fl)\cap R_k(\fl)$.
We
claim that a $\fd$-submodule $\fk\subset S(\fl)\cap R_k(\fl)$ is contained
in $p_k(S(\fd_k))$ if and only if there exists an element $\Phi\in
\Hom_\fd(\fk,\fd_k)$
such that $p_k\circ \Phi=\id$. Indeed, if $\fk\subset p_k(S(\fd_k))$, then
we can choose a $\fd$-invariant complement $\tilde \fk$ of $p_k^{-1}(0)\cap
S(\fd_k)$ in the semi-simple $\fd$-module $p_k^{-1}(\fk)\cap S(\fd_k)$. The projection $p_{k}$
maps $\tilde\fk$  isomorphically to $\fk$. Then we can take
$\Phi=(p_k|_{\tilde\fk})^{-1}$.
Vice versa, if
$\Phi$ as above exists, then by semi-simplicity of $\fk$ the module $\tilde \fk:=\Phi(\fk)$ is
semi-simple and
therefore $p_k(S(\fd_k))\supset p_k(\tilde \fk)=\fk$.

Let $i: \fk \rightarrow R_k(\fl)$ be the natural inclusion.
Observe that $[R(\fl),S(\fl)]=0$ implies $\fk\subset R_k(\fl)^{R(\fl)}\subset
\fz(R_k(\fl))$.
Using this fact and again the formulas (\ref{c3}) to (\ref{c8})
one shows that a homomorphism
$$ \Phi = (\Phi_2 , -\Phi_1 ,i ) \in \Hom(\fk,R_k(\fl)^*)\oplus \Hom(\fk,\fa)
\oplus \Hom(\fk, R_k(\fl))= \Hom(\fk,\fd_k) $$
is $\fd$-equivariant if and only if the equations
$(i), (ii)$ of condition $(A_k)$ are satisfied and $\im \Phi_1\in\fa^{R_k(\fl)}$.
The latter condition is vacuous since for $k\ge 1$ we have by semi-simplicity that
$\rho|_{R_k(\fl)}=0$. We conclude that $(a_k)$ is equivalent to $(A_k)$.

Concerning $(b_k)$, $k\ge 1$, the same reasoning as above yields
that a $\fd$-submodule $\fb\subset \fa$ is contained
in $pr_\fa(S(M_k))$ if and only if there exists an element $\Psi\in
\Hom_\fd(\fb,M_k)$
such that $pr_\fa\circ \Psi=\id$. An element $\Psi=(\Psi_1,i)\in
\Hom(\fb,R_k(\fl)^*)\oplus \Hom(\fb,\fa )=\Hom(\fb,M_k)$ is $\fd$-equivariant if
and only if
for all $L\in\fl$ and $B\in\fb$
$$ \Psi_1(\rho(L) B)=-\langle B, \alpha(L,\cdot)\rangle_\fa -\langle \Psi_1(B),
[L,\cdot]_\fl\rangle \in R_k(\fl)^*\ .$$
If $\Phi=(\Psi_1\circ i)^*\in \Hom(R_k(\fl),\fa)$, then this equation is equivalent
to
$$ \langle \alpha (L,K), B\rangle_\fa = \langle \rho(L)\Phi(K) -\Phi([L,K]_\fl),
B\rangle_\fa\ , \quad L\in\fl, K\in R_k(\fl), B\in\fb.$$
It follows that $(b_k)$ is equivalent to $(B_k)$ for all $k\ge 1$.
\qed

\begin{re}
{\rm Theorem \ref{bebe} shows in particular that in the case of abelian $\fl$ which
has been studied in \cite{KO03}
the extension $\dd$ is regularly balanced if and
only if it is regular in the sense of \cite{KO03} and the representation $\rho$ is
semi-simple.}
\end{re}

Each of the conditions $(a_k),(b_k)$, and $(b_0')$ only depends on the equivalence
class of the quadratic extension $\dd$. By Proposition \ref{P} this implies that
each
of the conditions $(A_k),(B_k)$, and $(B_0')$ only depends on the cohomology
class $[\alpha,\gamma]\in \cH^2_Q(\fl,\fa)$. It is an interesting exercise
to check this directly.

\begin{de}\label{zwiesel}
    Let $\fa$ be a semi-simple orthogonal $\fl$-module. Let $m\in\NN_0$ be such
that $R_{m+1}(\fl)=0$.
    A cohomology class $[\alpha,\gamma]\in \cH^2_Q(\fl,\fa)$ is called admissible
    if $(\alpha,\gamma)\in \cZ_Q^2(\fl,\fa)$ satisfies the conditions $(A_k)$,
$(B_k)$ for all $0\le k\le m$.
    We denote the subset of admissible classes in
    $\cH^2_Q(\fl,\fa)$ by $\cH^2_{Q}(\fl,\fa)_\sharp$. If $\fa$ is not semi-simple,
then we
    set $\cH^2_Q(\fl,\fa)_\sharp=\emptyset$. Furthermore, a Lie algebra $\fl$ is
called admissible if
    there exists a semi-simple orthogonal $\fl$-module $\fa$ such that
    $\cH^2_{Q}(\fl,\fa)_\sharp\not=\emptyset$.
\end{de}

Now we can reformulate Theorem \ref{bebe} in the following way:

\begin{co}\label{kohl}
A quadratic extension is balanced if and only if the quadratic
cohomology class assigned to it by (\ref{Ealpha}) and (\ref{Egamma}) belongs to
$\cH^2_Q(\fl,\fa)_\sharp$.
\end{co}

%%%%%%%%%%%%%%%%%%%%%%%%%%%%%%%%%%%%%%%%%%%%%%%%%%%%%%%%%%%%%%%%%%%%%%%%%%%%%%%%%%%
%%%%%%%%%%%%%%%%%%%%%%%%%%%%%%%%%%%%%%%%%%%%%%%%%%%%%%%%%%%%%%%%%%%%%%%%%%%%%

\section{Isomorphy and decomposability of metric Lie algebras}\label{isodec}

By the results of the previous sections (in particular Proposition \ref{schmidt},
Theorem \ref{brandt}, and Corollary \ref{kohl}) the metric Lie algebras $\dd$
associated
with semi-simple orthogonal modules $(\rho,\fa)$ of a Lie algebra $\fl$ and
$[\alpha,\gamma]\in \cH^2_Q(\fl,\fa)_\sharp$ exhaust all isomorphism classes
of metric Lie algebras without simple ideals. In order to approach
the classification of indecomposable metric Lie algebras we have to decide
which of these data lead to isomorphic or decomposable metric Lie algebras,
respectively. This is the first aim of the present section. We conclude the section
giving a classification scheme for isomorphism classes
of non-simple indecomposable metric Lie algebras.

Recall the definition of morphisms $(S,U): (\fl_1,\fa_1)\rightarrow (\fl_2,\fa_2)$
of pairs (of Lie algebras and orthogonal modules) from Section \ref{qc} around Equation
(\ref{opair}). If $(S,U)$ is an isomorphism of pairs and $(\fg,\fri,i,p)$
is a quadratic extension of
$\fl_1$ by $\fa_1$, then we observe that $(\fg,\fri,i\circ U,S\circ p)$ is a quadratic
extension of $\fl_2$ by $\fa_2$.

\begin{lm}\label{pieck}
Let $(\rho_j,\fa_j)$, $j=1,2$, be orthogonal modules of Lie algebras $\fl_j$,
and let $(\fg_j,\fri_j,i_j,p_j)$ be quadratic extensions of $\fl_j$ by $\fa_j$.
\begin{enumerate}
\item[{\rm (a)}] If there is an isomorphism of pairs
$(S,U): (\fl_1,\fa_1)\rightarrow (\fl_2,\fa_2)$ such that
the quadratic extensions $(\fg_1 ,\fri_1,i_1\circ
U,S\circ p_1)$ and $(\fg_2,\fri_2,i_2,p_2)$ are equivalent,
then $\fg_1$ and
$\fg_2$ are isomorphic as metric Lie
algebras.
\item[{\rm (b)}] If the quadratic extensions $(\fg_j,\fri_j,i_j,p_j)$, $i=1,2$, are
balanced, and the metric Lie algebras $\fg_1$ and $\fg_2$ are isomorphic, then
there
exists an isomorphism of pairs $(S,U): (\fl_1,\fa_1)\rightarrow
(\fl_2,\fa_2)$ such that
the extensions $(\fg_1 ,\fri_1,i_1\circ
U,S\circ p_1)$ and $(\fg_2,\fri_2,i_2,p_2)$ are equivalent.
\end{enumerate}
\end{lm}
\proof Part (a) is a triviality since any equivalence of quadratic extensions is by
definition
an isomorphism of metric Lie algebras.
Let us prove (b). Let $F:\fg_1\rightarrow
\fg_2$ be an isomorphism of metric Lie algebras.
Then $F(\fri(\fg_1))=\fri(\fg_2)$ and $F(\fj(\fg_1))=\fj(\fg_2)$.
The quadratic extensions are balanced, i.e.,
$\fri_j=\fri(\fg_j)$. Thus $F$ induces a Lie algebra isomorphism
$$ \bar F: \fg_1/\fri_1 \longrightarrow \fg_2/\fri_2 $$
such that $\bar F(\ker p_1)=\ker p_2$. We then define $(S,U)$ by
$$ S(p_1(X))=p_2(\bar F(X))\ ,\ X \in \fg_1/\fri_1,\quad
i_1(U(A))=\bar F^{-1}(i_2(A))\ ,\ A\in\fa_2\ .$$
We compute for $L\in \fl_1$, and $\tilde L\in \fg_1/\fri$ such that $p_1(\tilde L)=L$
and $A\in\fa_2$
\begin{eqnarray}\label{stalin}
 i_1\circ\rho_1(L)\circ U(A)&=&[\tilde L,i_1\circ U(A)]=
\bar F^{-1}[\bar F (\tilde L),i_2(A)]\\
&=& \bar F^{-1}\circ i_2\circ\rho_2(S(L))(A)=i_1\circ U \circ\rho_2(S(L))(A)\ .\nonumber
\end{eqnarray}
In the third step we have used that $p_2(\bar F(\tilde L))=S(L)$. Equation
(\ref{stalin}) shows that $(S,U)$ is an isomorphism of pairs. Now $F$ defines
an equivalence between $(\fg_1,\fri_1,i_1\circ U,S\circ p_1)$ and $(\fg_2 ,\fri_2,i_2,p_2)$.
\qed

We can now give a necessary and sufficient criterion of isomorphy of metric Lie algebras
in terms of the admissible quadratic cohomology classes associated with their structures as a
balanced
quadratic extensions (see Proposition \ref{schmidt}).

\begin{pr}\label{wiesel}
Let $(\rho_i,\fa_i)$, $i=1,2$, be orthogonal modules of Lie algebras $\fl_i$,
and let $(\alpha_i,\gamma_i)\in \cZ^2_Q(\fl_i,\fa_i)$.
\begin{enumerate}
\item[{\rm (a)}] If there exists an isomorphism of pairs
$(S,U): (\fl_1,\fa_1)\rightarrow (\fl_2,\fa_2)$ such that \linebreak[4]
$ (S,U)^*[\alpha_2,\gamma_2]=[\alpha_1,\gamma_1]\in \cH^2_Q(\fl_1,\fa_1)$,
then $\fd_{\alpha_1,\gamma_1}(\fl_1,\fa_1,\rho_1)$ and
$\fd_{\alpha_2,\gamma_2}(\fl_2,\fa_2,\rho_2)$ are isomorphic as metric Lie
algebras.
\item[{\rm (b)}] If $[\alpha_i,\gamma_i]\in \cH^2_Q(\fl_i,\fa_i)_\sharp$, $i=1,2$,
and
$\fd_{\alpha_1,\gamma_1}(\fl_1,\fa_1,\rho_1)$ and
$\fd_{\alpha_2,\gamma_2}(\fl_2,\fa_2,\rho_2)$ are isomorphic metric Lie algebras,
then there
is an isomorphism of pairs $(S,U): (\fl_1,\fa_1)\rightarrow
(\fl_2,\fa_2)$ such that
$ (S,U)^*[\alpha_2,\gamma_2]=[\alpha_1,\gamma_1]\in
\cH^2_Q(\fl_1,\fa_1)_\sharp$.
\end{enumerate}
\end{pr}
\proof We write the quadratic extensions
$\fd_{\alpha_j,\gamma_j}(\fl_j,\fa_j,\rho_j)$ as $(\fd_j, \fl_j^*,i_j,p_j)$.
Let $s: \fl_1\rightarrow \fd_1$ be
the embedding. Let $(S,U): (\fl_1,\fa_1) \rightarrow (\fl_2,\fa_2)$ be an isomorphism of pairs.
Then $\tilde s:=s\circ S^{-1}: \fl_2\rightarrow \fd_1$
is a section of $S\circ \tilde p_1$ with isotropic image. Now the cocycle associated
with the quadratic
extension $(\fd_1, \fl_1^*,i_1\circ U,S\circ p_1)$ of $\fl_2$ by $\fa_2$
and the section $\tilde s$ is given by
$$ \left( (S^{-1}, U^{-1})^*\alpha_1, (S^{-1})^*\gamma_1\right)$$
(see (\ref{Ealpha}) and (\ref{Egamma})).
Therefore the proposition is a consequence of Lemma \ref{pieck}, Proposition
\ref{P}, and Corollary \ref{kohl}.
\qed

There is a natural notion of a direct sum of quadratic extensions.
Namely, if
$(\fg_j,\fri_j,i_j,p_j)$, $j=1,2$, are quadratic extensions of Lie algebras $\fl_j$
by orthogonal modules $\fa_j$, then
$$ (\fg_1\oplus\fg_2,\fri_1\oplus\fri_2,i_1\oplus i_2,p_1\oplus p_2)$$
is a quadratic extension of $\fl_1\oplus\fl_2$ by $\fa_1\oplus\fa_2$.
A quadratic extension is called decomposable, if it can be written
as a non-trivial direct sum of two quadratic extensions.
If a quadratic extension is equivalent to a decomposable one, then it is
decomposable.
Of course, the decomposability of a quadratic extension $(\fg,\fri,i,p)$ implies
the decomposability of $\fg$ as a metric Lie algebra.
The opposite assertion is not true in general. However, we have

\begin{lm}\label{basel}
Let $(\fg,\fri,i,p)$ be a balanced quadratic extension of $\fl$ by $\fa$.
If the metric Lie algebra $\fg$ is decomposable, then the quadratic extension
$(\fg,\fri,i,p)$ is decomposable, too.
\end{lm}
\proof  If $\fg=\fg_1\oplus\fg_2$ and $(\fg,\fri,i,p)$ is balanced, then
$\fri=\fri(\fg)=\fri(\fg_1)\oplus\fri(\fg_2)$. In particular,
$\fg/\fri=\fg_1/\fri(\fg_1)\oplus\fg_2/\fri(\fg_2)$.
We set $\fl_j:=p(\fg_j/\fri(\fg_j))$, $p_j:=  p|_{\fg_j/\fri(\fg_j)}$,
$\fa_j=i^{-1}(\fg_j/\fri(\fg_j))$, $i_j:=i|_{\fa_j}$.
Then $(\fg_j,\fri_j,i_j,p_j)$ are quadratic extensions of $\fl_j$ by $\fa_j$, and
$$ (\fg,\fri,i,p)=
(\fg_1\oplus\fg_2,\fri(\fg_1)\oplus\fri(\fg_2),i_1\oplus i_2,p_1\oplus p_2)\ .$$
\qed

Recall the notion of a non-trivial direct sum of pairs from Defininition \ref{pro}.

\begin{de}\label{kosel}
Let $(\rho,\fa)$ be an orthogonal module of a Lie algebra $\fl$.
A cohomology class $c\in \cH^2_Q(\fl,\fa)$ is called decomposable if there is a decomposition
$$(\fl,\fa)=(\fl_1,\fa_1)\oplus (\fl_2,\fa_2)$$
into a non-trivial direct
sum of pairs such that in the notation of Lemma \ref{pru}
$$ c\in (q_1,j_1)^*\cH^2_Q(\fl_1,\fa_1)+ (q_2,j_2)^*\cH^2_Q(\fl_2,\fa_2)\ .$$
A cohomology class which is not decomposable is called indecomposable.
We denote the set of all indecomposable elements
in $\cH^2_Q(\fl,\fa)_\sharp$ by $\cH^2_Q(\fl,\fa)_0$.
\end{de}

\begin{pr}\label{diesel}
Let $(\rho,\fa)$ be an orthogonal module of a Lie algebra $\fl$,
and let $(\alpha,\gamma)\in \cZ^2_Q(\fl,\fa)$.
\begin{enumerate}
\item[{\rm (a)}] If $[\alpha,\gamma]\in \cH^2_Q(\fl,\fa)$ is decomposable,
then the quadratic extension $\dd$ is decomposable.
\item[{\rm (b)}] If $[\alpha,\gamma]\in \cH^2_Q(\fl_i,\fa_i)_\sharp$ and
$\dd$ is decomposable as a metric Lie algebra, then $[\alpha,\gamma]$
is decomposable.
\end{enumerate}
\end{pr}
\proof
Let $\fl_i$, $\fa_i$, $j_i$, and $q_i$ be as in Definition \ref{kosel}.
We first note that, if a quadratic extension $(\fg,\fri,i,p)$ is the direct sum
$$ (\fg_1,\fri_1,i_1,p_1)\oplus (\fg_2,\fri_2,i_2,p_2) $$
of quadratic extensions of $\fl_i$ by $\fa_i$ with associated cohomology classes
$c_i\in \cH^2_Q(\fl_i,\fa_i)$, then the cohomology class associated with
$(\fg,\fri,i,p)$ is given by
$(q_1,j_1)^*c_1+ (q_2,j_2)^*c_2\in \cH^2_Q(\fl,\fa)$.

Assume now that $[\alpha,\gamma]\in \cH^2_Q(\fl,\fa)$ is decomposable. Then
there exist elements $(\alpha_i,\gamma_i)\in \cZ^2_Q(\fl_i,\fa_i)$ such that
$[\alpha,\gamma]=(q_1,j_1)^*[\alpha_1,\gamma_1]+ (q_2,j_2)^*[\alpha_2,\gamma_2]$.
By the above and \linebreak[4] Theorem \ref{brandt} the quadratic extension $\dd$ is equivalent
to
the direct sum \linebreak[4]
$\fd_{\alpha_1,\gamma_1}(\fl_1,\fa_1,\rho_1)\oplus
\fd_{\alpha_2,\gamma_2}(\fl_2,\fa_2,\rho_2)$ and therefore decomposable.
This proves (a).

If $[\alpha,\gamma]\in \cH^2_Q(\fl_i,\fa_i)_\sharp$ and
$\dd$ is decomposable as a metric Lie algebra, then $\dd$ is balanced and thus
by Lemma \ref{basel} decomposable as a quadratic extension. Then the discussion
at the beginning of the proof shows that $[\alpha,\gamma]$ is decomposable.
\qed

We conclude this section with a classification scheme for isomorphism classes
of non-simple indecomposable metric Lie algebras.

Let us fix a Lie algebra $\fl$ and a semi-Euclidean vector space $(\fa,\ipa)$.
We consider the set $\Hom(\fl,\so(\fa,\ip_\fa))_{\rm ss}$ of all orthogonal
semi-simple representations
of $\fl$ on $\fa$. If $\rho\in \Hom(\fl,\so(\fa,\ip_\fa))_{\rm ss}$ is fixed
we denote the corresponding $\fl$-module by $\fa_\rho$.
The group $G:=\Aut(\fl)\times O(\fa,\ip_\fa)$ acts from the right on
$\Hom(\fl,\so(\fa,\ip_\fa))_{\rm ss}$
by
$$(S,U)^*\rho:=\Ad(U^{-1})\circ S^*\rho\ ,\qquad S\in \Aut(\fl),\  U\in
O(\fa,\ip_\fa) \ .$$
Then for any $\rho\in \Hom(\fl,\so(\fa,\ip_\fa))_{\rm ss}$ an element
$g=(S,U)\in G$ defines an isomorphism of pairs
$ \bar g:=(S,U^{-1}): (\fl,\fa_{g^*\rho})\rightarrow (\fl,\fa_\rho)$
and therefore induces a bijection
$$ \bar g^*: \cH^2_Q(\fl,\fa_\rho)\rightarrow \cH^2_Q(\fl,\fa_{g^*\rho})\
.$$
We
obtain a right action of $G$ on the disjoint union
$$ \coprod_{\rho\in \Hom(\fl,\so(\fa,\ipa))_{\rm ss}} \cH^2_Q(\fl,\fa_\rho)\ .$$
If it is clear from the context that a pair
$g=(S,U)$ is considered as an element of $G$, then its action on the
above space will be simply denoted by $g^*$ or $(S,U)^*$. Note the
slight difference of the meaning of $(S,U)^*$, if $(S,U)$ is a morphism of pairs.

As in Definition \ref{kosel} let $\cH^2_{Q}(\fl,\fa_\rho)_{0}\subset
\cH^2_Q(\fl,\fa_\rho)$ be the subset
of all admissible indecomposable elements (see Definition \ref{zwiesel} and Theorem
\ref{bebe} for the admissibility conditions). Then the set
$$ \coprod_{\rho\in \Hom(\fl,\so(\fa,\ipa))_{\rm ss}} \cH^2_Q(\fl,\fa_\rho)_{0}$$
is $G$-invariant.
Combining Proposition \ref{schmidt} and Corollary \ref{kohl} with Propositions
\ref{wiesel} and \ref{diesel} we obtain
\begin{theo}\label{Pwieder}
    Let $\fl$ be a Lie algebra, and let $(\fa,\ipa)$ be a semi-Euclidean vector
space. We consider the
    class $\cA(\fl,\fa)$ of non-simple indecomposable metric Lie algebras $\fg$
                       satisfying
    \begin{enumerate}
        \item The Lie algebras $\fg/\fj(\fg)$ and $\fl$ are isomorphic.
        \item $\fj(\fg)/\fri(\fg)$ is isomorphic to $(\fa,\ipa)$ as a
semi-Euclidean vector space.
    \end{enumerate}
    Then the set of isomorphism classes of $\cA(\fl,\fa)$ is in bijective
    correspondence with the orbit space of the action of $G=\Aut(\fl)\times
    O(\fa,\ip_\fa)$ on
    $$\coprod_{\rho\in \Hom(\fl,\so(\fa,\ipa))_{\rm ss}} \cH^2_Q(\fl,\fa_\rho)_{0}\
.$$
This orbit space can also be written as
$$\coprod_{[\rho]\in \Hom(\fl,\so(\fa,\ipa))_{\rm ss}/G}
\cH^2_Q(\fl,\fa_\rho)_{0}/G_\rho\ ,$$
where $G_\rho=\{g\in G\:|\: g^*\rho=\rho\}$ is the automorphism group of the pair
$(\fl,\fa_\rho)$.
\end{theo}

In view of Proposition \ref{inner} the automorphism group $G_\rho$ can be replaced by the group
of outer automorphisms of $(\fl,\fa_\rho)$.

Note that $\cA(\fl,\fa)$ is empty if the Lie algebra $\fl$ is not
admissible (see Definition \ref{zwiesel}).
We will use Theorem \ref{Pwieder} in order to provide
a classification of all indecomposable metric Lie algebras of index
3 (see Theorem \ref{T71}). In the way of this classification, in particular in
Section \ref{SM},
we shall see that there are  many isomorphism classes of
non-admissible Lie algebras $\fl$. In order to apply Theorem \ref{Pwieder} to more
general situations one would like to have (as a first step) a good description
of the class $\cM$ of all admissible Lie algebras. Up to
now we only know that $\cM$ contains all reductive Lie algebras as well
as some solvable, nilpotent and mixed ones, and that $\cM$ is closed under
forming direct sums.

%%%%%%%%%%%%%%%%%%%%%%%%%%%%%%%%%%%%%%%%%%%%%%%%%%%%%%%%%%%%%%%%%%%%%%%%%%%%%%%%%%%
%%%%%%%%%%%%%%%%%%%%%%%%%%%%%%%%%%%%%%%%%%%%%%%%%%%%%%%%%%%%%%%%%%%%%%%%%%%%%
\section{Solvable admissible Lie algebras with small radical}
\label{SM}
In this section we will classify all solvable admissible Lie algebras
whose radical of nilpotency is one- or two-dimensional. Recall that a Lie algebra $\fl$ is called admissible if there exists a balanced
quadratic extension $\dd$. On the one hand this
classification serves as an example which shows how one can handle the
admissibility conditions. On the other hand we will need the results in Section
\ref{Sindex3} in order to classify the metric Lie algebras of index 3.

\subsection{Weight spaces}
Throughout this section we will assume that $\fl$ is a solvable Lie algebra whose radical of nilpotency $R:=R(\fl)=\fl'$ is abelian. Moreover, let $(\rho,\fa)$ be a semi-simple orthogonal
$\fl$-module.
Since $R$ is abelian
$$\adR: \fl\longrightarrow \gl(R),\quad L\longmapsto \ad(L)|_{R}$$
induces a representation of the abelian Lie algebra $\fl/R$ on $R$, which we also
denote by $\adR$. Furthermore, since $\rho$ is semi-simple we have $\rho|_R=0$ and
$\fa$ can be considered as a semi-simple $\fl/R$-module.
Hence the complexification
$\fa_ {\Bbb C}$ of $\fa$ decomposes into $\fa_{\Bbb C}=
\bigoplus_{\lambda\in\Lambda}E_\lambda$, where
$$\Lambda=(\fl_{\Bbb C}/R_{\Bbb C})^*\cong
\{\lambda\in\fl_{\Bbb C}^*\mid \lambda|_{R_{\Bbb C}}=0\} \subset \fl_{\Bbb
C}^*$$
and
$$E_{\lambda}=\{A\in\fa_{\Bbb C}\mid \rho(L)(A)=\lambda(L)\cdot A \mbox{
for all } L\in \fl\}.$$
Let $p_{\lambda}:\fa_{\Bbb C}\rightarrow E_{\lambda}$ denote the
projection.

Let $\ip$ on $\fa_{\Bbb C}$ be the sesquilinear extension of $\ip$
on $\fa_ {\Bbb C}$. Then $E_{\lambda}\perp E_{\mu}$ if not $\mu=-\bar\lambda$.
Let us also define analogous spaces for the complexification $R_{\Bbb
C}$ of $\adR$ by
$$V_{\lambda}:=\{U\in R_{\Bbb C}\mid \adR(L)(U)=\lambda(L)\cdot U \mbox{
for all } L\in \fl\}.$$
\begin{lm}\label{fliege}
    Let $R$ be abelian and $\rho\in\Hom(\fl,\so(\fa,\ip_{\fa}))_{\rm ss}$. Assume
$\alpha\in Z^{2}(\fl,\fa)$ satisfies
    $\alpha|_{R\times R}=0$. Then there exists a cocycle $\tilde\alpha\in
    Z^{2}(\fl,\fa)$ such that $\tilde\alpha|_{R\times R}=0$,
$[\alpha]=[\tilde\alpha]\in
    H^{2}(\fl,\fa)$ and
    $ \tilde\alpha(\fl,V_\lambda)\subset E_{\lambda}$
    holds for all $\lambda\in \Lambda$.
    \end{lm}
\proof
\ Because of $\alpha|_{R\times R}=0$ and $\rho|_{R}=0$ the cocycle $\alpha$ defines
a
cocycle $\bar\alpha\in Z^{1}(\fl/R,C^1(R,\fa))$ by
$\bar\alpha(L+R)(U)=\alpha(L,U) \mbox{ for } L\in\fl,\ U\in R.$ Since
$\fl/R$ is abelian  we have $H^{1}(\fl/R,C^1(R,\fa))$ $=
H^{1}(\fl/R,C^1(R,\fa)^{(\fl/R)})$. Here $C^1(R,\fa)^{(\fl/R)}$ denotes the nilsubspace of $C^1(R,\fa)$ with respect to the action of $\fl/R$. Thus there exists a cochain $\bar
\tau \in C^{0}(\fl/R,  C^1(R,\fa))\cong  C^1(R,\fa)$ such that
$\bar\alpha+d\bar\tau$ has values only in $C^1(R,\fa)^{(\fl/R)}$.
Hence for all $L'\in\fl$ there exists a  $k\in\NN$ such that
$(L')^{k}\cdot((\bar\alpha+d\bar\tau)(L+R))=0$. For $U\in V_{\lambda}$ we
have
$$0=(\,(L')^{k}\cdot((\bar\alpha+d\bar\tau)(L+R))\,)(U)=(\rho(L')-
\lambda(L'))^{k}(\,(\bar\alpha+d\bar\tau)(L+R)(U)\,).$$ Since $\rho$
is semi-simple we get $(\bar\alpha+d\bar\tau)(L+R)(U)\in E_{\lambda}$
for all $L\in\fl$.
Now we can choose a linear extension $\tau\in C^{1}(\fl,\fa)$ of
$\bar\tau$ and obtain
$(\alpha+d\tau)(L,U)=(\bar\alpha+d\bar\tau)(L+R)(U)\in E_{\lambda}$.
\qed

\subsection{$\dim R(\fl)=1$}

In the following we will often describe a Lie algebra $\fl$ giving
only the non-trivial Lie bracket relations between vectors of a basis.
We will use this notation only for those Lie algebras for which all basis vectors appear in one of these relations.

In Section \ref{Sindex3} we will see that the Heisenberg algebra
$$\fh(1)=\{[X,Y]=Z\}$$
is admissible. On the other hand we can prove the following.
\begin{pr}\label{R1}
    \ If $\fl$ is a solvable admissible Lie algebra with $\dim R(\fl)=1$,
    then $$\fl\cong\fh(1)\oplus \RR^{k}.$$
\end{pr}
\proof
We choose $Z\in\fl$ such that $R=\RR\cdot Z$. Then we define
$\lambda\in\Lambda$ by
$[L,Z]=\lambda(L)Z$. Obviously, $\lambda$ is real. By Lemma  \ref{fliege} we may assume
$\alpha(Z,\fl)\subset
E_\lambda$.
Suppose $\lambda\not=0$. Then $E_{\lambda}$ is isotropic and therefore
$E_\lambda\subset \fb_1$ (see Theorem \ref{bebe} $(B_1)$,
$\Phi=0$ is a solution).
By $(B_1)$ the space $\fb_1$ is non-degenerate. Therefore also $E_{-\lambda}\subset
\fb_1$.
By definition of $\fb_1$ there exists  an element $\Phi\in\Hom(R,\fa)$ such that
$$\langle\alpha(L,Z),B\rangle=\langle \rho(L)(\Phi(Z))-\Phi([L,Z]),B\rangle=
\langle ( \rho(L)-\lambda(L))(\Phi(Z)),B\rangle=0$$
for all $B\in E_{-\lambda}$. Hence $\alpha(L,Z)=0$ for all $L\in\fl$.
Condition $(A_1)$ now implies $R=0$ which is a contradiction. Thus $\lambda=0$.
This implies $$\fl\cong\{ [X_{2i-1},X_{2i}] =Z\mid i=1,...,r\}\oplus\RR^k.$$
For $r>1$ the cocycle condition on $\alpha$ yields $\alpha(Z,\cdot)=0$. Hence 
$r=1$
and the assertion follows.
\qed

\subsection{$\dim R(\fl)=2$}
We consider the following Lie algebras with two-dimensional radical 
of nilpotency:
\begin{eqnarray*}
    &&\fn(2)=\{ [X,Y]=Z,\ [X,Z]=-Y\},\\
    &&\fr_{3,-1}=\{ [X,Y]=Y,\ [X,Z]=-Z\},\\
    &&\fr_{3,-2}=\{ [X,Y]=-2Y,\ [X,Z]=Z\}.
\end{eqnarray*}
All these Lie algebras are admissible. For $\fn(2)$ and $\fr_{3,-1}$ 
we will see this in Section \ref{Sindex3}. As for $\fr_{3,-2}$ consider 
$\fa=\RR^{1,1}$ spanned by the isotropic vectors $e_+$ and $e_-$. Let $\rho$ be 
given by
$$\rho(X)(e_+)=e_+,\  \rho(X)(e_-)=-e_-,\ \rho(Y)=\rho(Z)=0.$$
We define $(\alpha,\gamma)\in\cZ_Q^2(\fl,\fa_\rho)$ by
$$\alpha(Y,Z)=e_-,\ \alpha(X,Z)=e_+,\ \alpha(X,Y)=0,\ \gamma=0.$$
Then $[\alpha,0]\in\cH^2_Q(\fl,\fa_\rho)_\sharp$, hence $\fr_{3,-2}$ is admissible.

In this section we will prove the following classification result.
\begin{pr}\label{elefant}
\ If $\fl$ is a solvable, non-nilpotent admissible Lie algebra with 
\linebreak $\dim R(\fl)=2$,
then $\fl$ is isomorphic to one of the
Lie algebras
$$\fn(2)\oplus \RR^{k},\ \fr_{3,-1}\oplus \RR^{k},\ \fr_{3,-2}\oplus
\RR^{k}.$$
\end{pr}

In the following let $\fl$ be a solvable, non-nilpotent Lie algebra
with $\dim R=2$.  Since in this case $R$ is (two-dimensional and) nilpotent it must 
be abelian.

Let us first assume that the representation $\adR$
of $\fl_{\Bbb C}/R_{\Bbb C}$ on $R_{\Bbb C}$ is semi-simple. Then $\adR$ 
has two weights $\lambda^{1},\ \lambda^{2}$ and either $\lambda^1$ 
and  $\lambda^2$ are real weights or
$\lambda^1=\overline{\lambda^2}$ are complex weights. In both cases we may assume 
$\lambda^1\not=0$ since $\fl$ is not nilpotent.
\begin{lm}\label{muecke}
    Any cocycle $\alpha\in Z^{2}(\fl,\fa)$ satisfies
    $\alpha(R,R)\subset E_{\lambda^1+\lambda^2}$.
    There exists a cocycle $\tilde\alpha\in
    Z^{2}(\fl,\fa)$ such that $[\alpha]=[\tilde\alpha]\in
    H^{2}(\fl,\fa)$ and $\tilde \alpha(R,\fl)\subset 
    E_{\lambda^{1}+\lambda^{2}}+E_{\lambda^{1}}+E_{\lambda^{2}}$.
\end{lm}
\proof The first assertion follows from the cocycle condition. To 
prove the second one we look at the decomposition 
$\alpha=p_{\lambda^1+\lambda^{2}} \circ \alpha + \alpha'$. Then the 
first assertion implies 
$\alpha'|_{R\times R}=0$. Hence we may apply  Lemma \ref{fliege} to 
$\alpha'$ and 
the second assertion follows. 
\qed

\begin{lm} \label{Lk}  Let $\fl$ be admissible and $[\alpha,\gamma]\in
\cH^2_{Q}(\fl,\fa)_\sharp$.  Assume $\lambda^{1}+ \lambda^{2}\not=0$. Let 
$\fk\subset 
    R(\fl)\cap S(\fl)=R(\fl)$ be an $\fl$-ideal. If $\alpha(K,L)=0$ for all
    $K\in\fk$ and $L\in\fl$, then $\fk=0$.
\end{lm}
\proof
If $\alpha(K,L)=0$ for all $K\in\fk$ and $L\in\fl$, then
$(A_{1})\,(i)$ is satisfied for $\Phi_1=0$. We will show that $(A_{1})\,(ii)$ is 
also
satisfied. Suppose $L_{1},L_{2}\in\fl$, $K\in\fk$ and $U\in R$.
Because of $\alpha(K,\cdot)=0$ we have
$d\gamma(K,U,L_{1},L_{2})=\frac 12 \langle \alpha \wedge
\alpha\rangle (K,U,L_{1},L_{2})=0$. On the other hand
$ d\gamma (K,U,L_{1},L_{2}) = -(\lambda^{1}+
\lambda^{2})(L_{1})\gamma(K,U,L_{2}) + (\lambda^{1}+
\lambda^{2})(L_{2})\gamma(K,U,L_{1})$ holds where we used $\dim R=2$. Hence
\begin{equation}\label{gerade}
(\lambda^{1}+ \lambda^{2})(L_{1})\gamma(K,U,L_{2}) = (\lambda^{1}+
\lambda^{2})(L_{2})\gamma(K,U,L_{1})\,.
\end{equation}
Now let $L_{0}\in\fl$ be such
that $(\lambda^{1}+ \lambda^{2})(L_{0})\not=0$ and define
$\Phi_{2}\in\Hom(\fk,R(\fl)^{*})$ by
$$\langle \Phi_{2}(K),U\rangle = \textstyle{\frac 1{(\lambda^{1}+
\lambda^{2})(L_{0})}}\gamma(L_{0},K,U)\,.$$
Then we have
\begin{eqnarray*}
    \langle \Phi_{2}(K),[L,U]\rangle + \langle
    \Phi_{2}([L,K]),U\rangle &=&   \textstyle{\frac 1{(\lambda^{1}+
\lambda^{2})(L_{0})}}\Big(
\gamma(L_{0},K,[L,U])+\gamma(L_{0},[L,K],U)\Big)\\
&=& \textstyle{\frac {(\lambda^{1}+
\lambda^{2})(L)}{(\lambda^{1}+ \lambda^{2})(L_{0})}}\gamma(L_{0},K,U)\ =\
\gamma(L,K,U)
\end{eqnarray*}
by (\ref{gerade}). Hence $(A_{1})\,(ii)$ is satisfied for
$\Phi_{1}$, $\Phi_{2}$ as above. Consequently $\fk=0$.
\qed

\begin{lm}
    If $\fl$ is admissible, then 
$\lambda^{1}+\lambda^{2}\in\{0,-\lambda^{1},-\lambda^{2}\}$.
\end{lm}
\proof Since $\fl$ is admissible we can choose a cohomology class
$[\alpha,\gamma]\in\cH^2_Q(\fl,\fa)_\sharp$ for a suitable orthogonal $\fl$-module $\fa$. By 
Lemma \ref{muecke} we may 
suppose $\alpha(\fl,R)\subset
E_{\lambda^1}+E_{\lambda^2}+E_{\lambda^1+\lambda^{2}}$. As above we
assume $\lambda^{1}\not=0$ since $\fl$ is not nilpotent. Assume now
$\lambda^{1}+\lambda^{2}\notin\{0,-\lambda^{1},-\lambda^{2}\}$. 
Then we have 
\begin{eqnarray}
&\{\lambda^{1}, \overline{\lambda^{1}},
\lambda^{1}+\lambda^{2}\}\ \cap \ \{-\overline{\lambda^{1}},
-\overline{\lambda^{2}},
-(\lambda^{1}+\lambda^{2})\}=\emptyset&\label{cap1}\\
&\{\lambda^{1}, \overline{\lambda^{1}},
\lambda^{1}+\lambda^{2}\}\ \cap \ \{-\overline{\lambda^{1}},
-{\lambda^{1}},
-(\lambda^{1}+\lambda^{2})\}=\emptyset&.\label{cap2}
\end{eqnarray}
{}From (\ref{cap1}) we obtain 
$$E_{\Bbb C}:= 
E_{\lambda^1}+E_{\overline{\lambda^1}}+E_{\lambda^1+\lambda^{2}}
\perp   E_{\lambda^1}+E_{\lambda^2}+E_{\lambda^1+\lambda^{2}}.$$
The vector space $E_{\Bbb C}$ is invariant  under conjugation,
thus it is the complex span of a real subspace $E\subset \fa$. By (\ref{cap2}) 
$E$ is totally isotropic.
We have $E\subset \fb_1$ since $\alpha(\fl,R)\subset
E_{\lambda^1}+E_{\lambda^2}+E_{\lambda^1+\lambda^{2}}$ implies
$\langle \alpha(\fl,R),B \rangle = 0$ for all $B\in E$, hence $\Phi=0$ is a 
solution.
By $(B_1)$ the space $\fb_1$ is non-degenerate, thus we have
$$E':=\fa\cap   (\, E_{-\overline{\lambda^1}}+ 
E_{-\lambda^1}+E_{-(\lambda^1+\lambda^{2})} \,)
\subset \fb_1.$$ Therefore there exists a homomorphism $\Phi\in \Hom(R,\fa)$ 
such that 
$$\langle\alpha (L,U),B\rangle=\langle 
\rho(L)\Phi(U)-\Phi([L,U]),B\rangle=0$$
for all $L\in\fl$, $U\in R$ and $B\in E'$.
In particular, we have
$\langle\alpha (Y,Z),B\rangle=\langle 
\rho(Y)\Phi(Z)-\Phi([Y,Z]),B\rangle=0$
for all $Y,Z\in R$ and
$B\in E_{-(\lambda^1+\lambda^{2})}\subset E'$.
Since $\alpha|_{R\times R}\subset E_{\lambda^1+\lambda^{2}} $ this implies 
$\alpha|_{R\times R}=0$.
According to Lemma \ref{Lk} there are elements $U_1\in 
V_{\lambda^{1}}$ and $L\in\fl$ such that $\alpha(U_1,L)\not= 0$.
By Lemma \ref{fliege} we may assume $\alpha(U_1,L)\in E_{\lambda^1}$. Hence
there exists an element $B\in E_{-\overline{\lambda^1}}$ such that
$\langle \alpha(L,U_1),B\rangle\not= 0$. Then 
$$\langle \alpha(L,U_{1i}),B_j\rangle
=\langle\rho(L)\Phi(U_{1i})-\Phi([L,U_{1i}],B_j \rangle,\quad i,j\in\{1,2\}$$
for $U_{11}=\Re U_1$, $U_{12}=\Im U_1$, $B_1=\Re B\in E'\subset\fb_1$, 
$B_2=\Im B\in E'\subset\fb_1$.
Extending $\Phi$ complex linearly we obtain
$$0\not=\langle\alpha(X,U_1),B\rangle = \langle 
\rho(L)\Phi(U_1)-\Phi([L,U_1]),B\rangle
=\langle (\rho(L)-\lambda^1(L))\Phi(U_1),B\rangle =0$$
since   $B \in E_{-\overline{\lambda^1}}$, which is a contradiction.
\qed
 
Finally we will assume that $\adR$ is not semi-simple. Then $\adR$ 
has a real weight $\lambda\not=0$ with generalised
weight space, i.e.
    $(\adR(L)-\lambda(L))^2=0$ for all $L\in\fl$ but 
    $\adR-\lambda\Id\not=0$. We will prove that $\fl$ 
is not admissible. Assume $\fl$ is admissible and $[\alpha,\gamma]\in 
\cH^2_{Q}(\fl,\fa)_\sharp$ for a suitable orthogonal $\fl$-module 
$\fa$. Let $U\in R$ be a weight vector for $\lambda$. Then 
$R_{2}:=R_{2}(\fl)=\RR\cdot U$. As above one proves that any cocycle $\alpha\in 
Z^{2}(\fl,\fa)$ 
satisfies $\alpha(R,R)\subset E_{2\lambda}$. Furthermore, there exists a cocycle 
$\tilde\alpha\in Z^{2}(\fl,\fa)$ such that $[\alpha]=[\tilde\alpha]\in
H^{2}(\fl,\fa)$ and $\tilde \alpha(R_{2},\fl)\subset 
E_{2\lambda}+E_{\lambda}$. Therefore we may assume $\alpha(R_{2},\fl)\subset 
E_{2\lambda}+E_{\lambda}$ and hence $E_{2\lambda}+E_{\lambda}\subset \fb_{2}$. 
Since $\fb_{2}$ 
is non-degenerate and $\lambda\not=0$ we obtain 
$E_{-2\lambda}+E_{-\lambda}\subset \fb_{2}$. Because of $E_{-2\lambda}\subset 
\fb_{2}$ now Condition $(B_{2})$ gives $\alpha|_{R\times R}=0$. By 
Lemma \ref{fliege} we may assume $\alpha(R_{2},\fl)\subset 
E_{\lambda}$.  Now $E_{-\lambda}\subset \fb_{2}$ implies 
$\alpha(R_{2},\fl)=0$, which is a contradiction to~$(A_{1})$.

{\sl Proof of Proposition \ref{elefant}.} By the above we know that 
$\adR$ is semi-simple and 
the weights $\lambda^1$ and $\lambda^2$ of $\adR$ satisfy $\lambda^{1}+\lambda^{2}=0$ or 
(possibly after change 
of the numbering) 
$\lambda^{1}+2\lambda^{2}=0$. In particular, 
$\lambda^{1},\lambda^{2}\not=0$, which implies that $\adR$ does not 
have invariants. Therefore we have $H^{2}(\fl/R,R)=0$. Hence there 
exists an abelian subalgebra $\fl_1\subset\fl$ such that $\fl= 
\fl_{1}{\;}_{\adR |_{\fl_{1}}}\hspace{-0.5em}\ltimes R$. Since 
$\lambda^1$ and $\lambda^2$ are linearly dependent the codimension of
$\ker (\adR |_{\fl_{1}})$ in $\fl_1$ is one. Hence we may choose an 
element $X_1\in\fl_1\setminus \ker (\adR |_{\fl_{1}})$  such that 
$\fl= (\RR \cdot X_{1}{\;}_{\adR}\hspace{-0.3em}\ltimes 
R)\oplus \RR^{k}$ and \\[1ex]
${}\qquad\quad\adR(X_{1})= ${\small $
\left(
\begin{array}{cc}
    0 &-1  \\
    1 &  0
\end{array}
\right)
 $} \,or\,
 $\adR(X_{1})= ${\small $
\left(
\begin{array}{cc}
    1 &0  \\
    0 &  -1
\end{array}
\right)
 $} \,or\,
 $\adR(X_{1})= ${\small $
\left(
\begin{array}{cc}
    -2 &0  \\
    0 &  1
\end{array}
\right)
 $}\\[1ex]
with respect to a suitable basis of $R$.
\qed

For the sake of completeness we will also give the classification result 
in the nilpotent case. However, we will omit the proof since on the one 
hand we will not use the result in this paper and on the other hand 
our proof is tricky and not very enlightening.

\begin{pr}\label{R2n}
If $\fl$ is an admissible nilpotent Lie algebra with $\dim R(\fl)=2$,  then $\fl$ 
is isomorphic to one of the
following Lie algebras
\begin{eqnarray*}
    &&\{ [X_1,Z]=Y,\ [X_1,X_2]=Z\}\oplus \RR^{k},\\
    &&\{ [X_1,X_2]=Y,\ [X_1,X_3]=Z\}\oplus \RR^{k},\\
    &&\{ [X_1,X_2]=Y,\ [X_3,X_4]=Z\}\oplus \RR^{k}\,=\,\fh(1)\oplus\fh(1)\oplus\RR^k,\\
    &&\{ [X_1,X_2]=Y,\ [X_1,X_3]=Z,\ [X_3,X_4]=Y\}\oplus \RR^{k},\\
    &&\{ [X_1,X_2]=Y,\ [X_1,X_3]=Z,\ [X_2,X_4]=Z,\ [X_3,X_4]=\pm Y\}\oplus \RR^{k}.
\end{eqnarray*}
\end{pr}

\section{Metric Lie algebras of index 3}\label{i3}
The aim of this section is the classification of all indecomposable metric Lie algebras of 
index 3. 
\subsection{Preliminaries}
\begin{pr}\label{jonas}
If $(\fg,\ip)$ is a non-simple indecomposable metric Lie algebra of index 3, then
$\fg/\fj(\fg)$ is isomorphic to one of the Lie algebras
$\fn(2)$, $\fr_{3,-1}$, $\fh(1)$, $\fsl(2,\RR)$, $\fsu(2)$ or $\RR^k$, $k=1,2,3$.
\end{pr}
\proof 
Let $(\fg,\ip)$ be a non-simple indecomposable metric Lie algebra of index 3.
By Proposition \ref{schmidt} and Proposition \ref{Pequiv} $(\fg,\ip)$  has the 
structure of
a balanced quadratic extension $\dd$, where $\fl=\fg/\fj(\fg)$. In particular, 
$\fl=\fg/\fj(\fg)$ is admissible.
First we will determine all
admissible $\fl$ which can appear in such a quadratic extension.
Obviously $\dim \fl\le 3$ because $\ip$ has index $3$.
In particular, $\dim R(\fl)\in\{0,1,2\}$. If $\dim R(\fl)=0$, then $\fl$ is  abelian or simple. If $\dim R(\fl)=1$, then $\fl\cong \fh(1)$ by Proposition~\ref{R1}. Finally, suppose $\dim R(\fl)=2$. Then $\dim \fl=3$. In particular, the codimension of $R(\fl)=\fl'$ in $\fl$ is one and therefore $[\fl,\fl']=[\fl,\fl]$, hence $\fl$ is not nilpotent. Now Proposition \ref{elefant} implies that $\fl$ is one of the 3-dimensional  Lie
algebras
$\fn(2)$, $\fr_{3,-1}$,
$ \fr_{3,-2}$.

On the other hand we must take into consideration that $\fa$ must be
Euclidean if $\fl$ is three-dimensional. This excludes
$\fl=\fr_{3,-2}$. Indeed, the weights $\lambda^1$ and $\lambda^{2}$ of
the representation $\adR$ of $\fr_{3,-2}$ are given by
$\lambda^{1}(X)=-2$ and $\lambda^{2}(X)=1$. In particular they are real
and satisfy $\lambda^{1},\ \lambda^{2}\not=0$ and $\lambda^{1}+\lambda^{2}\not=0$.
Hence, if $\fa$ is Euclidean, then
$E_{\lambda^{1}+\lambda^{2}}+E_{\lambda^{1}}+E_{\lambda^{2}}=0$.
If we now assume that $\dd$ is balanced, then this together with
Lemma \ref{muecke} implies $[\alpha]=0$ in $H^{2}(\fl,\fa)$. Hence
$R(\fl)=0$ by Lemma \ref{Lk}, which is a contradiction.
\qed

It remains to classify all indecomposable metric Lie algebras $(\fg,\ip)$ of index
3 for which $\fg/\fj(\fg)$ is isomorphic to one of the Lie algebras $\fsl(2,\RR),\
\fsu(2)$, $\fn(2),\
\fr_{3,-1},\ \fh(1)$ or $\RR^{k}$ for $k\le 3$. Before we will start we
prove the following fact on 3-dimensional Lie algebras.
\begin{lm}\label{maus}
    Let $\fl$ be a 3-dimensional Lie algebra such that $R:=R(\fl)$
    is 2-dimensional and $\adR$ is semi-simple. We denote the weights of the
    complexified representation $\adR$ by $\lambda^{1}$
    and $\lambda^{2}$. Let $\fa$ be a semi-simple orthogonal $\fl$-module and
    $V_{\lambda}$ be defined as above. Choose an element $X\in
    \fl\setminus R$. Then
    \begin{eqnarray*} &\{\alpha\in
    C^2(\fl,\fa)\mid\alpha(R,R)\subset E_{\lambda^{1}+\lambda^{2}},\
    \alpha(X,V_{\lambda^{i}})\subset E_{\lambda^{i}},\ i=1,2\}
    \longrightarrow
    H^2(\fl,\fa)&\\
    &\alpha \longmapsto [\alpha]&
\end{eqnarray*}
 is well-defined and an isomorphism.
\end{lm}
\proof
The map is well defined since $\alpha(R,R)\subset E_{\lambda^{1}+\lambda^{2}}$ 
implies
$\alpha\in Z^{2}(\fl,\fa)$. Let us prove that the map is surjective.
For a given cohomology class
$a=[\alpha]\in H^2(\fl,\fa)$ we define cocycles $\alpha_{1}$,
$\alpha_{2}\in Z^2(\fl,\fa)$ by
$$  \alpha=\alpha_{1}+\alpha_{2},\quad  \alpha_{1}(X,\cdot)=0,\quad
    \alpha_{2}|_{R\times R}=0 \,.
$$
Then $\alpha_{1}(R,R)\subset E_{\lambda^{1}+\lambda^{2}}$ and by
Lemma
\ref{fliege} there is a cocycle $\tilde\alpha_{2}\in
Z^{2}(\fl,\fa)$ such that $[\alpha_{2}]=[\tilde\alpha_{2}]\in
H^{2}(\fl,\fa)$ and $\tilde\alpha_{2}|_{R\times R}=0$,
$\tilde\alpha_{2}(\fl,V_{\lambda^{i}})\subset E_{\lambda^{i}}$ for $i=1,2$.
Then $\alpha_{1}+\tilde\alpha_{2}$ is a preimage of $a$.

It remains to show that the map is injective. Assume that
$\alpha=d\tau$ satisfies $d\tau(X,U_{i})\in E_{\lambda^{i}}$
for $U_{i}\in V_{\lambda^{i}}$, $i=1,2$. Since
$d\tau(X,U_{i})=(\rho(X)-\lambda^i)(\tau(U_{i}))$ we have on the
other hand $p_{\lambda^{i}}(d\tau(X,U_{i}))=0$. Hence,
$d\tau(X,U_{i})=0$ and therefore $\alpha=d\tau=0$.
\qed

\begin{re}\label{HS}
{\rm
    The assertion of Lemma \ref{maus} can also be verified using the Hochschild--Serre spectral 
    sequence associated with the ideal $R\subset\fl$. For
    dimensional reasons this spectral sequence leads to the exact sequence
    $$0\longrightarrow H^1(\fl/R,H^1(R,\fa))\longrightarrow H^2(\fl,\fa) \longrightarrow
    H^0(\fl/R,H^2(R,\fa))\longrightarrow 0.$$
    Since $H^1(R,\fa)\cong C^1(R,\fa^R)=C^1(R,\fa)$ we have
    \begin{eqnarray*}
    H^1(\fl/R,H^1(R,\fa)) &\cong&
    C^1(\RR\cdot X, C^1(R,\fa)^{\fl/R})
    \\
    &\cong&\{\alpha\in C^2(\fl,\fa)\mid\alpha(X,V_{\lambda^i})\subset E_{\lambda^i},
    i=1,2\}\,.
    \end{eqnarray*}
    On the other hand
    \begin{eqnarray*}
    H^0(\fl/R,H^2(R,\fa)) &\cong&H^2(R,\fa)^{\fl/R}\cong C^2(R,\fa)^{\fl/R}\\
    &\cong& \{\alpha\in C^2(R,\fa)\mid \alpha(R,R)\subset E_{\lambda^1+\lambda^2}\}
    \end{eqnarray*}
    and we obtain the assertion of Lemma \ref{maus}.
    }
\end{re}

\begin{lm}\label{buten} Let $\fl$ be a 3-dimensional unimodular Lie
algebra.
For an orthogonal $\fl$-module
$\fa$ the map
\begin{eqnarray*}\iota_{Q}:\  \cH^2_Q(\fl,\fa) &\longrightarrow&
(\,H^2(\fl,\fa)\setminus\{0\}\,) \ \cup\ C^3(\fl)\\
\ [\alpha,\gamma] &\longmapsto & \left\{\begin{array}{ll}
[\alpha]\in H^2(\fl,\fa) &\mbox{if } [\alpha]\not= 0\\
\gamma  \in C^3(\fl) &\mbox{if } [\alpha]= 0
\end{array} \right.
\end{eqnarray*}
is a bijection.
\end{lm}

\proof To prove this we will use Proposition
\ref{club}. Let us determine the vector space $H^{3}(\fl)/a\cup H^{1}(\fl,\fa)$ for
$a=[\alpha]\in
H^2_\cup(\fl,\fa)=H^{2}(\fl,\fa)$. Since $\fl$ is unimodular
$d\sigma(X,Y,Z)=-\tr (\adR(X))\cdot \sigma(Y,Z)=0$ holds for all $\sigma\in
C^{2}(\fl,\fa)$. Hence $H^{3}(\fl)\cong C^{3}(\fl)$ is one-dimensional. By (\ref{poincare}) we 
know that $\cup: H^{2}(\fl,\fa)\otimes H^{1}(\fl,\fa)\rightarrow
H^3(\fl)$ is a non-degenerate pairing. If now $a\not=0$, then this implies
$H^{3}(\fl)/a\cup
H^{1}(\fl,\fa)=0$. Hence for $a\not= 0$ the set $p^{-1}(a)$ only
consists of one element, namely $[\alpha,0]\in\cH^2_Q(\fl,\fa)$. Now
suppose $a=0$. Then
we have $H^{3}(\fl)/0\cup
H^{1}(\fl,\fa)=H^{3}(\fl)\cong C^{3}(\fl)$. Hence, $p^{-1}(0)\ni[0,\gamma]\mapsto
\gamma\in
C^{3}(\fl)$ is a bijection.
\qed

Let $\RR^{p,q}$ be the standard pseudo-Euclidean space of
dimension $n=p+q$. As usual, $\RR^{n}:=\RR^{0,n}$. We identify  $\RR^{p_1,q_1}\oplus\RR^{p_2,q_2}$ with $\RR^{p_1+p_2,q_1+q_2}$.

\begin{de}
Let $\fl_{0}$ be an abelian Lie algebra and $\lambda\in\fl_{0}^{*}$.
We define orthogonal representations $\rho^+_{\lambda}$ of $\fl_{0}$ on $\RR^2$,
$\rho^-_{\lambda}$ of $\fl_{0}$ on $\RR^{2,0}$ and $\rho'_{\lambda}$ of $\fl_{0}$ 
on $\RR^{1,1}$ by
\begin{eqnarray*}
&\rho^\pm_{\lambda}(L)=\small{
\left(\begin{array}{cc}
   0&-\lambda(L)\\
   \lambda(L)&0
   \end{array}
\right)},\quad
\rho'_{\lambda}(L)=\small{
\left(\begin{array}{cc}
   0&\lambda(L)\\
   \lambda(L)&0
   \end{array}
\right)}&
\end{eqnarray*}
w.\,r.\,t. an orthonormal basis of $\RR^2$, $\RR^{2,0}$, and $\RR^{1,1}$, 
respectively.

Moreover, for $\mu,\,\nu\in\fl_{0}^{*}$ we define an orthogonal representation 
$\rho''_{\mu,\nu}$ of
$\fl_0$ on $\RR^{2,2}$ by
$$
\rho''_{\mu,\nu}(L)=\small{
        \left(
        \begin{array}{cccc}
            0 & -\nu(L) & \mu(L) & 0  \\
            \nu(L) & 0 & 0 & \mu(L) \\
            \mu(L) & 0 & 0 & -\nu (L) \\
            0 & \mu(L) & \nu(L) & 0
        \end{array}
            \right)}
$$
w.\,r.\,t. an orthonormal basis of $\RR^{2,2}$.

For 
$\lambda=(\lambda^{1},\ldots,\lambda^{m}),\,\mu=(\mu^{1},\ldots,\mu^{m}),\,\nu=(\nu
^{1},\ldots,\nu^{m})\in (\fl_{0}^{*})^{m}$ we define semi-simple
orthogonal representations $\rho^+_{\lambda}$ of $\fl_{0}$ on $\RR^{2m}$,
$\rho^-_{\lambda}$ of $\fl_{0}$ on $\RR^{2m,0}$, $\rho'_{\lambda}$ of $\fl_{0}$ on 
$\RR^{m,m}$ and $\rho''_{\mu,\nu}$ of $\fl_{0}$ on $\RR^{2m,2m}$ by
$$\rho^\pm_\lambda=\bigoplus_{i=1}^m \rho^\pm_{\lambda^i}\,,\quad
\rho'_\lambda=\bigoplus_{i=1}^m \rho'_{\lambda^i}\,,\quad
\rho''_{\mu,\nu}=\bigoplus_{i=1}^m \rho''_{\mu^i,\nu^i}.$$

Now let $\fl$ be a solvable Lie algebra, $\fl_{0}=\fl/R(\fl)$ and let $\rho^+_\lambda$, $\rho^-_\lambda$, $\rho'_\lambda$ and $\rho''_{\mu,\nu}$ be
the above defined representation of $\fl_0$. Composing the projection $\fl\rightarrow \fl_0$ with these representations we obtain representations of $\fl$, which we denote by the same symbol.

Finally we denote by $\rho_{0}$ the trivial representation of $\fl$
on $\fa=\RR^{p,q}$.
\end{de}
The symmetric group $\frak S_{m}$ acts on $(\fl_{0}^{*})^{m}$ by permuting 
coordinates and on
$(\fl_{0}^{*})^{m}\oplus(\fl_{0}^{*})^{m}$ by permuting pairs of coordinates. The 
group $(\ZZ_{2})^{m}$ acts on $(\fl_{0}^{*})^{m}$ by changing the signs of the 
coordinates. We define the orbit spaces
$\Lambda_m:=(\fl_0^*\setminus  0)^{m}/\frak S_{m}\ltimes (\ZZ_{2})^{m}$ and
$\Lambda''_m:=\Big( ((\fl_0^*\setminus 
0)^{m}/(\ZZ_{2})^{m})\oplus((\fl_0^*\setminus 0)^{m}/(\ZZ_{2})^{m}) \Big)/\frak 
S_{m}$. Finally we define an action of $\Aut(\fl)$ on $\Lambda_m$ and $\Lambda''_m$ 
by
$S^*[\lambda]:=[S^*\lambda]$ and $S^*[\mu,\nu]:=[S^*\mu,S^*\nu]$.
\begin{pr} \label{PHom} For a solvable Lie algebra $\fl$ we consider the map
 \begin{eqnarray*}
      \bigcup _{\scriptsize{
      \begin{array}{c}
      (m_1,\ldots,m_4,p_0,q_0)\\
      2m_1+m_3+2m_4+p_0=p\\
      2m_2+m_3+2m_4+q_0=q
      \end{array}}}
      \Lambda_{m_1}\times\Lambda_{m_2}\times 
\Lambda_{m_3}\times\Lambda''_{m_4}&\longrightarrow&
      \Hom(\fl,\so(p,q))_{\rm ss} /O(p,q)\\
      ([\lambda_1],[\lambda_2],[\lambda_3],[\mu,\nu])&\longmapsto&
      [\rho^+_{\lambda_1}\oplus \rho^-_{\lambda_2} \oplus \rho'_{\lambda_3} \
      \oplus    \rho''_{\mu,\nu}\oplus\rho_0]\,,
\end{eqnarray*}
where $\rho_0$ is the trivial representation of $\fl$ on $\RR^{p_0,q_0}$. This map is a 
bijection. It is equivariant with respect to the action of $\Aut(\fl)$.
\end{pr}

We call an
orthogonal basis $A_{1},\ldots,A_{n}$ of $\RR^{p,q}$, $p+q=n$, orthonormal if $\langle
A_{k},A_{k}\rangle =-1$ for $k=1,\ldots,p$ and $\langle
A_{k},A_{k}\rangle =1$ for $k=p+1,\ldots,n$.
For later use we fix an orthonormal basis $A_1^+,\dots,A_{2m}^+$ of the $\fl$-module $(\rho_\lambda^+,\RR^{2m})$ and an orthonormal basis $A_1^0,\dots,A_{p+q}^0$ of the trivial $\fl$-module $(\rho_0,\RR^{p,q})$.

\label{Sindex3}

\subsection{The case $\fl=\fn(2)$}
\label{S71}
First we consider $\fl=\fn(2)$.
We will write elements of $\Aut(\fl)$ as matrices w.~r.~t.~the basis $X,Y,Z$ of 
$\fl$.
$$\Aut(\fl)=\Big\{S(u,v,w,a,b):=
{\small
\left( \begin{array}{ccc}
u&0&0\\
v&a&-b\\
w&ub&ua
\end{array}\right)
}
\ \Big|\ u=\pm1,\, a,b,v,w\in\RR,\ a^{2}+b^{2}\not=0\Big\}$$
Let $\lambda$ and $\bar \lambda=-\lambda$ be the weights of the representation 
$\adR$ of $\fl_{\Bbb C}/R_{\Bbb C}$ on $R_{\Bbb C}$. Then 
$\lambda(X)=i$. We use again the notation
$E=(E_\lambda\oplus E_{\bar \lambda})\cap\fa\subset \fa_{\Bbb C}$. Next 
we will determine $H^2(\fl,\fa)$, $\cH^2_Q(\fl,\fa)$, and  
$\cH^2_Q(\fl,\fa)_\sharp$ for a Euclidean semi-simple orthogonal $\fl$-module
$(\rho,\fa)$.
By Lemma
\ref{maus} the map
\begin{eqnarray*} &\{\alpha\in C^2(\fl,\fa)\mid\alpha(Y,Z)\in\fa^\fl,
 \alpha(X,Y)\in E,
\alpha(X,Z)=\rho(X)\alpha(X,Y)\}\longrightarrow
H^2(\fl,\fa)&\\
&\alpha \longmapsto [\alpha]&
\end{eqnarray*} is an isomorphism. We denote the inverse of this isomorphism by 
$\iota$.
\begin{lm}If $\fa$ is a Euclidean orthogonal $\fl$-module, then we have
$$\cH^2_Q(\fl,\fa)_\sharp=  \cH^2_Q(\fl,\fa)\setminus 
\{[0,0]\}=\iota_Q^{-1}(\,(H^2(\fl,\fa)\setminus{0})\cup(C^3(\fl)\setminus{0})\,).$$
\end{lm}
\proof
We have to check under which conditions a 
cohomology class $[\alpha,\gamma]\in  \cH^2_Q(\fl,\fa)$ is admissible.
First we note that $(A_0)$ and $(A_k)$, $k>1$ are satisfied for all 
$[\alpha,\gamma]$ since
$\fz(\fl)=0$ and $R_k(\fl)=0$ for $k>1$. Conditions $(B_0)$ and $(B_k)$,
$k\ge1$ are also satisfied
since $\fa$ is Euclidean. It remains to check $(A_1)$ for $\fk=R$ since $R$
is the only non-vanishing ideal in $R$. We may assume
$\alpha=\iota([\alpha])$. In particular, $\alpha(X,Y-iZ)\in E_\lambda$.
On the other hand  $(A_1)(i)$  is satisfied if and only if $\alpha(Y,Z)=0$ and
$\alpha(X,Y-iZ)=(\rho(X)-i)\Phi_1(X-iY)$  for a homomorphism
$\Phi_1\in \Hom(R,\fa)$, thus if and only if $\alpha=0$. Hence all 
$[\alpha,\gamma]$
with $[\alpha]\not=0$ are admissible. 
Now we assume $\alpha=0$. Obviously, $[0,0]$ is not admissible. 
Suppose that $\gamma\not=0$. Assume that $[0,\gamma]$ satisfies the assumption 
$(ii)$ of 
$(A_1)$. Then there is a homomorphism $\Phi_2\in \Hom(\fk,R^*)$ such that 
$$\gamma(X,Y,Z)=\langle \Phi_2(Y),[X,Z]\rangle+
\langle\Phi_2([X,Y]),Z\rangle = -\langle \Phi_2(Y),Y\rangle+
\langle\Phi_2(Z),Z\rangle$$
and    
$$\gamma(X,Z,Y)=\langle \Phi_2(Z),[X,Y]\rangle+
\langle\Phi_2([X,Z]),Y\rangle = \langle \Phi_2(Z),Z\rangle-
\langle\Phi_2(Y),Y\rangle.$$  
Since on the other hand $\gamma(X,Y,Z)=-\gamma(X,Z,Y)$ these 
equations imply $\gamma=0$, a contradiction. Thus $[0,\gamma]$ does not satisfy the 
assumption $(ii)$ of 
$(A_1)$. Consequently, the implication $(A_{1})$ is true in this case. 
Hence, $[0,\gamma]$ is admissible for all $\gamma\not=0$.
\qed

Since here $\fl_{0}=\fl/R(\fl)=\RR\cdot
X$ we identify $\lambda\in
(\fl_{0}^{*})^{m}$ with $\lambda(X)\in\RR^{m}$.
\begin{pr}\label{Pn2}
    If $(\fg,\ip)$ is an indecomposable metric Lie algebra of index 3 such that
    $\fg/\fj(\fg)\cong\fl:=\fn(2)$, then $\fg$
    is isomorphic to exactly one of the following indecomposable 
    metric Lie algebras $\dd$: 
    \begin{itemize}
        \item[(Ia)] $\fa=\RR^{2m}$, $m\ge 0$, $\rho=\rho^+_{\lambda}$,\\ where
        $\lambda=(\lambda^{1},\ldots,\lambda^{m})$, $0<\lambda^{1}\le\lambda^{2}\le 
\ldots \le\lambda^{m}$,\\
        $\alpha=0$,\\
        $\gamma(X,Y,Z)=1;$
        \item[(Ib)] as above but $\gamma(X,Y,Z)=-1;$
        \item[(II)]$\fa=\RR^{2m+1}=\RR^{2m}\oplus\RR^{1}$, $m\ge 0$, 
        $\rho=\rho^+_{\lambda}\oplus\rho_{0}$,\\ where
        $\lambda=(\lambda^{1},\ldots,\lambda^{m})$, 
        $0<\lambda^{1}\le\lambda^{2}\le \ldots \le\lambda^{m}$,\\
        $\alpha(Y,Z)=A^{0}_{1}$, $\alpha(X,Y)=\alpha(X,Z)=0$,\\
        $\gamma=0$;
        \item[(III)] $\fa=\RR^{2m+3}=\RR^{2m+2}\oplus\RR^{1}$, $m\ge 0$,
        $\rho=\rho^+_{\lambda'}\oplus\rho_{0}$,\\
        where $\lambda'=(\lambda^{1},\ldots,\lambda^{m},1)$,
        $0<\lambda^{1}\le\lambda^{2}\le \ldots \le\lambda^{m}$,\\
        $\alpha(Y,Z)=A_{1}^{0}$, $\alpha(X,Y)=r A^+_{2m+1}$, $\alpha(X,Z)=r
        A^+_{2m+2}$, $r>0$;\\
        $\gamma=0$;
        \item[(IV)]$\fa=\RR^{2m+2}$, $m\ge 0$, $\rho=\rho^+_{\lambda'}$, \\ where
        $\lambda'=(\lambda^{1},\ldots,\lambda^{m},1)$,
        $0<\lambda^{1}\le\lambda^{2}\le \ldots \le\lambda^{m}$,\\
        $\alpha(Y,Z)=0$, $\alpha(X,Y)= A^+_{2m+1}$, $\alpha(X,Z)=
        A^+_{2m+2}$;\\
        $\gamma=0$.
    \end{itemize}
\end{pr}
\proof We already know that $(\fg,\ip)$ is a quadratic extension of $\fl$ by a suitable 
orthogonal $\fl$-module $\fa$. Since $\ip$ has index 3, $\fa$ must be Euclidean, i.e. 
$\fa=\RR^{n}$. In
particular we have $\Hom(\fl,\so(\fa,\ip_\fa))_{\rm ss}=\Hom(\fl,\so(n))$. If 
$\rho\in \Hom(\fl,\so(n))$ and $S\in\Aut(\fl)$, then $\rho\circ
S=\pm\rho$ holds. On the other hand, if $\rho\in \Hom(\fl,\so(n))$, then there 
exists a map $U\in O(n)$ such that $\Ad (U)\rho=-\rho$. This implies
$$\Hom(\fl,\so(n))/G=\Hom(\fl,\so(n))/(\ZZ_{2}\times
O(n))=\Hom(\fl,\so(n))/O(n).$$ By Proposition \ref{PHom} the map
\begin{eqnarray*}
    \bigcup_{2m\le n} \Lambda_m
    &\longrightarrow
     &\Hom(\fl,\so(n))/G\\
     \,[\lambda]&\longmapsto&[\rho^+_{\lambda}\oplus\rho_{0}]
\end{eqnarray*}
is a bijection. Since we can identify $\lambda\in\Lambda_m$  with 
$\lambda(X)\in\RR^m$, we can also identify
$$\Lambda_m=\{\lambda\in\RR^{m}\mid  0<\lambda^{1}\le\lambda^{2}\le \ldots
     \le\lambda^{m}\}.$$
Hence each $G$-orbit in $\Hom(\fl,\so(n))$ has a canonical representative 
$\rho^+_{\lambda}\oplus\rho_{0}$ with $\lambda\in\Lambda_m$. We fix a 
representation
$\rho:=\rho^+_{\lambda}\oplus\rho_{0}$. Next we will describe the action of
$G_{\rho}$ on $\cH^2_Q(\fl,\fa_{\rho})$ identifying $\cH^2_Q(\fl,\fa_{\rho})$
with $(\,H^2(\fl,\fa_{\rho})\setminus\{0\}\,) \ \cup\ C^3(\fl)$ via
$\iota_{Q}$ and $H^2(\fl,\fa_{\rho})\setminus\{0\}$ with 
$\iota(H^2(\fl,\fa_{\rho})\setminus\{0\})$.
We claim that for $(S,U)\in G_{\rho}$ with
$S=S(u,v,w,a,b)$ and
$\alpha\in \iota(\,H^2(\fl,\fa_{\rho})\setminus\{0\}\,)$, $\gamma\in
C^{3}(\fl)$ the following holds:
\begin{eqnarray}
\iota\circ \iota_{Q}((S,U)^* \iota_{Q}^{-1}([\alpha])) &=&\tilde 
\alpha\quad\mbox{with }\nonumber\\
\tilde\alpha(Y,Z)&=&u(a^2+b^{2})U^{-1}(\alpha(Y,Z))\,\in\, \fa^{\fl} 
\label{tilde1}\\ 
\tilde\alpha(X,Y)&=&u(a+b\rho(X))U^{-1}(\alpha(X,Y))\,\in\, E \label{tilde2}\\ 
\tilde\alpha(X,Z)&=&\rho(X)(\tilde\alpha(X,Y))\label{tilde3}\\[1ex]
\iota_{Q}((S,U)^*
\iota_{Q}^{-1}(\gamma))&=&(a^{2}+b^{2})\gamma.\label{gamma}
\end{eqnarray}
Let us verify this. The condition 
$(S,U)\in G_{\rho}$ says that $\Ad(U^{-1})(\rho(S(X)))=\rho(X)$ holds, which 
here is equivalent to $U^{-1}\rho(uX) U=\rho(X)$. We note that on the one hand
$$S(u,v,w,a,b)=S(1,uv,uw,1,0)\cdot S(u,0,0,a,b)$$
and that on the other hand by Proposition \ref{inner} the subgroup
$$G_\rho':=\{(S(1,v,w,1,0),\Id)\mid v,w\in\RR\}\subset
\{(e^{\ad(L)},e^{\rho(L)})\mid L\in\fl\}\subset G_\rho$$
acts trivially on $\cH^2_Q(\fl,\fa_{\rho})$. Therefore it suffices to prove
$$\iota\circ\iota_{Q}((S_0,U)^* \iota_{Q}^{-1}([\alpha]))=\tilde \alpha,\quad 
\iota_{Q}((S_0,U)^*
\iota_{Q}^{-1}(\gamma))=(a^{2}+b^{2})\gamma$$
for $S_0=S(u,0,0,a,b)$.
We have
$$\iota\circ\iota_Q((S_0,U)^*\iota_{Q}^{-1}([\alpha]))
=\iota_{Q}(U_{*}^{-1}S_0^{*}[\alpha,0])=\iota_{Q}([U_{*}^{-1}S_0^{*}\alpha,0])
=\iota(U_{*}^{-1}S_0^{*}\alpha)$$
and we calculate
\begin{eqnarray*}
    (U_{*}^{-1}S_0^{*}\alpha)(Y,Z) &=& U^{-1}(\alpha(S_0Y,S_0Z))\\
    &=& u(a^2+b^{2})U^{-1}(\alpha(Y,Z))\\[1ex]
    (U_{*}^{-1}S_0^{*}\alpha)(X,Y) &=& U^{-1}(\alpha(S_0X,S_0Y))\\
    &=& U^{-1}(au\,\alpha(X,Y)+u^{2}b\,\alpha(X,Z))\\
    &=& U^{-1}(u(a+b\rho(uX))\alpha(X,Y))\\
    &=& u(a+b\rho(X))U^{-1}(\alpha(X,Y))\\[1ex]
    (U_{*}^{-1}S_0^{*}\alpha)(X,Z)
    &=&u(-b+a\rho(X))U^{-1}(\alpha(X,Y)).
\end{eqnarray*}
Then $\tilde \alpha:=U_{*}^{-1}S_0^{*}\alpha$ satisfies
(\ref{tilde1}), (\ref{tilde2}), and (\ref{tilde3}).
Hence, $\tilde\alpha=\iota(U_{*}^{-1}S_0^{*}\alpha)$.  Finally,
$$\iota_{Q}((S_0,U)^* \iota_{Q}^{-1}(\gamma))=S_0^{*}\gamma=\det
S_0\cdot\gamma=(a^{2}+b^{2})\gamma,
$$ 
which proves the claim. 

Using this description of the $G_{\rho}$-action we can  
distinguish between the following types of $G_{\rho}$-orbits in 
$\cH^2_Q(\fl,\fa_{\rho})_{0}$, which we characterise by properties of 
their elements $[\alpha,\gamma]$, where we may assume
that $\alpha=\iota([\alpha])$ and, moreover, that $\gamma=0$ if $\alpha\not=0$:
\begin{tabbing}
   Type (Ia)\,:\quad \=$\alpha=0,\ \gamma(X,Y,Z)> 0$,\\
   Type (Ib)\,:\quad \> $\alpha=0,\ \gamma(X,Y,Z)< 0$,\\
   Type (II)\,:\>$\alpha(Y,Z)\not=0,\ \alpha(X,Y)=\alpha(X,Z)=0,\
                \gamma=0$,\\
   Type (III)\,:\>$\alpha(Y,Z)\not=0,\ \alpha(X,\cdot)\not=0 $ on $R$,\
                $\gamma=0$,\\
   Type (IV) \,:\>$\alpha(Y,Z)=0,\ \alpha(X,\cdot)\not=0 $ on $R$,\
                $\gamma=0.$
\end{tabbing}
Next we will classify the $G_{\rho}$-orbits of each type. The result will
give the Lie algebras in the corresponding item of the proposition.

Since for $\fl=\fn(2)$ each decomposition $\fl=\fl_1\oplus\fl_{2}$ of 
Lie  algebras is trivial Definition \ref{kosel} says that a 
cohomology class $[\alpha,\gamma]\in \cH^2_Q(\fl,\fa_{\rho})_{\sharp}$ 
is decomposable if and only if $\fa_{\rho}^{\fl}\cap 
(\alpha(\fl,\fl))^{\perp}=0$. Therefore $[\alpha,\gamma]\in
\cH^2_Q(\fl,\fa_{\rho})_{\sharp}$ is indecomposable if and only if
$\fa_{\rho}^{\fl}=\RR\cdot \alpha(Y,Z)$. In particular, we may assume
$n=2m$ or $n=2m+1$ since otherwise  $\cH^2_Q(\fl,\fa_{\rho})_0$ is empty.

We start with orbits of type (Ia) and type (Ib). Here we may assume 
$n=2m$.
If $S=S(1,0,0,a,b)\in\Aut(\fl)$ then $(S,\Id)\in G_{\rho}$. This together with (\ref{gamma}) 
implies that two
elements
$[0,\gamma_{1}]$ and $[0,\gamma_{2}]$ of  $\cH^2_Q(\fl,\fa_{\rho})_0$ are 
in the same $G_{\rho}$-orbit of type (Ia)  (or of type (Ib)) if and only 
if $\gamma_{1}(X,Y,Z)$ and $\gamma_{2}(X,Y,Z)$ have the same sign. This yields a
classification of orbits of type (Ia) and type (Ib).

Now we consider orbits of type (II) or (III). Here we may assume 
$n=2m+1$.
Besides $(S,\Id)\in G_{\rho}$ for
$S=S(1,0,0,a,b)\in\Aut(\fl)$ we also have $(S,-\Id)\in G_{\rho}$. Now
(\ref{tilde1}) 
implies that each orbit of type 
(II) or (III) contains an element $[\alpha,0]$ with $\alpha=\iota([\alpha])$
and $\alpha(Y,Z)=A_{1}^0$. For $G_{\rho}$-orbits of type (II) this yields
the  claimed classification. Now consider elements $[\alpha_{1},0]$,
$[\alpha_{2},0]$ which belong to $G_{\rho}$-orbits of type (III) and 
satisfy $\alpha_{i}=\iota([\alpha_{i}])$ and 
$\alpha_{i}(Y,Z)=A_{1}^0$, $i=1,2$. Assume $[\alpha_{1},0]$ and
$[\alpha_{2},0]$ are in the same $G_{\rho}$-orbit, i.e. there is an 
element $(S,U)\in G_{\rho}$, $S=S(u,v,w,a,b)$ such that 
$\iota_{Q}((S,U)^* \iota_{Q}^{-1}(\alpha_{1})) =\alpha_{2}$. 
Then $a^{2}+b^{2}=1$ by (\ref{tilde1}). Now (\ref{tilde2}) 
implies
\begin{equation}\label{penny}
    |\alpha_{1}(X,Y)|=|\alpha_{2}(X,Y)|
\end{equation}
since $\rho$ is 
orthogonal and $\rho^{2}=-\Id$ on $E$. Hence, (\ref{penny}) is a
necessary condition for $[\alpha_{1},0]$ and 
$[\alpha_{2},0]$ being in the same $G_{\rho}$-orbit. We will show 
that it is also sufficient. Assume (\ref{penny}) is satisfied. Since, 
furthermore, $\alpha_{i}(X,Z)=\rho(X)\alpha_{i}(X,Y)$ for $i=1,2$ we can define an 
orthogonal map $U$ on $E$ which commutes with $\rho(X)$ such that 
$U(\alpha_{2}(X,Y))=\alpha_{1}(X,Y)$ and
$U(\alpha_{2}(X,Z))=\alpha_{1}(X,Z)$. We extend $U$ to a map 
$U_{0}\in O(n)$ such that $U_{0}|_{E^{\perp}}=\Id$. Then $(\Id,U_{0})\in G_{\rho}$ 
and $\iota_{Q}((\Id,U_{0})^* \iota_{Q}^{-1}(\alpha_{1})) =\alpha_{2}$, 
hence $[\alpha_{1},0]$ and $[\alpha_{2},0]$ are in the same $G_{\rho}$-orbit.

Now consider elements $[\alpha_{1},0]$, 
$[\alpha_{2},0]$ which belong to $G_{\rho}$-orbits of type (IV). As 
usual we may assume  $\alpha_{i}=\iota([\alpha_{i}])$ for $i=1,2$. 
We set $r_0:=|\alpha_{2}(X,Y)|/|\alpha_{1}(X,Y)|= 
|\alpha_{2}(X,Z)|/|\alpha_{1}(X,Z)|$ and 
define as above a map $U_{0}\in O(n)$ commuting with $\rho(X)$ such that 
$ U_{0}(\alpha_{2}(X,Y))=r_{0}\alpha_{1}(X,Y)$ and 
$U_{0}(\alpha_{2}(X,Z))=r_{0}\alpha_{1}(X,Z).$
If we choose $S=S(1,0,0,\sqrt{r_{0}},0)$, then $(S,U_{0})\in G_{\rho}$ 
and $\iota_{Q}((S,U_{0})^* \iota_{Q}^{-1}(\alpha_{1})) =\alpha_{2}$. 
Hence $[\alpha_{1},0]$ and $[\alpha_{2},0]$ are in the same $G_{\rho}$-orbit.
\qed

\subsection{The case $\fl=\fr_{3,-1}$}

Now we suppose $\fl=\fr_{3,-1}$.
Let again $\lambda^{1}$ and $\lambda^{2}$ be the weights of the representation 
$\adR$ of $\fl_{\Bbb C}/R_{\Bbb C}$ on $R_{\Bbb C}$. They are given 
by $\lambda^{1}(X)=1$ and $\lambda^{2}(X)=-1$.  In particular both weights are 
real. Since $\fa$ is Euclidean this implies 
$E_{\lambda^{1}}=E_{\lambda^{2}}=0$. Moreover, we have 
$E_{\lambda^{1}+\lambda^{2}}=\fa^{\fl}$. Suppose $\rho\in\Hom(\fl,\so(\fa,\ip_{\fa}))_{\rm 
ss}$. By Lemma
\ref{maus} the map
\begin{eqnarray}
   &\{\alpha\in C^2(\fl,\fa)\mid\alpha(Y,Z)\in\fa^\fl,
 \alpha(X,\cdot)=0\}  \longrightarrow H^2(\fl,\fa)&\\
 &\alpha \longmapsto [\alpha] &\nonumber
\end{eqnarray}
 is an isomorphism.  Let $\iota$ be the inverse of this isomorphism. Obviously
 \begin{eqnarray*}
 &\iota_1:\quad \{\alpha\in C^2(\fl,\fa)\mid\alpha(Y,Z)\in\fa^\fl,
 \alpha(X,\cdot)=0\}\longrightarrow \fa^\fl&\\
&\alpha\longmapsto\alpha(Y,Z)\,.&
\end{eqnarray*}
 is an isomorphism. Since $d\tau(Y,Z)=0$ for all $\tau\in C^1(\fl,\fa)$ the map
\begin{eqnarray*}
\iota_0:\quad H^2(\fl,\fa)&\longrightarrow& \fa^\fl\\
 \,[\alpha] &\longmapsto &\alpha(Y,Z) \nonumber
\end{eqnarray*}
is well-defined. We have $\iota_1^{-1}\circ\iota_0\circ\iota^{-1}=\Id$, thus
$\iota_0=\iota_1\circ\iota.$ In particular, $\iota_0$ is an isomorphism.

\begin{lm}If $\fa$ is a Euclidean orthogonal $\fl$-module, then we have
$$\cH^2_Q(\fl,\fa)_\sharp=  \cH^2_Q(\fl,\fa)\setminus 
\{[0,0]\}=\iota_Q^{-1}(\,(H^2(\fl,\fa)\setminus{0})\cup(C^3(\fl)\setminus{0})\,).$$
\end{lm}
\proof We have to check which cohomology classes $[\alpha,\gamma]\in 
\cH^2_Q(\fl,\fa)$ are admissible. We assume 
$\alpha=\iota([\alpha])$.  As in the case of $\fl=\fn(2)$ all
conditions but $(A_{1})$ are trivially satisfied. Note that here
$(A_{1})(i)$ is equivalent to $\alpha(Y,Z)=0$ and thus to $\alpha=0$. 
However in case 
$\alpha=0$ Condition $(A_1)(ii)$ is satisfied 
if and only if
$$\gamma(X,Y,Z)=\langle \Phi_2(Y),[X,Z]\rangle+\langle \Phi_2([X,Y]),Z\rangle
=  \langle \Phi_2(Y),-Z\rangle+\langle \Phi_2(Y),Z\rangle =0.$$ Hence 
$(A_{1})$ holds for all $[\alpha,\gamma]\not=[0,0]$. Obviously 
$[0,0]$ is not admissible.
\qed

Next we will describe the automorphism group of $\fl=\fr_{3,-1}$.
Let $v,w,a,b,c,d\in\RR$, such that $ad\not=0$ and $bc\not=0$. We define automorphisms 
$S'(v,w,a,d)$ and $S''(v,w,b,c)$ of $\fl$
by the following matrices with respect to the basis $X,\,Y,\,Z$ of $\fl$:
$$S'(v,w,a,d)=
    \small{
    \left( \begin{array}{ccc}
    1&0&0\\
    v&a&0\\
    w&0&d
    \end{array}\right)}\,,\quad
    S''(v,w,b,c)=
    \small{
    \left( \begin{array}{ccc}
    -1&0&0\\
    v&0&b\\
    w&c&0
    \end{array}\right)}.$$
Then we have
  $$  \Aut(\fl)=\{S'(v,w,a,d)\mid v,w,a,d\in\RR,\  ad\not=0\}\ \cup \
      \{S''(v,w,b,c) \mid v,w,b,c\in\RR,\ bc\not=0 \}. $$

Since $\fl_{0}=\fl/R(\fl)=\RR\cdot
X$ we identify again $\lambda\in 
(\fl_{0}^{*})^{m}$ with $\lambda(X)\in\RR^{m}$.

\begin{pr} \label{Pr}
If $(\fg,\ip)$ is an indecomposable metric Lie algebra of index 3 such that
$\fg/\fj(\fg)^{\perp}\cong\fl:=\fr_{3,-1}$, then $\fg$
    is isomorphic to exactly one of the following indecomposable metric Lie algebras $\dd$
    with
    \begin{itemize}
        \item[(I)] $\fa=\RR^{2m+1}=\RR^{2m}\oplus\RR^{1}$,\ $m\ge0$, \
        $\rho=\rho^+_{\lambda}\oplus\rho_{0}$,\\ where
        $0<\lambda^{1}\le\lambda^{2}\le \ldots \le\lambda^{m}$,\\
        $\alpha(Y,Z)=A_{1}^{0}$, $\alpha(X,Y)=\alpha(X,Z)=0$,\\
        $\gamma=0$,
        \item[(II)]$\fa=\RR^{2m}$,\ $m\ge0$,  \
        $\rho=\rho^+_{\lambda}$,\\ where
        $0<\lambda^{1}\le\lambda^{2}\le \ldots \le \lambda^{m}$,\\
        $\alpha=0$, \\
        $\gamma(X,Y,Z)=1$.
\end{itemize}
\end{pr}
\proof Again $(\fg,\ip)$ is a quadratic extension of $\fl$ by $\fa=\RR^{n}$. By the same 
reasons as in the proof of
Proposition \ref{Pn2} each element in $\Hom(\fl,\so(n))/G$ has a canonical 
representative $\rho^+_{\lambda}\oplus\rho_{0}$, where
$\lambda\in\RR^{m}$, $2m\le n$, $0<\lambda^{1}\le\lambda^{2}\le \ldots 
\le\lambda^{m}$ and
$\rho_{0}$ is the trivial representation on $\RR^{n-2m}$. We fix such a 
representation $\rho=\rho^+_{\lambda}\oplus\rho_{0}$.
As in the case of $\fn(2)$ a cohomology class $[\alpha,\gamma]\in
\cH^2_Q(\fl,\fa_{\rho})_{\sharp}$ is indecomposable if and only if
$\fa_{\rho}^{\fl}=\RR\cdot \alpha(Y,Z)$. In particular, we may assume
$n=2m$ or $n=2m+1$ since otherwise  $\cH^2_Q(\fl,\fa_{\rho})_0$ is empty.

Let $S:=S'(0,0,a,d)$ be as
above. Then we have $(S,\Id)\in G_{\rho}$ and
\begin{eqnarray*}
\iota\circ\iota_{Q}((S,\Id)^* \iota_{Q}^{-1}([\alpha])) 
&=&\iota([S^*\alpha])=\iota_1^{-1}\circ\iota_0([S^*\alpha)])
=\iota_1^{-1}(ad\,\alpha(Y,Z))\\
\iota_{Q}((S,\Id)^*
\iota_{Q}^{-1}(\gamma))&=&ad\cdot\gamma.
\end{eqnarray*}
Therefore $G_\rho$ acts transitively on
$\iota_Q^{-1}(H^2(\fl,\fa_\rho)\setminus{0})$ and on $\iota_Q^{-1}(C^3(\fl)\setminus{0})$. It follows that we can represent each element in $\cH^2_Q(\fl,\fa_\rho)_0/G_\rho$ as in {\it(I)} or {\it (II)}.
\qed

\subsection{The case $\fl=\fh(1)$}
We consider now $\fl=\fh(1)$. Here we have $R=\RR\cdot Z$.
\begin{lm}
For $\rho\in\Hom(\fl,\so(\fa,\ip_{\fa}))_{\rm ss}$ the map
\begin{eqnarray*} & \{\alpha\in C^2(\fl,\fa)\mid\alpha(X,Y)=0,\ 
\alpha(Z,\fl)\subset\fa^\fl\}
    {\longrightarrow}H^2(\fl,\fa)&\\
    &\alpha\longmapsto [\alpha]&
\end{eqnarray*}
is an isomorphism.
\end{lm}
\proof The map is well defined and injective. We will prove that it is surjective. Let $a\in 
H^2(\fl,\fa)$ be given. By Lemma
\ref{fliege} we know that $a=[\alpha]$ with $\alpha(Z,\fl)\subset\fa^\fl$.
We write $\alpha=\alpha_{1}+\alpha_{2}$ with $\alpha_1(X,Y)=0$ and 
$\alpha_{2}(Z,\cdot)=0$. Then $\alpha_{2}$ belongs to the subcomplex 
of $C^{*}(\fl,\fa)$ which consists of all cochains $\sigma\in C^{*}(\fl,\fa)$ which
satisfy $\sigma(Z,\cdot)=0$. This subcomplex is equivalent to
$C^{*}(\fl/R,\fa)$. Since $\fl/R$ is abelian we obtain
$[\alpha_{2}]=[\tilde\alpha_{2}]$ for a cocycle $\tilde\alpha_{2}\in 
Z^2(\fl,\fa)$ satisfying $\tilde\alpha_{2}(X,Y)=:A\in\fa^{\fl}$ and 
$\tilde\alpha_{2}(Z,\cdot)=0$. Now we define a cocycle $\tau\in 
Z^1(\fl,\fa)$ by $\tau(X)=\tau(Y)=0$, $\tau(Z)=A$. Then 
$\tilde\alpha_{2}+d\tau=0$. Hence, $[\alpha]=[\alpha_{1}]$ and we 
obtain $[\alpha_{1}]=a$.
\qed

Again the lemma can be considered as a consequence of the Hochschild-Serre spectral sequence 
associated with $R\subset\fl$ (compare Remark \ref{HS}).

Let $\iota$ denote the inverse of the isomorphism defined in the above lemma. 
Obviously
\begin{eqnarray*}
 &\iota_1:\{\alpha\in C^2(\fl,\fa)\mid\alpha(X,Y)=0,\ 
\alpha(Z,\fl)\subset\fa^\fl\}\ \longrightarrow (\fl/R)^*\otimes\fa^\fl&\\
&\alpha\longmapsto \alpha(Z,\cdot) &
\end{eqnarray*}
 is an isomorphism. Since 
$\proj_{\fa^\fl}(d\tau(Z,\cdot))=-\proj_{\fa^\fl}(\rho(\cdot)\tau(Z))=0$ for all
$\tau\in C^1(\fl,\fa)$ the map
\begin{eqnarray*}
\iota_0:\quad H^2(\fl,\fa)&\longrightarrow& (\fl/R)^*\otimes\fa^\fl\\
 \,[\alpha] &\longmapsto &\proj_{\fa^\fl}(\alpha(Z,\cdot)) \nonumber
\end{eqnarray*}
is well-defined. We have $\iota_1^{-1}\circ\iota_0\circ\iota^{-1}=\Id$, thus
$\iota_0=\iota_1\circ\iota.$ In particular, $\iota_0$ is an isomorphism.

\begin{lm}If $\fa$ is a Euclidean orthogonal $\fl$-module, then we have
$$\cH^2_Q(\fl,\fa)_\sharp= \{ [\alpha,\gamma]\in \cH^2_Q(\fl,\fa)\mid 
[\alpha]\not=0\}=\iota_Q^{-1}(H^2(\fl,\fa)\setminus{0}).$$
\end{lm}
\proof
All cohomology classes $[0,\gamma]\in \cH^2_Q(\fl,\fa)$ are not admissible. Neither
$(A_{0})$ nor $(A_{1})$ is satisfied.  For instance, if we assume 
that  $(A_{0})$ holds and if we consider 
$L_{0}=Z$, $A_{0}=0$, and $Z_{0}\in\fl^{*}$ defined by 
$Z_{0}(X)=Z_{0}(Y)=0$, $Z_{0}(Z)=-\gamma(X,Y,Z)$ we get a contradiction.
If $[\alpha]\not=0$ neither assumption {\it (i)} of $(A_{0})$ nor  assumption 
{\it (i)} of $(A_{1})$ is satisfied. Hence $(A_{0})$ and $(A_{1})$
hold. Since also $(B_{0})$ and $(B_{1})$ hold (because $\fa$ is 
Euclidean) $[\alpha,\gamma]$ is admissible for $[\alpha]\not=0$. 
\qed

We describe automorphisms of $\fl$ by matrices with respect to the 
basis $X,\,Y,\,Z$ of $\fl$. The automorphism group of $\fl=\fh(1)$ equals
$$\Aut(\fl)=\Big\{
    S(A,u,x)=
\small{
\left( \begin{array}{cc}
A&0\\
x^{\top}&u
\end{array}\right)
}
\ \Big|\ A\in GL(2,\RR),\ \det A=u,\
x\in\RR^{2}
\Big\}.$$

\begin{pr}\label{Ph1}
    Let $\fg$ be an indecomposable metric Lie algebra of index 3 such that
    $\fg/\fj(\fg)\cong\fl:=\fh(1)$. Then $\fg$
        is isomorphic to one of the indecomposable metric Lie algebras
        $\fd_{\alpha,0}(\fl,\fa,\rho)$ with
    \begin{itemize}
        \item[(I)]
        $\fa=\RR^{2m+1}=\RR^{2m}\oplus\RR^1$, $m\ge0$, 
        $\rho=\rho^+_{\lambda}\oplus\rho_0$, $\lambda\in ((\fl/R)^{*}\setminus 
        0)^{m}$\\
        $\alpha(X,Y)=0,\ \alpha(X,Z)=A_{1}^0,\ \alpha(Y,Z)=0$.

        Two such Lie algebras for $\lambda\in ((\fl/R)^*\setminus 0)^m$ and
        $\bar \lambda\in ((\fl/R)^*\setminus 0)^{\bar m}$ are isomorphic if and 
        only if $m=\bar m$ and either
        \begin{itemize}
        \item[(a)]
        both $\Span \{\lambda(X),\lambda(Y)\}$ and 
        $\Span\{\bar\lambda(X),\bar\lambda(Y)\}$
        are one-dimensional and
        $$(\Span \{\lambda(X),\lambda(Y)\},\RR\cdot\lambda(Y)) =
          (\Span\{\bar\lambda(X),\bar\lambda(Y)\},\RR\cdot\bar\lambda(Y))\mod\frak 
        S_{m}\ltimes
        (\ZZ_{2})^{m}$$ or
        \item[(b)]
        both $\Span \{\lambda(X),\lambda(Y)\}$ and 
        $\Span\{\bar\lambda(X),\bar\lambda(Y)\}$
        are two-dimensional and
        $$(\bar\lambda(X)\wedge\bar\lambda(Y),\bar\lambda(Y))=(r\lambda(X)\wedge
        \lambda(Y),r^2\lambda(Y))\mod\frak S_{m}\ltimes 
        (\ZZ_{2})^{m}$$
        for a real number $r\not=0$.
        \end{itemize}
        \item[(II)]
        $\fa=\RR^{2m+2}=\RR^{2m}\oplus\RR^2$, $m\ge 0$, 
        $\rho=\rho^+_{\lambda}\oplus\rho_0$, $\lambda\in ((\fl/R)^{*}\setminus 
        0)^{m}$,\\
        $\alpha(X,Y)=0,\ \alpha(X,Z)=A_{1}^0,\ \alpha(Y,Z)=A_{2}^0$.\\
        Two such Lie algebras are isomorphic if and only if
        the $(m\times 2)$-matrices $M_\lambda:=(\lambda(X),\lambda(Y))$,
        $M_{\bar\lambda}:=(\bar\lambda(X),\bar\lambda(Y))$    satisfy
        $$ M_\lambda M_\lambda^\top = M_{\bar\lambda}M_{\bar\lambda}^\top \ \mod\
\frak S_{m}\ltimes
        (\ZZ_{2})^{m}.$$
\end{itemize}
\end{pr}
\proof We know that $(\fg,\ip)$ is a quadratic extension of $\fl$ by $\fa=\RR^n$. We first 
describe the action of $\Aut(\fl)$ on
$$\coprod_{\rho\in \Hom(\fl,\so(n))}\cH^2_Q(\fl,\fa_\rho)_\sharp=\coprod_{\rho\in 
\Hom(\fl,\so(n))}\iota_Q^{-1}(H^2(\fl,\fa_\rho)\setminus{0}).$$
Suppose $\rho \in \Hom(\fl,\so(n))$, 
$[\alpha,0]\in\iota_Q^{-1}(H^2(\fl,\fa_\rho)\setminus{0})$ and assume 
$\alpha=\iota([\alpha])$. Then we have for $(S,U)\in G$
\begin{eqnarray*}
\iota\circ\iota_Q((S,U)^*[\alpha,0])&=&\iota ([U_*^{-1}S^*\alpha])\
=\ \iota_1^{-1}\circ \iota_0 
([U_*^{-1}S^*\alpha])\in\iota(H^2(\fl,\fa_{\rho'})\setminus \{0\}),
\end{eqnarray*}
where $\rho'=(S,U)^*\rho$.
Since $\alpha=\iota([\alpha])$ we obtain for $S=S(A,u,x)$
\begin{eqnarray}
&\iota\circ\iota_Q (S,U)^*([\alpha,0])=\tilde\alpha\in H^2(\fl,\fa_{\rho'})&\nonumber\\
&\tilde\alpha(X,Y)=0,\quad  (\tilde\alpha(X,Z),\,\tilde\alpha(Y,Z))=
u\cdot(U^{-1}\alpha(X,Z),\,U^{-1}\alpha(Y,Z))\cdot A.& \label{cafe}
\end{eqnarray}

Now suppose that $\rho\in\Hom(\fl,\so(n))$ is given such that 
$\cH^2_Q(\fl,\fa_\rho)_\sharp\not=\emptyset$. By Proposition \ref{PHom} we may 
assume that $\rho=\rho^+_\lambda\oplus\rho_0$ for some $\lambda \in 
((\fl/R)^*\setminus 0)^m$, $2m\le n$. Now we consider the $G$-orbit through an 
element $[\alpha,0] \in \cH^2_Q(\fl,\fa_\rho)_\sharp$ . We assume 
$\alpha=\iota([\alpha])$, i.e. $\alpha(X,Y)=0$ and 
$\alpha(Z,\fl)\subset\fa_\rho^\fl$. Since $\fl$ does not decompose into the direct 
sum of two non-trivial Lie algebras $[\alpha,0]$ is indecomposable if and only if 
$\alpha(Z,\fl)=\fa_\rho^\fl$. Hence, we may identify $\fa_\rho^\fl$ with $\RR^1$ or 
$\RR^2$ spanned by the orthonormal basis $A^0_1$ or $A^0_1$, $A^0_2$, respectively.
We may choose a map $A\in \GL(2,\RR)$ such that
$$(\alpha(X,Z),\,\alpha(Y,Z))\cdot A=\left\{
\begin{array}{lll}
(A^0_1,0)& \mbox{ if }  &\dim \alpha(Z,\fl)=1 \\
(A^0_1,A^0_2)& \mbox{ if }  &\dim \alpha(Z,\fl)=2
\end{array}
\right.
$$
We set $u:=(\det A)^{1/3}$ and $S=S(u^{-1}A,u,0)$. Then $S^*\rho$ and
$S^*\alpha$ satisfy the conditions in $(I)$ if $\dim \alpha(Z,\fl)=1$
or in $(II)$ if $\dim \alpha(Z,\fl)=2$.

Now consider two representations $\rho=\rho^+_{\lambda}\oplus\rho_0$ and
$\bar\rho=\rho^+_{\bar\lambda}\oplus\rho_0$ for
$\lambda\in ((\fl/R)^*\setminus 0)^m$ and $\bar \lambda\in ((\fl/R)^*\setminus
0)^{\bar m}$. Let the 2-form $\alpha$ be defined by $\alpha(X,Y)=0,\
\alpha(X,Z)=A_{1}^0,\ \alpha(Y,Z)=0$. Then $a=[\alpha]\in H^2(\fl,\fa_\rho)$ and
$\bar a=[\alpha]\in H^2(\fl,\fa_{\bar \rho})$. When $a$ and $\bar a$ are in the
same $G$-orbit? We have to check under which conditions there is an element 
$(S,U)\in G$ such that
$(S,U)^*\rho=\bar\rho$ and $(S,U)^* a=\bar a$, which is equivalent to
$\iota\circ\iota_Q (S,U)^*\iota_Q^{-1}([\alpha])=\alpha$. By (\ref{cafe}) we can find
such an element $(S,U)\in G$ if and only if
$m=\bar m$ and there are maps $U_0\in O(1)=\pm1$ and $A\in \GL(2,\RR)$ such that
$$(\lambda(X),\lambda(Y)) \cdot A=
(\bar\lambda(X),\bar\lambda(Y))\ \mod \frak S_{m}\ltimes
(\ZZ_{2})^{m},$$ and
$$\det A\cdot(U_0^{-1} A_1^0,0)\cdot A=(A_1^0,0).$$
The last equation is satisfied if and only if
$A=\left( \begin{array}{cc} \delta&0\\c&d\end{array}\right)$ with
$\delta,c,d\in\RR$ and $U_0(A_1^0)=\delta^2d \cdot A_1^0$. In particular, $\delta^2d=\pm1$.
Hence, $a$ and $\bar a$ are in the same $G$-orbit if and only if there are real 
numbers $c,\delta\in\RR$, $\delta\not=0$, such that
$$(\delta\lambda(X)+c\cdot\lambda(Y),\pm\frac1{\delta^2}\cdot\lambda(Y)) =
(\bar\lambda(X),\bar\lambda(Y))\ \mod \frak S_{m}\ltimes
(\ZZ_{2})^{m}.$$ This proves the isomorphy condition in $(I)$.

For metric Lie algebras of type (II) we proceed in a similar way. Here we obtain 
that
two such Lie algebras for
$\lambda\in ((\fl/R)^*\setminus 0)^m$ and $\bar \lambda\in ((\fl/R)^*\setminus 
0)^{\bar m}$ are isomorphic if and only if $m=\bar m$ and
        $$(\lambda(X),\lambda(Y)) =
        (\bar\lambda(X),\bar\lambda(Y))\
        \mod O(2)\times\frak S_{m}\ltimes
        (\ZZ_{2})^{m}.$$

Now we use the following general fact, which is true for arbitrary
homomorphisms $A,\bar A:\RR^N\rightarrow \RR^m$:
$$\exists U\in O(N): \bar A=A\cdot U \quad \Leftrightarrow \quad A\, A^*=\bar A \, 
\bar
A^*:\RR^m\longrightarrow \RR^m. $$
Applying this for $\lambda,\bar\lambda:\fl/R\rightarrow \RR^m$,
where we identify $\fl/R$ with $\RR^2$ using the basis $X+R,\,Y+R$, we obtain the 
isomorphy condition in $(II)$.
\qed
\subsection{The case $\fl=\fsl(2,\RR)$}

\begin{pr}\label{Psl}
    Let $(\fg,\ip)$ be an indecomposable metric Lie algebra of index 3 such that
    $\fg/\fj(\fg)\cong\fsl(2,\RR)$.
    Then $(\fg,\ip)$ is isomorphic to exactly one of the indecomposable metric Lie
    algebras $(\fsl(2,\RR)\ltimes\fsl(2,\RR)^{*},\ip_{c})$, $c\in\RR$,
    where $\ip_{c}$ is defined by
    $$\langle L_1+Z_1,L_2+Z_2 \rangle = Z_1(L_2)+Z_2(L_1)+ c B_{\fl}(L_1,L_2),$$
    for all $L_1,L_2\in \fsl(2,\RR)$ and $Z_1,Z_2\in \fsl(2,\RR)^*$. Here 
    $B_{\fl}$ denotes
    the Killing form of $\fl:=\fsl(2,\RR)$.
\end{pr}
\proof Let $\fa$ be an orthogonal $\fl$-module. Since $\fl$ is semi-simple (and, in particular, 
unimodular) we have
$H^2(\fl,\fa)=0$ and $H^3(\fl)=C^3(\fl)$. Therefore
$$\cH^2_Q(\fl,\fa)_\sharp=\cH^2_Q(\fl,\fa)=C^3(\fl)$$ by Proposition \ref{club}. In particular, 
$\cH^2_Q(\fl,\fa)_0=C^3(\fl)$ if $\fa^\fl=0$ and $\cH^2_Q(\fl,\fa)_0=\emptyset$ if 
$\fa^\fl\not=0$. Since $C^3(\fl)$ is
one-dimensional any $\gamma\in C^3(\fl)$ is a multiple of the non-vanishing 3-form
$B_\fl(\lb,\cdot)$, which is $\Aut(\fl)$-invariant.
On the other hand each orthogonal representation of $\fl$ on a Euclidean space
is trivial.
Therefore, by Theorem \ref{Pwieder}, the  metric Lie algebra $(\fg,\ip)$ is isomorphic to  exactly one of the balanced
quadratic extensions $\fd_{0,\gamma}(\fl,0,0)$ for $\gamma\in C^3(\fl)$.   Now it follows from Remark \ref{dprime} that $(\fg,\ip)$ is isomorphic to
$\fd'(\fl,cB_\fl,0,0)$ for a unique $c\in\RR$.
\qed
\subsection{The case $\fl=\fsu(2)$}
For $m\in\ZZ$, $m\ge 0$ we define a set $K_m$ by
$$K_m=\left\{k=((k^1,\ldots,k^r),(k_1,\ldots,k_s))\left| \begin{array}{l}
k_i,k^i\in\NN, \\
0<k^1\le\ldots\le k^r,\ 0<k_1\le\ldots\le k_s,\\
\textstyle{
\sum_{i=1}^r (2k^i+1) +  \sum_{i=1}^s 4k_i =m
}
\end{array}
 \right.\right\}\,.$$
Here we allow that $r=0$ or $s=0$, e.g. $K_0=\{(\emptyset,\emptyset)\}$.

For $k\in\NN$ let $\sigma_k:\fsu(2)\rightarrow\so(2k+1)$ and
$\sigma'_k:\fsu(2)\rightarrow\so(4k)$ be the non-trivial irreducible real
representations of $\fsu(2)$. For $k=((k^1,\ldots,k^r),(k_1,\ldots,k_s))\in K_m$
let the representation $\rho_k$ of $\fsu(2)$ on $\RR^m$ be the direct sum
$$\rho_k=\bigoplus_{i=1}^r \sigma_{k^i}\oplus \bigoplus_{i=1}^s\sigma'_{k_i}.$$
In particular, if $k=(\emptyset,\emptyset)\in K_0$, then $\rho_k$ is the zero representation.

Let $B_\fl$ denote the Killing form of $\fl=\fsu(2)$. In the following proposition we will use 
the modified quadratic extensions defined in Remark \ref{dprime}.
\begin{pr}\label{Psu}
    If $(\fg,\ip)$ is an indecomposable metric Lie algebra of index 3 such that
    $\fg/\fj(\fg)\cong\fl=\fsu(2)$, then $(\fg,\ip)$ is isomorphic to exactly one of 
    the metric Lie algebras $\fd'_{0,0}(\fl,cB_\fl,\RR^m,\rho_k),$  where $m>0$, 
    $k\in K_m$ and $c\in\RR$.
   
    In particular, if $\dim \fg =6$, then $(\fg,\ip)$ is
    isomorphic to exactly one of the
    metric Lie algebras $\fd'_{0,0}(\fl,cB_\fl,0,0)=\fsu(2)\ltimes\fsu(2)^*$ for
    $c\in\RR$.
\end{pr}
\proof As in the case of $\fsl(2,\RR)$ we have $\cH^2_Q(\fl,\fa)_0=C^3(\fl)$ if $\fa^\fl=0$ and 
$\cH^2_Q(\fl,\fa)_0=\emptyset$ if $\fa^\fl\not=0$. Each orthogonal representation $(\rho,\fa)$ 
of $\fl$ with $\fa^\fl=0$ is equivalent to a unique representation $\rho_k$, $k\in K_m$, 
$m\ge0$.

Furthermore, we have $C^3(\fl)=\RR\cdot\gamma_0$ for $\gamma_0:=B_\fl(\lb,\cdot)$. As in the 
case of $\fsl(2,\RR)$ we can now use Theorem \ref{Pwieder} and Remark \ref{dprime} to finish the proof.
\qed

\subsection{The case $\fl=\RR^k$, $k=1,2$}

If $\fl=\RR^k$, $k=1,2$, then we can identify
$$ \cH^2_Q(\fl,\fa)=H^2(\fl,\fa)=C^2(\fl,\fa^\fl).$$
\begin{lm} For $\fl=\RR^k$, $k=1,2$ we have
    $$\cH^2_Q(\fl,\fa)_0=\{\alpha\in C^2(\fl,\fa^\fl)\mid \alpha(\fl,\fl)\subset
\fa^\fl\ \mbox{\rm non-degenerate}, \
      \alpha \ \mbox{\rm indecomposable}\}.$$
\end{lm}
\proof Condition $(B_0)$ implies that $\cH^2_Q(\fl,\fa)_0$ is contained in the set 
on the r.\,h.\,s.
Now let $\alpha\in C^2(\fl,\fa^\fl)$ be indecomposable and such that 
$\alpha(\fl,\fl)$ is non-degenerate.
Assume that $\alpha$ is not admissible. Then Condition $(A_0)$ cannot be satisfied.
Hence there are elements $L_0\in\Ker \rho$, $L_0\not=0$ and $A_0\in\fa$ such that:
\begin{eqnarray}
    \alpha(\cdot,L_0)=\rho(\cdot)(A_0) \label{i}.
\end{eqnarray}
Since $\alpha(\fl,L_0)\subset\fa^\fl$ Equation (\ref{i}) implies  $\alpha(L_0,\cdot)=0$. Since 
on the other hand
$L_0\in\Ker\rho$ and $L_0\not=0$ this is a contradiction to the indecomposability
of $\alpha$.
\qed

A pair $(\fl,\fa)$ is called decomposable if it is a non-trivial direct sum of two pairs. Otherwise the pair is called indecomposable.
\begin{co}\label{co}
    For $\fl=\RR^k$, $k=1,2$ we have
    $$\cH^2_Q(\fl,\fa)_0=
    \left\{
        \begin{array}{lll}
            \emptyset &\mbox{ if }& \dim \fa^\fl>1\\
            C^2(\fl,\fa^\fl)\setminus \{0\} &\mbox{ if }& \dim \fa^\fl=1\\
            \{0\} &\mbox{ if }& \dim \fa^\fl=0, \ (\fl,\fa) \mbox{ is indecomposable}\\
            \emptyset &\mbox{ if }&\dim \fa^\fl=0, \ (\fl,\fa) \mbox{ is decomposable}
        \end{array}
    \right..$$
\end{co}

First we consider the case $\fl=\RR^{2}$. We fix a basis $\{Y,Z\}$ of $\RR^{2}$.
\begin{pr} \label{PR2}
If $\fg$ is an indecomposable metric Lie algebra of index 3 such that
$\fg/\fj(\fg)\cong\fl:=\RR^{2}$, then $\fg$
    is isomorphic to one of the following indecomposable Lie algebras
$\fd_{\alpha,0}(\fl,\fa,\rho)$
    with
    \begin{itemize}
        \item[(I)] $\fa=\RR^{1,2m}=\RR^{2m}\oplus\RR^{1,0},\ m\ge 0,\
        \rho=\rho^+_{\lambda}\oplus\rho_{0},$ where $\lambda\in( \fl^{*}\setminus
        0)^{m}$,\\
        $\alpha(Y,Z)=A_{1}^{0}$.\\
        Two such Lie algebras for $\lambda\in ((\fl/R)^*\setminus 0)^m$ and
        $\bar \lambda\in ((\fl/R)^*\setminus 0)^{\bar m}$ are isomorphic if and
        only if $m=\bar m$ and
        \begin{eqnarray*}
        &(\,\Span \{\lambda(Y),\lambda(Z)\},\lambda(Y)\wedge\lambda(Z)\,)=
        (\,\Span
        \{\bar\lambda(Y),\bar\lambda(Z)\},
	\pm\bar\lambda(Y)\wedge\bar\lambda(Z)\,)&\\&
        \hspace{7cm}\mod \frak S_{m}\ltimes
        (\ZZ_{2})^{m}.&
        \end{eqnarray*}
        \item[(II)]
        $\fa=\RR^{1,2m+2}=\RR^{2m}\oplus\RR^{1,1}\oplus\RR^{1},\ m\ge0, \
        \rho=\rho^+_{\lambda}\oplus \rho_{\mu}'\oplus\rho_{0},$ \\
        where  $\lambda\in (\fl^{*}\setminus 0)^{m}$, $\mu\in\fl^{*}$ such
        that\\
        $\mu(Y)=1$, $\mu(Z)=0$,\\
        $\alpha(Y,Z)=A^{0}_{1}$.\\
        Two such Lie algebras for $\lambda\in ((\fl/R)^*\setminus 0)^m$ and
        $\bar \lambda\in ((\fl/R)^*\setminus 0)^{\bar m}$ are isomorphic if
	and only if $m=\bar m$ and
        either
        \begin{itemize}
        \item[(a)] $\Span \{\lambda(Y),\lambda(Z)\}$ and $\Span
\{\bar\lambda(Y),\bar\lambda(Z)\}$
        are one-dimensional, and
        $$\lambda(Z)\not=0  \mbox{ and }
        \lambda(Z)= \bar\lambda(Z)
        \ \mod \frak S_{m}\ltimes(\ZZ_{2})^{m}$$
        or
        $$ \lambda(Z)=\bar\lambda(Z)=0  \mbox{ and }
        \lambda(Y)= \bar\lambda(Y)
        \ \mod \frak S_{m}\ltimes(\ZZ_{2})^{m}$$
        or
        \item[(b)] $\Span \{\lambda(Y),\lambda(Z)\}$ and $\Span
\{\bar\lambda(Y),\bar\lambda(Z)\}$
        are two-dimensional and
        $$(\lambda(Y)\wedge\lambda(Z)\,,\lambda(Z)))=
        (\pm\bar\lambda(Y)\wedge\bar\lambda(Z)\,,\lambda(Z))
        \ \mod \frak S_{m}\ltimes
        (\ZZ_{2})^{m}.$$
        \end{itemize}
        \item[(III)] $\fa=\RR^{1,2m+1}=\RR^{2m}\oplus\RR^{1,1}, \ m\ge2,\
        \rho=\rho^+_{\lambda}\oplus \rho_{\mu}',$ \\
        where  $\lambda\in (\fl^{*}\setminus 0)^{m}$ is such
        that
        the set $\{(1,0)\}\cup\{(\lambda^{i}(Y),\lambda^{i}(Z))\mid i=1,\ldots,m\}$
        is not contained in the union of two 1-dimensional subspaces of
        $\RR^{2}$, \\$\mu\in\fl^{*}$ is given by
        $\mu(Y)=1$, $\mu(Z)=0$,\\
        $\alpha=0$.\\
        Two such Lie algebras for $\lambda\in ((\fl/R)^*\setminus 0)^m$ and
        $\bar \lambda\in ((\fl/R)^*\setminus 0)^{\bar m}$ are isomorphic if and
        only if $m=\bar m$ and
        \begin{eqnarray*}
        &&\lambda(Y)+\RR\cdot\lambda(Z) =
        \bar\lambda(Y)+\RR\cdot\bar\lambda(Z)
        \ \mod \frak S_{m}\ltimes
        (\ZZ_{2})^{m}.
        \end{eqnarray*}
\end{itemize}
\end{pr}
\proof
Let $(\rho,\fa)$ be such that $\cH^2_Q(\fl,\fa)_0\not=\emptyset$. By Corollary
\ref{co} either $\dim \fa^\fl=1$ or $\dim \fa^\fl=0$. If $\dim \fa^\fl=1$, then
either $\fa^\fl=\RR^{1,0}$ or $\fa^\fl=\RR^1$.

Let us first consider the case $\fa^\fl=\RR^{1,0}$. This will lead to Lie algebras
of type $(I)$ in the proposition. By Proposition \ref{PHom} we may assume
$\fa=\RR^{1,2m}$ and $\rho=\rho^+_\lambda\oplus \rho_0$, $\lambda\in(\fl^*\setminus
0)^m$. Suppose $\alpha\in C^2(\fl,\fa^\fl)\setminus \{0\}\, (=\cH^2_Q(\fl,\fa)_0)$.
Let $A_1^0$ be a fixed unit vector in $\fa^\fl=\RR^{1,0}$. It is easy to find a map
$S\in\Aut(\fl)=\GL(2,\RR)$ such that $S^*\alpha(Y,Z)=A^0_1$. Then $S^*\rho$,
$S^*\alpha$ satisfy the conditions in $(I)$. Now let
$\rho=\rho^+_\lambda\oplus\rho_0$ and $\bar\rho=\rho^+_{\bar\lambda}\oplus\rho_0$
be representations on $\RR^{2m}\oplus\RR^{1,0}$ for different
$\lambda,\bar\lambda\in(\fl^*\setminus 0)^m$ and let $\alpha$ be defined by
$\alpha(Y,Z)=A^0_1$. Then $\alpha$ defines cohomology classes $a\in
\cH^2_Q(\fl,\fa_\rho)$ and
$\bar a\in \cH^2_Q(\fl,\fa_{\bar \rho})$. We have to check under which conditions
$a=\bar a\ \mod\ G$ holds.
First we note that $(S,U)^*\rho=\bar\rho$ implies $U(\fa^\fl)=\fa^\fl$ (thus
$U|_{\fa^\fl}=\pm1$) and
$(S,U|_{(\fa^\fl)^\perp})^*\rho_\lambda=\rho_{\bar\lambda}$.
Using this it is easy to see that $a=\bar a\ \mod\ G$ if and only if there is a map
$S\in\GL(2,\RR)$ such that $\alpha(SY,SZ)=\pm\alpha(Y,Z)$ and
$[S^*\lambda]=[\bar\lambda]\in\Lambda_m$. This is the case if and only if there is
a map $S\in\SL^\pm(2,\RR):=\{S\in\GL(2,\RR)\mid \det S=\pm 1\}$ satisfying
$[S^*\lambda]=[\bar\lambda]\in\Lambda_m$,
which is equivalent to the condition in $(I)$.

Now let us suppose $\fa^\fl=\RR^1$. This will lead to Lie algebras of type $(II)$.
Here we may assume $\fa=\RR^{1,2m+2}$ and $\rho=\rho^+_\lambda\oplus
\rho'_\mu\oplus\rho_0$, $\lambda\in(\fl^*\setminus 0)^m$,
$\mu\in\fl^*\setminus 0$. If $\alpha\in C^2(\fl,\fa^\fl)\setminus
0(=\cH^2_Q(\fl,\fa)_0)$ and $A_1^0$ is a fixed unit vector in $\fa^\fl=\RR^1$, then
one can easily find a map $S\in\GL(2,\RR)$ such that $S^*\alpha(Y,Z)=A^0_1$,
$S^*\mu(Y)=1$ and $S^*\mu(Z)=0$. Then $S^*\rho$, $S^*\alpha$ satisfy the conditions
in $(II)$. Now let $\alpha$ and $\mu$ be as in $(II)$ and consider
$\rho=\rho^+_\lambda \oplus \rho'_\mu\oplus\rho_0$ and
$\bar\rho=\rho^+_{\bar\lambda} \oplus \rho'_\mu\oplus\rho_0$ for
$\lambda,\bar\lambda\in(\fl^*\setminus 0)^m$. One proves in a similar way as above
that $a=[\alpha]\in \cH^2_Q(\fl,\fa_\rho)$ and $\bar a=[\alpha]\in
\cH^2_Q(\fl,\fa_{\bar \rho})$ are in the same $G$-orbit if and only if there is a
map $S\in\GL(2,\RR)$ such that $\alpha(SY,SZ)=\pm\alpha(Y,Z)$, $S^*\mu=\pm\mu$ and
$[S^*\lambda]=[\bar\lambda]\in\Lambda_m$. This is the case if and only if there
exists a $c\in\RR$ such that
$$(\pm\lambda(Y)+c\lambda(Z), \lambda(Z)))=(\bar\lambda(Y),\bar\lambda(Z))\ \mod
\frak S_{m}\ltimes (\ZZ_{2})^{m}.$$
This yields the isomorphism condition in {\it (II)}.

Finally we consider the case $\fa^\fl=0$. Here we have $\alpha=0$ and we may assume
$\fa=\RR^{1,2m}$ and $\rho=\rho^+_\lambda\oplus \rho'_\mu$,
$\lambda\in(\fl^*\setminus 0)^m$, $\mu\in\fl^*\setminus 0$.

The indecomposability of $\fg$ is equivalent to the line condition
on $\lambda$.

We have $S^*\mu=\pm\mu$ and $[S^*\lambda]=[\bar\lambda]\in\Lambda_m$ if and only if
there exists $c,d\in\RR$ such that $d\not=0$ and
$$(\pm\lambda(Y)+c\lambda(Z), d\lambda(Z)))=(\bar\lambda(Y),\bar\lambda(Z))\ \mod
\frak S_{m}\ltimes (\ZZ_{2})^{m}.$$
\qed

If $\fl=\RR=\RR\cdot X$, then it is easy to prove the following classification result. We
identify $\lambda\in
(\fl_{0}^{*})^{m}$ with $\lambda(X)\in\RR^{m}$ and $\mu\in
(\fl_{0}^{*})^{r}$ with $\mu(X)\in\RR^{r}$.
\begin{pr}\label{PR1}
    If $\fg$ is an indecomposable metric Lie algebra of index 3 such that
    $\fg/\fri(\fg)^{\perp}\cong\fl:=\RR^{1}$, then $\fg$
    is isomorphic to exactly one of the following indecomposable Lie algebras
    $\fd_{0,0}(\fl,\fa,\rho)$
    with
    \begin{itemize}
        \item[(I)] $\fa=\RR^{2,2m+2}=\RR^{2m}\oplus\RR^{2,2}, \ m\ge0,\
        \rho=\rho^+_{\lambda}\oplus \rho_{(\mu,\nu)}'',$ \\where $\lambda\in
        \RR^{m}$,
        $0<\lambda^{1}\le\lambda^{2}\le \ldots  \le\lambda^{m}$, $\mu=1,\
\nu\in\RR,\ \nu\not=0$;\\
        \item[(II)] $\fa=\RR^{2,2m}=\RR^{2m}\oplus\RR^{2,0}, \ m\ge0,\
        \rho=\rho^+_{\lambda}\oplus\rho^-_{\lambda^0}$, \\ where
        $\lambda^0=1,\ \lambda=(\lambda^{1},\ldots,\lambda^{m})\in
        \RR^{m}$,\
        $0<\lambda^{1}\le\lambda^{2}\le \ldots  \le\lambda^{m}$;\\
        \item[(III)] $\fa=\RR^{2,2m+2}=\RR^{2m}\oplus \RR^{2,2}, \ m\ge0,\
        \rho=\rho^+_{\lambda}\oplus\rho'_{\mu}$, \\ where
        $\lambda=(\lambda^{1},\ldots,\lambda^{m})\in
        \RR^{m}$,\
        $0<\lambda^{1}\le\lambda^{2}\le \ldots  \le\lambda^{m}$,\\
        $\mu=(\mu^{1},\mu^{2})\in\RR^{2}$, $1=\mu^{1}\le\mu^{2}$.\\
\end{itemize}
\end{pr}

\subsection{Summary}
Propositions \ref{jonas} and \ref{Pn2} -- \ref{PR2} yield a classification of
indecomposable
non-simple metric Lie algebras of index 3 up to the case where
$\fg/\fj(\fg)$ is isomorphic to $\RR^{3}$.
In this case $\fa$ is Euclidean and $(\fg,\ip)$ is an
indecomposable metric Lie algebra of index 3 with maximal isotropic
centre. For a classification of these algebras see \cite{KO03}, Theorem 5.1.

It remains to determine all simple metric Lie algebras of index 3.
This is done by checking the list of all simple Lie algebras.

Finally we obtain the following classification result for metric Lie
algebras of index 3.

\begin{theo}\label{T71}
If $(\fg,\ip)$ is a simple metric Lie algebra of index 3, then $\fg$
is isomorphic to $\fsu(2)$ or $\fsl(3,\RR)$ and $\ip$ is a positive
multiple of the Killing form or it is isomorphic to $\fsl(2,\CC)$
and $\ip$ is a non-zero multiple of the Killing form.

If $(\fg,\ip)$ is a non-simple indecomposable metric Lie algebra of index 3, then
$(\fg,\ip)$ is isomorphic to exactly one Lie algebra $\dd$ with the following data:
\end{theo}
\newpage
\thispagestyle{empty}
\begin{sideways}
{\rm
\small{
\begin{tabular}{|c|c|c|l|c|l|l|}
\hline
&&&&&&\\[-2ex]
$\fl$ & $\fa$ & $\rho$ &$\alpha$ (characteristic & $\gamma$&  parameters &detailed 
de-\\
&&& property)&&&scription in\\[1ex]
\hline \hline
&&&&&&\\[-2.0ex]
$\fn(2)$ &$\RR^{2m}$& $\rho_\lambda^{+}$&$\alpha=0$ & $\gamma_\kappa\not=0$  &$m\ge0,\ 
\kappa=\pm1,\ [\lambda]\in \Lambda_m$ &Prop.\,\ref{Pn2} $(Ia,b)$\\[0.5ex]
\cline{2-7}&&&&&& \\[-2ex]
&$\RR^{2m}\oplus\RR^1$& $\rho_\lambda^{+}\oplus\rho_0$&$\alpha(\fl,\fl)=\RR^1$ & 
$\gamma=0$  & $m\ge0,\ [\lambda]\in \Lambda_m$ &Prop.\,\ref{Pn2} {\it(II)}\\[0.5ex]
\cline{2-7}&&&&&& \\[-2ex]
&$\RR^{2m+2}\oplus\RR^1$& 
$\rho_{\lambda'}^{+}\oplus\rho_0$&$\alpha_r(\fl,\fl)=\RR^2\oplus\RR^1$ & $\gamma=0$  
&$m\ge0,\ \lambda'=(\lambda,1),\ [\lambda]\in \Lambda_m$, $r\in\RR,r>0$ 
&Prop.\,\ref{Pn2} {\it(III)}\\[0.5ex]
\cline{2-7}&&&&&& \\[-2ex]
&$\RR^{2m+2}$& $\rho_{\lambda'}^{+}$&$\alpha(\fl,\fl)=\RR^2$ & $\gamma=0$  &$m\ge0,\ 
\lambda'=(\lambda,1),\ [\lambda]\in \Lambda_m$ &Prop.\,\ref{Pn2} {\it (IV)}\\[0.5ex]
\hline&&&&&& \\[-2ex]
$\fr_{3,-1}$&$\RR^{2m}\oplus\RR^1$&$\rho_\lambda^{+}\oplus\rho_0$&$\alpha(\fl,\fl)=\RR^
1$&$\gamma=0$&$m\ge0,\ [\lambda]\in \Lambda_m$&Prop.\,\ref{Pr} $(I)$ \\[0.5ex]
\cline{2-7}&&&&&& \\[-2ex]
&$\RR^{2m}$&$\rho_\lambda^{+}$&$\alpha=0$&$\gamma\not=0$&$m\ge0,\ [\lambda]\in 
\Lambda_m$& Prop.\,\ref{Pr} {\it(II)}\\[0.5ex]
\hline&&&&&& \\[-2ex]
$\fh(1)$&$\RR^{2m}\oplus\RR^1$&$\rho_\lambda^{+}\oplus\rho_0$&$\alpha(\fl,\fl)=\RR^1$ 
&$\gamma=0$ &$m\ge0,\ [\lambda]\in 
\Lambda_m/(\RR^*\ltimes\RR)$&Prop.\,\ref{Ph1} {\it (I)} \\[0.5ex]
\cline{2-7}&&&&&& \\[-2ex]
&$\RR^{2m}\oplus\RR^2$&$\rho_\lambda^{+}\oplus\rho_0$&$\alpha(\fl,\fl)=\RR^2$ &$\gamma=0$       
&$m\ge0,\ [\lambda]\in \Lambda_m/O(2)$&Prop.\,\ref{Ph1} {\it (II)} \\[0.5ex]
\hline&&&&&& \\[-2ex]
$\fsl(2,\RR)$&$0$&--&$\alpha=0$&$\gamma_c$&$c\in\RR$&Prop.\,\ref{Psl} 
\\[0.5ex]
\hline&&&&&& \\[-2ex]
$\fsu(2)$
&$\RR^{m}$&$\rho_k$&$\alpha=0$&$\gamma_c$&$m\ge0,\ k\in 
K_{m},\ c\in\RR$& Prop.\,\ref{Psu}\\[0.5ex]
\hline&&&&&& \\[-2ex]
$\RR^1$&$\RR^{2m}\oplus\RR^{2,2}$&$\rho_{\lambda}^{+}\oplus\rho''_{(1,\nu)}$&$\alpha=0$
&$\gamma=0$ &$m\ge0,\ [\lambda]\in \Lambda_m,\ \nu\in\RR\setminus\{0\}$&
Prop.\,\ref{PR1} {\it (I)} \\[0.5ex]
\cline{2-7}&&&&&& \\[-2ex]
&$\RR^{2m}\oplus\RR^{2,0}$&$\rho_{\lambda}^{+}\oplus\rho_{1}^{-}$&$\alpha=0$&$\gamma=0$
&$m\ge 0,\ [\lambda]\in \Lambda_m$&Prop.\,\ref{PR1} {\it (II)}  \\[0.5ex]
\cline{2-7}&&&&&& \\[-2ex]
&$\RR^{2m}\oplus\RR^{2,2}$&$\rho_{\lambda}^{+}\oplus\rho'_{(1,\mu)}$&$\alpha=0$&$\gamma=0$
&$m\ge 0,\ [\lambda]\in \Lambda_m,\ 
\mu\in\RR, \ \mu\ge 1$&Prop.\,\ref{PR1} {\it (III)}   \\[0.5ex]
\hline&&&&&& \\[-2ex]
$\RR^2$&$\RR^{2m}\oplus\RR^{1,0}$&$\rho_{\lambda}^{+}\oplus\rho_{0}$&
$\alpha(\fl,\fl)=\RR^{1,0}$&$\gamma=0$&$m\ge0,\ [\lambda]\in 
\Lambda_m/SL^\pm(2,\RR)$& Prop.\,\ref{PR2} {\it (I)}\\[0.5ex]
\cline{2-7}&&&&&& \\[-2ex]
&$\RR^{2m}\oplus\RR^{1,1}\oplus\RR^{1}$&
$\rho_{\lambda}^{+}\oplus\rho'_{\mu}\oplus\rho_{0}$&
$\alpha(\fl,\fl)=\RR^1$&$\gamma=0$&$m\ge0,\ [\lambda]\in 
\Lambda_m/(\ZZ_{2}\ltimes\RR)$&  Prop.\,\ref{PR2} {\it (II)}\\[0.5ex]
\cline{2-7}&&&&&& \\[-2ex]
&$\RR^{2m}\oplus\RR^{1,1}$&$\rho_{\lambda}^{+}\oplus\rho'_{\mu}$&$\alpha=0$&$\gamma=0$
&$m\ge 2,\ [\lambda]\in {\cal O},\ {\cal O} \subset\Lambda_m/(\RR^*\ltimes\RR)$ open
&  Prop.\,\ref{PR2} {\it (III)}\\[0.5ex]
\hline&&&&&& \\[-2ex]
$\RR^3$&$\RR^{2m}$&$\rho_{\lambda}^{+}$&$\alpha=0$&$\gamma=0$&$m\ge 4,\ [\lambda]\in 
{\cal O}_{1},\ {\cal O}_{1} \subset\Lambda_m/GL(3,\RR)$ open &\cite{KO03}, \\[0.5ex]
\cline{2-6}&&&&&& \\[-2ex]
&$\RR^{2m}$&$\rho_{\lambda}^{+}$&$\alpha=0$&$\gamma\not=0$&$m\ge 0,\ [\lambda]\in 
\Lambda_m/SL(3,\RR)$&Theorem 5.3 \\[0.5ex]
\cline{2-6}&&&&&& \\[-2ex]
&$\RR^{2m}\oplus\RR^{1}$&$\rho_{\lambda}^{+}\oplus\rho_{0}$&$\alpha(\fl,\fl)=\RR^1$&$\gamma=0$
&$m\ge 2,\ [\lambda]\in {\cal O}_{2},\ {\cal O}_{2} 
\subset\Lambda_m/((SL^\pm(2,\RR)\times \RR^*)\ltimes\RR^{2})$ open& \\[0.5ex]
\cline{2-6}&&&&&& \\[-2ex]
&$\RR^{2m}\oplus\RR^{2}$&$\rho_{\lambda}^{+}\oplus\rho_{0}$&$\alpha(\fl,\fl)=\RR^2$&$\gamma=0$
&$m\ge0,\ [\lambda]\in \Lambda_m/( CO(2)\ltimes
\RR^{2})$& \\[0.5ex]
\cline{2-6}&&&&&& \\[-2ex]
&$\RR^{2m}\oplus\RR^{3}$&$\rho_{\lambda}^{+}\oplus\rho_{0}$&$\alpha(\fl,\fl)=\RR^3$&
$\gamma=0$
&$m\ge0,\ [\lambda]\in \Lambda_m/O(3)$& \\[0.5ex]
\hline
\end{tabular}
}}
\end{sideways}

With some effort, it should be possible to classify metric Lie algebras of index $4, 5, \dots$ in the same way. But very soon the method will reach its limits. On the one hand, there is the problem
of classification of admissible Lie algebras. On the other hand, also the explicit determination of the
orbit spaces
$\cH^2_Q(\fl,\fa_\rho)_0/G_\rho$ sometimes
leads to unsolved classification problems in multilinear algebra
(an example is discussed in \cite{KO03}, Remark 5.3).

\vspace{0.8cm}
{\footnotesize

Ines Kath\\
Max-Planck-Institut f\"ur Mathematik in den Naturwissenschaften\\
Inselstra{\ss}e 22-26, D-04103 Leipzig, Germany\\
email: ikath@mis.mpg.de\\
\\
Martin Olbrich\\
Mathematisches Institut
der Georg-August-Universit\"at G\"ottingen\\
Bunsenstra{\ss}e 3-5, D-37073 G\"ottingen, Germany\\
email: olbrich@uni-math.gwdg.de\\
}


\begin{thebibliography}{MMMM}

\bibitem[BK\,02]{BK02} H.\,Baum, I.\,Kath, {\sl Doubly Extended Lie Groups -- Curvature, Holonomy
and Parallel Spinors.} To appear in J. Diff. Geom. Appl., see also Preprint arXiv:math.DG/0203189.

\bibitem[BB1]{BB1}
L.~Berard Bergery,
{\sl D\'ecomposition de Jordan-H\"older d'une repr\'esentation de dimension finie,
adapt\'ee \`a une forme r\'eflexive.} Handwritten notes.


\bibitem[BB2]{BB2}
L.~Berard Bergery,
{\sl Structure des espaces symetriques pseudo-riemanniens.}
Handwritten notes.

\bibitem[Bor\,97]{bordemann}
M.~Bordemann,
{\sl Nondegenerate invariant bilinear forms on nonassociative algebras.}
Acta Math. Univ. Comenian. {\bf 66} (1997), 151-201.

\bibitem[Bou\,71]{bourbaki}
N.~Bourbaki,
{\sl Groupes et alg\`ebres de Lie. Chap. I: Alg\`ebres de Lie.}
Hermann, Paris, 1971.

\bibitem[CP\,80]{CP80} M.~Cahen, M.~Parker, {\sl Pseudo-Riemannian symmetric spaces.} Mem. Amer.
Math. Soc. {\bf 24} (1980), no. 229.

\bibitem[Gr\,98]{grishkov}
A.~N.~Grishkov,
{\sl Orthogonal modules and nonlinear cohomologies.}
Algebra and Logic {\bf 37} (1998), 294--306.

\bibitem[KO\,02]{KO03} I.\,Kath, M.\,Olbrich,
{\sl Metric Lie algebras with maximal isotropic centre.} 
To appear in Math. Z., see also Preprint arXiv:math.DG/0209366, 2002.

\bibitem[M 85]{M85} A.~Medina, {\sl Groupes de Lie munis de m\'etriques bi-invariantes.} Tohoku
Math. J. (2) {\bf 37} (1985), no. 4, 405--421.

\bibitem[MR\,85]{MR85} A.~Medina, Ph.~Revoy, {\sl Alg\`ebres de Lie et produit scalaire invariant.}
Ann. Sci. \`Ecole Norm. Sup. (4) {\bf 18} (1985), no. 3, 553--561.

\bibitem[N\,03]{N03} Th.~Neukirchner, {\sl Solvable Pseudo-Riemannian Symmetric
Spaces.}
arXiv:math.DG/0301326, 2003.

\end{thebibliography}
\end{document}